\magnification=1100
\overfullrule0pt

\input amssym.def
\input prepictex
\input pictex
\input postpictex



\def\CC{{\Bbb C}}
\def\FF{{\Bbb F}}

\def\KK{{\Bbb K}}

\def\QQ{{\Bbb Q}}
\def\RR{{\Bbb R}}
\def\ZZ{{\Bbb Z}}

\def\fg{\frak{g}}
\def\fh{\frak{h}}

\def\Hom{\hbox{Hom}}


\font\smallcaps=cmcsc10
\font\titlefont=cmr10 scaled \magstep1

\font\sectionfont=cmbx10
\font\tinyrm=cmr10 at 8pt


\newcount\sectno
\newcount\subsectno
\newcount\resultno

\def\section #1. #2\par{
\sectno=#1
\resultno=0
\bigskip\noindent{\sectionfont #1.  #2}~\medbreak}

\def\subsection #1\par{\bigskip\noindent{\it  #1} \medbreak}


\def\prop{ \global\advance\resultno by 1
\medskip\noindent{\bf Proposition \the\sectno.\the\resultno. }\sl}
\def\lemma{ \global\advance\resultno by 1
\medskip\noindent{\bf Lemma \the\sectno.\the\resultno. }
\sl}
\def\fact{ \global\advance\resultno by 1
\medskip\noindent{\bf Fact \the\sectno.\the\resultno. }
\sl}

\def\remark{ \global\advance\resultno by 1
\medskip\noindent{\bf Remark \the\sectno.\the\resultno. }}
\def\example{ \global\advance\resultno by 1
\medskip\noindent{\bf Example \the\sectno.\the\resultno. }}
\def\cor{ \global\advance\resultno by 1
\medskip\noindent{\bf Corollary \the\sectno.\the\resultno. }\sl}
\def\thm{ \global\advance\resultno by 1
\medskip\noindent{\bf Theorem \the\sectno.\the\resultno. }\sl}
\def\defn{ \global\advance\resultno by 1
\medskip\noindent{\it Definition \the\sectno.\the\resultno. }\slrm}
\def\endthm{\rm\medskip}

\def\endlemma{\rm\medskip}

\def\pf{\rm\smallskip\noindent{\it Proof. }}
\def\endpf{\qed\hfil\medskip}

\newbox\strutAbox
\setbox\strutAbox=\hbox{\vrule height 12pt depth6pt width0pt}

\newbox\strutBbox
\setbox\strutBbox=\hbox{\vrule height 10pt depth5pt width0pt}

\newbox\strutDbox
\setbox\strutDbox=\hbox{\vrule height 11pt depth5pt width0pt}

\newbox\strutHbox
\setbox\strutHbox=\hbox{\vrule height 11pt depth1pt width0pt}

\newbox\strutLbox
\setbox\strutLbox=\hbox{\vrule height 1pt depth5pt width0pt}

\def\ignore#1{\relax}


\def\qed{\hbox{\hskip 1pt\vrule width4pt height 6pt depth1.5pt \hskip 1pt}}

\def\sqr#1#2{{\vcenter{\vbox{\hrule height.#2pt
\hbox{\vrule width.#2pt height#1pt \kern#1pt
\vrule width.2pt}
\hrule height.2pt}}}}
\def\square{\mathchoice\sqr54\sqr54\sqr{3.5}3\sqr{2.5}3}


\def\formula{\global\advance\resultno by 1
\eqno{(\the\sectno.\the\resultno)}}
\def\formulano{\global\advance\resultno by 1 (\the\sectno.\the\resultno)}
\def\tableno{\global\advance\resultno by 1
\the\sectno.\the\resultno. }
\def\lformula{\global\advance\resultno by 1
\leqno(\the\sectno.\the\resultno)}


\def\mapright#1{\smash{\mathop
        {\longrightarrow}\limits^{#1}}}

\def\monthname {\ifcase\month\or January\or February\or March\or April\or
May\or June\or
July\or August\or September\or October\or November\or December\fi}

\newcount\mins  \newcount\hours  \hours=\time \mins=\time
\def\now{\divide\hours by60 \multiply\hours by60 \advance\mins by-\hours
     \divide\hours by60         
     \ifnum\hours>12 \advance\hours by-12
       \number\hours:\ifnum\mins<10 0\fi\number\mins\ P.M.\else
       \number\hours:\ifnum\mins<10 0\fi\number\mins\ A.M.\fi}


\nopagenumbers
\def\runningtitle{\smallcaps kostka-foulkes polynomials}
\headline={\ifnum\pageno>1\eoheadline\else\firstheadline\fi}
\def\names{\smallcaps k.\ nelsen and a.\ ram}
\def\firstheadline{}
\def\eoheadline{\ifodd\pageno\oddheadline\else\evenheadline\fi}
\def\oddheadline{\tenrm\hfil\runningtitle\hfil\folio}
\def\evenheadline{\tenrm\folio\hfil{\names}\hfil}

\vphantom{$ $}  
\vskip.75truein
\centerline{\titlefont Kostka-Foulkes polynomials and Macdonald
spherical functions}
\bigskip
\centerline{\rm Kendra Nelsen}
\centerline{Department of Mathematics}
\centerline{University of Wisconsin, Madison}
\centerline{Madison, WI 53706 USA}
\centerline{{\tt kanelsen@students.wisc.edu}}
\bigskip
\centerline{\rm Arun Ram}
\centerline{Department of Mathematics}
\centerline{University of Wisconsin, Madison}
\centerline{Madison, WI 53706 USA}
\centerline{{\tt ram@math.wisc.edu}}
\medskip

\footnote{}{\tinyrm 
\hskip-.3in
Research partially supported by the National Science Foundation
(DMS-0097977) the National Security Agency (MDA904-01-1-0032)
and by EPSRC Grant GR K99015 at the Newton Institute for
Mathematical Sciences.}
\footnote{}{\tinyrm
}

\bigskip

\noindent{\bf Abstract.}
Generalized Hall-Littlewood polynomials (Macdonald spherical functions)
and generalized Kostka-Foulkes polynomials ($q$-weight multiplicities)
arise in many places in combinatorics, representation theory,
geometry, and mathematical physics.
This paper attempts to organize the different definitions
of these objects and prove the fundamental combinatorial
results from ``scratch'', in a presentation which, hopefully,
will be accessible and useful for both the nonexpert and
researchers currently working in this very active field.
The combinatorics of the affine Hecke algebra plays a central role.
The final section of this paper can be read independently of
the rest of the paper.  It presents, with proof, Lascoux and
Sch\"utzenberger's positive formula for the Kostka-Foulkes poynomials
in the type A case.

\section 0. Introduction

The classical theory of Hall-Littlewood polynomials
and the Kostka-Foulkes polynomials appears in 
the monograph of I.G.\ Macdonald [Mac].
The Hall-Littlewood polynomials form a basis
of the ring of symmetric functions and the Kostka-Foulkes
polynomials are the entries of the transition matrix
between the Hall-Littlewood polynomials and the Schur functions.

This theory enters in many different places in
algebra, geometry and combinatorics.  Many of these 
connections appear in [Mac]:
\smallskip\noindent
\itemitem{(a)} [Mac, Ch.\ II] explains how this theory
describes the structure of the Hall algebra of finite 
${\frak{o}}$-modules, where ${\frak{o}}$ is a discrete
valuation ring.
\smallskip\noindent
\itemitem{(b)} [Mac, Ch.\ IV] explains how the Hall-Littlewood
polynomials enter into the representation theory of $GL_n(\FF_q)$
where $\FF_q$ is a finite field with $q$ elements.
\smallskip\noindent
\itemitem{(c)} [Mac, Ch,\ V] shows that the Hall-Littlewood
polynomials arise as spherical functions for $GL_n(\QQ_p)$
where $\QQ_p$ is the field of $p$-adic numbers.
\smallskip\noindent
\itemitem{(d)} [Mac, Ch.\ III \S 6 Ex.\ 6] explains how the
Kostka-Foulkes polynomials relate to the intersection
cohomology of unipotent orbit closures for $GL_n(\CC)$
and [Mac, Ch.\ III \S 8 Ex.\ 8] explains how the Kostka-Foulkes
polynomials describe the graded decomposition of the 
representations of the symmetric groups $S_n$ on the
cohomology of Springer fibers.
\smallskip\noindent
\itemitem{(e)} [Mac, Ch.\ App.\ A \S 8 and Ch.\ III \S 6] gives
that the Kostka-Foulkes polynomials are $q$-analogues of the 
weight multiplicities for representations of $GL_n(\CC)$.
\smallskip\noindent
\itemitem{(f)} [Mac, Ch.\ III (6.5)] explains how the Kostka-Foulkes
polynomials encode a subtle statistic on column strict Young tableaux.

\smallskip
Macdonald [Mac2, (4.1.2)] showed that there is a formula
for the spherical functions for the Chevalley group $G(\QQ_p)$
which generalizes the formula for Hall-Littlewood symmetric 
functions.  This combinatorial formula is in terms of 
the root system data of the Chevalley group $G$.
In [Lu] Lusztig showed that Macdonald's spherical function formula
can be seen in terms of the affine Hecke algebra and that
the ``$q$-weight multiplicities'' or generalized Kostka-Foulkes
polynomials coming from these spherical functions are Kazhdan-Lusztig
polynomials for the affine Weyl group.  Kato [Kt] proved 
the ``partition function formula'' for the $q$-weight
multiplicities which was conjectured by Lusztig.
The partition function formula has led to continuing analysis
of the connection between the $q$-weight multiplicities, 
functions on nilpotent orbits, filtrations of
weight spaces by the kernels of powers of a regular
nilpotent element, and degrees in harmonic polynomials
(see [JLZ] and the references there).

The connection between Hall-Littlewood polynomials and 
${\frak{o}}$-modules has seen generalizations in the
theory of representations of quivers, the classical
case being the case where the quiver is a loop consisting
of one vertex and one edge.  This theory has been generalized
extensively by Ringel, Lusztig, Nakajima and many others and
is developing quickly; fairly recent references are [Nak1] and [Nak2].

The connection to Springer representations of Weyl groups
and the representations of Chevalley groups over finite fields
has been developed extensively by Lusztig, Shoji and others;
a good survey of the current theory is in [Shj1] and the
recent papers [Shj2] show how this theory is beginning to extend
its reach outside Lie theory into the realm of complex
reflection groups.

Since the theory of Macdonald spherical functions
(the generalization of Hall-Littlewood polynomials)
and $q$-weight multiplicities (the generalization of
Kostka-Foulkes polynomials) appears in so many important
parts of mathematics it seems appropriate to give a 
survey of the basics of this theory.
This paper is an attempt to collect together the fundamental
combinatorial results analogous to those which are
found for the type A case in [Mac].
The presentation here centers on the
role played by the affine Hecke algebra.  Hopefully this
will help to illustrate how and why these objects
arise naturally from a combinatorial point of view
and, at the same time, provide enough underpinning to the
algebra of the underlying algebraic groups to be useful
to researchers in representation theory.

Using the terms {\it Hall-Littlewood polynomial} and
{\it Macdonald spherical function} interchangeably, and using
the words {\it Kostka-Foulkes polynomial} and 
{\it $q$-weight multiplicity} interchangeably,
the results that we prove in this paper are:
\smallskip\noindent
\itemitem{(1)}  The interpretation of the Hall-Littlewood
polynomials as elements of the affine Hecke algebra (via the 
Satake isomorphism),
\smallskip\noindent
\itemitem{(2)}  Macdonald's spherical function formula,
\smallskip\noindent
\itemitem{(3)}  The expansion of the Hall Littlewood polynomial
in terms of the standard basis of the affine Hecke algebra,
\smallskip\noindent
\itemitem{(4)}  The triangularity of transition matrices
between Macdonald spherical functions and other bases of
symmetric functions,
\smallskip\noindent
\itemitem{(5)}  The straightening rules for Hall-Littlewood polynomials,
\smallskip\noindent
\itemitem{(6)}  The orthogonality of Macdonald spherical functions,
\smallskip\noindent
\itemitem{(7)}  The raising operator formula for Kostka-Foulkes polynomials,
\smallskip\noindent
\itemitem{(8)}  The partition function formula for $q$-weight multiplicities,
\smallskip\noindent
\itemitem{(9)}  The identification of the Kostka-Foulkes polynomial
as a Kazhdan-Lusztig polynomial.
\smallskip\noindent
All of these results are proved here in general Lie type.  They
are all previously known, spread throughout various
parts of the literature.  The presentation here is a unified
one; some of the proofs may (or may not) be new.

Section 4 is designed so that it can be read independently
of the rest of the paper.  In Section 4
we give the proof of Lascoux-Sch\"utzenberger's positive combinatorial 
formula [LS] (see also [Mac, Ch.\ III (6.5)]) for Kostka-Foulkes 
polynomials in type A.  Versions of this proof have appeared previously 
in [Sch] and in [Bt].  This proof has a reputation for being difficult 
and obscure.  After finally getting the courage to attack the literature, 
we have found, in the end, that the proof is not so difficult after all.  
Hopefully we have been able to explain it so that others will also find it so.

\medskip\noindent
{\bf Acknowledgements.}  
A portion of this paper was written during a stay
of A.\ Ram at the Newton Institute for the Mathematical 
Sciences at Cambridge University.  A.\ Ram thanks them for 
their hospitality and support during Spring 2001.
The preparation of this paper has been greatly
aided by handwritten lecture notes of I.G.\ Macdonald
from lectures he gave at the University of California, San Diego, in 
Spring 1991.  In several places we have copied rather
unabashedly from them.  Over many years Professor Macdonald
has generously given us lots of handwritten notes. We cannot
thank him enough, these notes have opened our eyes to many 
beautiful things and shown us the ``right way'' many times when
we were going astray.
\bigskip

\section 1. Weyl groups, affine Weyl groups, and the affine Hecke algebra

This section sets up the definitions and notations.  Good references
for this preliminary material are [Bou], [St] and [Mac4].

\subsection The root system and the Weyl group

Let $\fh_\RR^*$ be a real vector space 
with a nondegenerate symmetric bilinear form $\langle\,,\,\rangle$.
The basic data is a reduced irreducible root system $R$ (defined below)
in $\fh_\RR^*$.  Associated to $R$ are the {\it weight lattice} 
$$P = \{ \lambda\in \fh_\RR^*\ |\ 
\langle\lambda,\alpha^\vee\rangle
\in \ZZ \hbox{\ for all $\alpha\in R$}\}
\qquad\hbox{where}\qquad
\alpha^\vee = {2\alpha\over \langle \alpha,\alpha\rangle},
\formula
$$
and the {\it Weyl group}
$$W = \langle s_\alpha\ |\ \alpha\in R\rangle
\qquad\hbox{generated by the reflections}\qquad
\matrix{
s_\alpha\colon &\fh_\RR^* &\longrightarrow &\fh_\RR^* \cr
&\lambda &\longmapsto &\lambda-\langle \lambda,\alpha^\vee\rangle\alpha
\cr
}
\formula$$
in the hyperplanes
$$H_\alpha=\{ x\in \fh_\RR^*\ |\ \langle x,\alpha^\vee\rangle=0\},
\qquad \alpha\in R.
\formula$$
With these definitions $R$ is a reduced irreducible 
root system if it is a subset of $\fh_\RR^*$
such that 
\smallskip\noindent
\itemitem{(a)} $R$ is finite, $0\not\in R$ and 
$\fh_\RR^* = \hbox{$\RR$-span}(R)$,
\smallskip\noindent
\itemitem{(b)} $W$ permutes the elements of $R$, 
i.e. $w\alpha\in R$ for $w\in W$ and $\alpha\in R$, 
\smallskip\noindent
\itemitem{(c)} $W$ is finite, 
\smallskip\noindent
\itemitem{(d)} $R\subseteq P$,
\smallskip\noindent
\itemitem{(e)} if $\alpha\in R$ then the only other multiple of 
$\alpha$ in $R$ is $-\alpha$,
\smallskip\noindent
\itemitem{(f)} $\fh_\RR^*$ is an irreducible $W$-module.
\smallskip\noindent
The choice of a fundamental region $C$ for the
action of $W$ on $\fh_\RR^*$ is equivalent to a choice
of {\it positive roots} $R^+$ of $R$,
$$R^+=\{ \alpha\in R\ |\ \langle x,\alpha^\vee\rangle>0
\hbox{\ for all $x\in C$}\}
\qquad\hbox{and}\quad
C = \{ x\in \fh_\RR^*\ |\ \langle x,\alpha^\vee\rangle>0
\hbox{\ for all $\alpha\in R^+$}\}.$$

\example
If $\fh_\RR^* = \RR^2$ with orthonormal basis
$\varepsilon_1=(1,0)$ and $\varepsilon_2= (0,1)$,
$P = \hbox{$\ZZ$-span}\{\varepsilon_1,\varepsilon_2\}$,
and $W = \{ 1, s_1, s_2, s_1s_2, s_2s_1, s_1s_2s_1, s_2s_1s_2, 
s_1s_2s_1s_2\}$ is the group of order $8$ generated by the reflections
$s_1$ and $s_2$ in the hyperplanes $H_{\alpha_1}$ and $H_{\alpha_2}$,
respectively, where
$$
\matrix{
\alpha_1 = 2\varepsilon_1,\hfill &\qquad 
&\alpha_1^\vee = \varepsilon_1,\hfill \cr
\alpha_2 = \varepsilon_2-\varepsilon_1, \hfill
&&\alpha_2^\vee = \alpha_2,\hfill \cr}
\qquad\hbox{then}\qquad
R = \{\pm \alpha_1, \pm \alpha_2, \pm(\alpha_1+\alpha_2),
\pm(\alpha_1+2\alpha_2)\}.
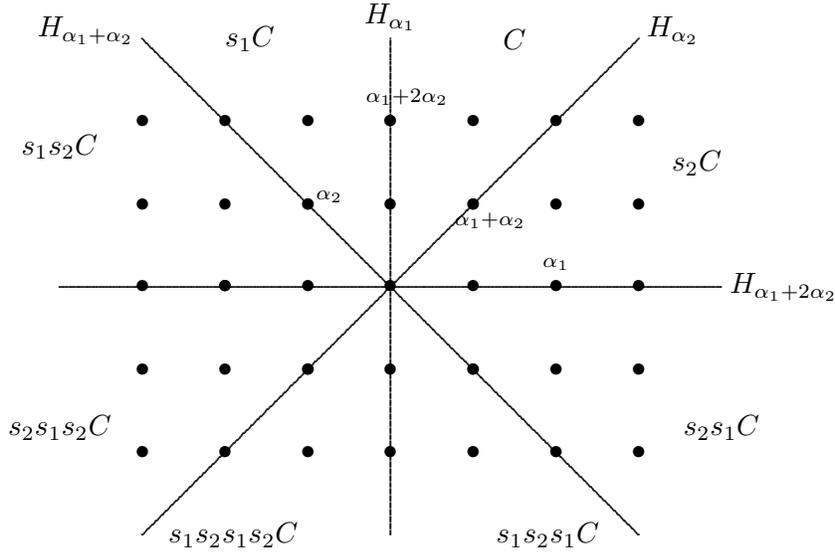
$$
$$
\beginpicture
\setcoordinatesystem units <1cm,1cm>         
\setplotarea x from -4 to 4, y from -4 to 4    
\put{$H_{\alpha_1}$}[b] at 0 3.1
\put{$H_{\alpha_2}$}[l] at 3.1 3.1
\put{$H_{\alpha_1+\alpha_2}$}[r] at -3.1 3.1
\put{$H_{\alpha_1+2\alpha_2}$}[l] at 4.1 0
\put{$C$} at 1.5 3
\put{$s_1C$} at -1.7 3
\put{$s_2C$} at  3.7 1.5
\put{$s_1s_2C$} at  -4 1.7
\put{$s_2s_1C$} at   4  -1.7 
\put{$s_1s_2s_1C$} at  1.9 -3
\put{$s_2s_1s_2C$} at  -4 -1.7
\put{$s_1s_2s_1s_2C$} at  -1.9 -3
\put{$\bullet$} at -3 -2
\put{$\bullet$} at -2 -2
\put{$\bullet$} at -1 -2
\put{$\bullet$} at  0 -2
\put{$\bullet$} at  1 -2
\put{$\bullet$} at  2 -2
\put{$\bullet$} at  3 -2
\put{$\bullet$} at -3 -1
\put{$\bullet$} at -2 -1
\put{$\bullet$} at -1 -1
\put{$\bullet$} at  0 -1
\put{$\bullet$} at  1 -1
\put{$\bullet$} at  2 -1
\put{$\bullet$} at  3 -1
\put{$\bullet$} at -3 0
\put{$\bullet$} at -2 0
\put{$\bullet$} at -1 0
\put{$\bullet$} at  0 0
\put{$\bullet$} at  1 0
\put{$\bullet$} at  2 0
\put{$\bullet$} at  3 0
\put{$\bullet$} at -3 1
\put{$\bullet$} at -2 1
\put{$\bullet$} at -1 1
\put{$\bullet$} at  0 1
\put{$\bullet$} at  1 1
\put{$\bullet$} at  2 1
\put{$\bullet$} at  3 1
\put{$\bullet$} at -3 2
\put{$\bullet$} at -2 2
\put{$\bullet$} at -1 2
\put{$\bullet$} at  0 2
\put{$\bullet$} at  1 2
\put{$\bullet$} at  2 2
\put{$\bullet$} at  3 2
\put{$\bullet$} at 0 2
\put{$\bullet$} at -2 0
\put{$\bullet$} at 1 1
\put{$\bullet$} at -1 -1
\put{$\bullet$} at 1 -1
\put{$\bullet$} at -1 1
\put{$\scriptstyle{\alpha_2}$}[l] at -0.9 1.1
\put{$\scriptstyle{\alpha_1}$}[b] at 2.0  0.2  
\put{$\scriptstyle{\alpha_1+\alpha_2}$}[t] at 1.2  0.9  
\put{$\scriptstyle{\alpha_1+2\alpha_2}$}[bl] at -0.3  2.2  
\plot -3 -3   3 3 /
\plot  3 -3  -3 3 /
\plot  0  3   0 -3 /
\plot  4  0  -4  0 /
\endpicture
$$
This is the root system of {\it type $C_2$}.
\endpf

For each $\alpha\in R^+$ define the {\it raising operator} 
$R_\alpha\colon P\to P$ by $R_\alpha\mu=\mu+\alpha.$
The {\it dominance order} on $P$ is given by
$$\mu\le \lambda
\qquad\hbox{if}\qquad
\lambda = R_{\beta_1}\cdots R_{\beta_\ell}\mu\formula$$
for some sequence of positive roots $\beta_1,\ldots, \beta_\ell\in R^+$.

The various fundamental chambers for the action of $W$ on $\fh_{\RR}^*$
are the $w^{-1}C$, $w\in W$.  The {\it inversion set} of an element 
$w\in W$ is
$$
R(w) = \{ \alpha\in R^+\ |\ 
\hbox{$H_\alpha$ is between $C$ and $w^{-1}C$}\}
\qquad\hbox{and}\qquad
\ell(w) = {\rm Card}(R(w))
\formula
$$
is the {\it length} of $w$.  
If $R^- = -R^+ = \{ -\alpha\ |\ \alpha\in R^+\}$ then
$$R = R^+\cup R^-
\qquad\hbox{and}\qquad 
R(w) = \{ \alpha\in R^+\ |\ w\alpha\in R^-\},
\qquad\hbox{for $w\in W$.}$$

The weight lattice, the set of {\it dominant integral weights},
and the set of {\it strictly dominant integral weights}, are
$$
\eqalign{
P~~ &= 
\{ \lambda\in \fh_\RR^*\ |\ 
\hbox{$\langle \lambda,\alpha^\vee\rangle\in \ZZ$ for all
$\alpha\in R$} \}, \cr
P^+ = P\cap \overline C &=
\{ \lambda\in \fh_\RR^*\ |\ 
\hbox{$\langle \lambda,\alpha^\vee\rangle\in \ZZ_{\ge 0}$ for all
$\alpha\in R^+$} \}, \cr
P^{++} = P\cap C &= 
\{ \lambda\in \fh_\RR^*\ |\ 
\hbox{$\langle \lambda,\alpha^\vee\rangle\in \ZZ_{>0}$ for all
$\alpha\in R^+$} \}, \cr
}\formula$$
where 
$\overline C = \{ x\in \fh_\RR^*\ |\ 
\hbox{$\langle x,\alpha^\vee\rangle \ge 0$
for all $\alpha\in R^+$}\}$
is the closure of the fundamental chamber $C$.

The {\it simple roots} are the positive roots
$\alpha_1,\ldots,\alpha_n$ such that the hyperplanes
$H_{\alpha_i}$, $1\le i\le n$, are the {\it walls} of $C$.
The {\it fundamental weights}, $\omega_1,\ldots, \omega_n\in P$,
are given by $\langle \omega_i,\alpha_j^\vee\rangle = \delta_{ij}$,
$1\le i,j\le n$, and 
$$P= \sum_{i=1}^n \ZZ \omega_i,
\qquad
P^+= \sum_{i=1}^n \ZZ_{\ge 0} \omega_i, \qquad\hbox{and}\qquad
P^{++} = \sum_{i=1}^n \ZZ_{>0} \omega_i.
\formula$$
The set $P^+$ is an integral cone with vertex $0$,
the set $P^{++}$ is a integral cone with vertex
$$\rho = \sum_{i=1}^n \omega_i 
= \hbox{$1\over 2$}\sum_{\alpha\in R^+} \alpha,
\qquad\hbox{and the map}\qquad
\matrix{P^+ &\longrightarrow &P^{++} \cr
\lambda &\longmapsto &\lambda+\rho \cr
}\formula$$
is a bijection.  In Example 1.4, with the root system of
type $C_2$, the picture is
$$\matrix
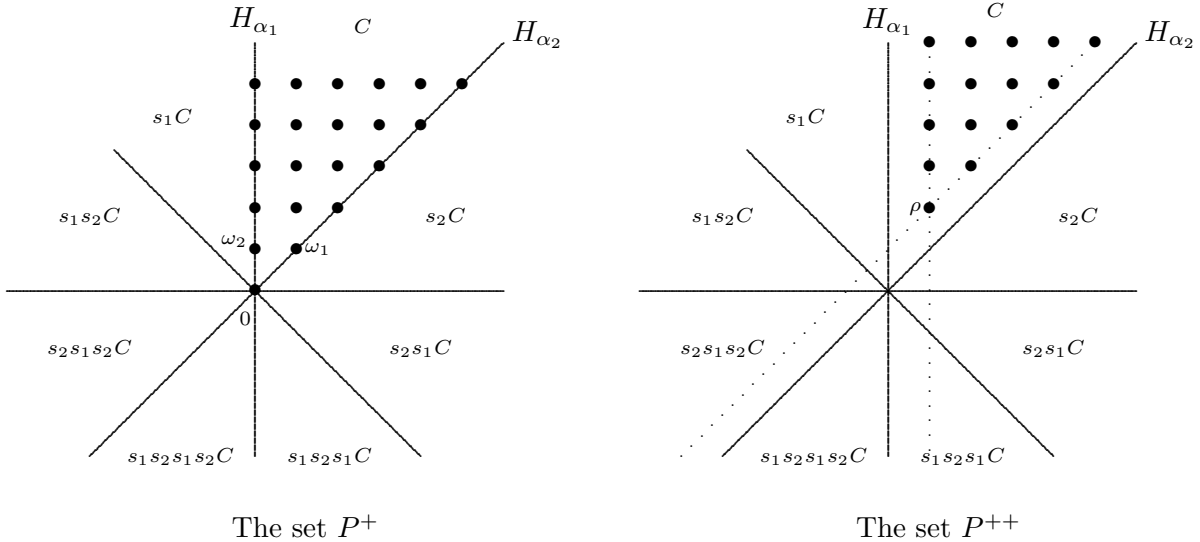
{
\beginpicture
\setcoordinatesystem units <1cm,1cm>         
\setplotarea x from -2 to 3.5, y from -2 to 3.5    
\put{$H_{\alpha_1}$}[b] at 0 3.1
\put{$H_{\alpha_2}$}[l] at 3.1 3.1
\put{$\scriptstyle{0}$}[tr] at -0.05 -0.25
\put{$\scriptstyle{\omega_1}$}[l] at 0.6 0.5
\put{$\scriptstyle{\omega_2}$}[r] at -0.1 0.6
\put{$\scriptstyle{C}$} at 1.3 3.2
\put{$\scriptstyle{s_1C}$} at -1.0 2.1
\put{$\scriptstyle{s_2C}$} at  2.3 0.9
\put{$\scriptstyle{s_1s_2C}$} at  -2 0.9
\put{$\scriptstyle{s_2s_1C}$} at   2  -0.7 
\put{$\scriptstyle{s_1s_2s_1C}$} at  0.9 -2
\put{$\scriptstyle{s_2s_1s_2C}$} at  -2.0 -0.7
\put{$\scriptstyle{s_1s_2s_1s_2C}$} at  -0.9 -2
\put{$\bullet$} at  0 0
\put{$\bullet$} at  0 0.5
\put{$\bullet$} at  0.5 0.5
\put{$\bullet$} at  0 1
\put{$\bullet$} at  0.5 1
\put{$\bullet$} at  1 1
\put{$\bullet$} at  0 1.5
\put{$\bullet$} at  0.5 1.5
\put{$\bullet$} at  1 1.5
\put{$\bullet$} at  1.5 1.5
\put{$\bullet$} at  0 2
\put{$\bullet$} at  0.5 2
\put{$\bullet$} at  1 2
\put{$\bullet$} at  1.5 2
\put{$\bullet$} at  2 2
\put{$\bullet$} at  0 2.5
\put{$\bullet$} at  0.5 2.5
\put{$\bullet$} at  1 2.5
\put{$\bullet$} at  1.5 2.5
\put{$\bullet$} at  2 2.5
\put{$\bullet$} at  2.5 2.5
\plot -2 -2   3 3 /
\plot  2 -2  -1.7 1.7 /
\plot  0  3   0 -2 /
\plot  3  0  -3  0 /
\endpicture
&\qquad
&
\beginpicture
\setcoordinatesystem units <1cm,1cm>         
\setplotarea x from -2 to 3.5, y from -2 to 3.5    
\put{$H_{\alpha_1}$}[b] at 0 3.1
\put{$H_{\alpha_2}$}[l] at 3.1 3.1
\put{$\scriptstyle{\rho}$}[r] at 0.4 1.0
\put{$\scriptstyle{C}$} at 1.3 3.4
\put{$\scriptstyle{s_1C}$} at -1.0 2.1
\put{$\scriptstyle{s_2C}$} at  2.3 0.9
\put{$\scriptstyle{s_1s_2C}$} at  -2 0.9
\put{$\scriptstyle{s_2s_1C}$} at   2  -0.7 
\put{$\scriptstyle{s_1s_2s_1C}$} at  0.9 -2
\put{$\scriptstyle{s_2s_1s_2C}$} at  -2.0 -0.7
\put{$\scriptstyle{s_1s_2s_1s_2C}$} at  -0.9 -2
\put{$\bullet$} at  0.5 1
\put{$\bullet$} at  0.5 1.5
\put{$\bullet$} at  1 1.5
\put{$\bullet$} at  0.5 2
\put{$\bullet$} at  1 2
\put{$\bullet$} at  1.5 2
\put{$\bullet$} at  0.5 2.5
\put{$\bullet$} at  1 2.5
\put{$\bullet$} at  1.5 2.5
\put{$\bullet$} at  2 2.5
\put{$\bullet$} at  0.5 3
\put{$\bullet$} at  1 3
\put{$\bullet$} at  1.5 3
\put{$\bullet$} at  2 3
\put{$\bullet$} at  2.5 3
\plot -2 -2   3 3 /
\plot  2 -2  -1.7 1.7 /
\plot  0  3   0 -2 /
\plot  3  0  -3  0 /
\setdots
\plot  0.5  3   0.5 -2 /
\plot -2.5 -2   2.5 3 /
\endpicture
\cr
\cr
\cr
\hbox{The set $P^+$}
&&\hbox{The set $P^{++}$} \cr
}
$$
The {\it simple reflections} are $s_i = s_{\alpha_i}$, for $1\le i\le n$.
The Weyl group $W$ has a presentation by generators
$s_1, \ldots, s_n$ and relations
$$\matrix{
s_i^2=1, &\qquad &\hbox{for $1\le i\le n$}, \cr
\underbrace{s_is_js_j\cdots}_{m_{ij}\ {\rm factors}}
=\underbrace{s_js_is_j\cdots}_{m_{ij}\ {\rm factors}},
&&i\ne j, \cr
}
\formula
$$
where $\pi/m_{ij}$ is the angle between the hyperplanes
$H_{\alpha_i}$ and $H_{\alpha_j}$.
A {\it reduced word} for $w\in W$ is an expression
$w=s_{i_1}\cdots s_{i_p}$ for $w$ as a product of simple
reflections which has $p$ minimal.  
The following lemma describes the inversion
set in terms of the simple roots and the simple reflections
and shows that if
$w=s_{i_1}\cdots s_{i_p}$ is a reduced expression 
for $w$ then $p=\ell(w)$.

\lemma [Bou VI \S 1 no.\ 6 Cor.\ 2 to Prop.\ 17]
{\sl Let $w=s_{i_1}\cdots s_{i_p}$ be a reduced word for $w$.
Then
$$R(w) = \{ \alpha_{i_p}, s_{i_p}\alpha_{i_{p-1}},
\ldots, s_{i_p}\cdots s_{i_2}\alpha_{i_1}\}.$$
}
\endthm

\noindent
The {\it Bruhat order}, or {\it Bruhat-Chevalley order},  
(see [St, \S 8 App., p.\ 126]) is the partial order on $W$ such that
$v\le w$
if there is a reduced word for $v$, $v = s_{j_1}\cdots s_{j_k}$,
which is a subword of a reduced word for $w$, $w=s_{i_1}\cdots s_{i_p}$,
(i.e. $s_{j_1}, \ldots, s_{j_k}$ is a subsequence of the sequence
$s_{i_1},\ldots, s_{i_p}$). 

\bigskip
\subsection The affine Weyl group

\medskip
For $\lambda\in P$, the {\it translation in $\lambda$} is
$$\matrix{ t_\lambda\colon &\fh_\RR^* &\longrightarrow &\fh_\RR^* \cr
&x &\longmapsto &x+\lambda. \cr
}
\formula$$
The {\it extended affine Weyl group} $\tilde W$ is the group
$$\tilde W = \{ wt_\lambda\ |\ w\in W, \lambda\in P\},
\formula$$
with multiplication determined by the relations
$$
t_\lambda t_\mu = t_{\lambda+\mu},
\qquad\hbox{and}\qquad
wt_\lambda = t_{w\lambda} w, 
\formula$$
for $\lambda,\mu\in P$ and $w\in W$.
The group $\tilde W$ is the group
of transformations of $\fh_\RR^*$ generated by the $s_\alpha$,
$\alpha\in R^+$, and $t_\lambda$, $\lambda\in P$.
The {\it affine Weyl group} $W_{\rm aff}$ is the subgroup
of $\tilde W$ generated by the reflections
$$s_{\alpha,k}\colon \fh_\RR^* \to \fh_\RR^*
\qquad\hbox{in hyperplanes}\qquad
H_{\alpha,k} = \{ x\in \fh_\RR^*\ |\ 
\langle x, \alpha^\vee\rangle = k\},
\qquad\alpha\in R^+, k\in \ZZ.
\formula$$
The reflections $s_{\alpha,k}$ can be written as elements of 
$\tilde W$ via the formula
$$s_{\alpha,k} = t_{k\alpha^\vee}s_\alpha = s_\alpha t_{-k\alpha^\vee}.
\formula$$

The {\it highest root} of $R$ is the unique element $\varphi\in R^+$
such that the {\it fundamental alcove}
$$A = C\cap \{ x\in \fh_\RR^*\ |\ \langle x,\varphi^\vee\rangle<1\}
\formula$$
is a fundamental region for the action of $W_{\rm aff}$ on $\fh_\RR^*$.
The various fundamental chambers for the action of $W_{\rm aff}$ on 
$\fh_\RR^*$ are $w^{-1}A$, $w\in W_{\rm aff}$.
The {\it inversion set} of $w\in \tilde W$ is 
$$R(w)= \{ H_{\alpha,k}\ |\ 
\hbox{$H_{\alpha,k}$ is between $A$ and $w^{-1}A$}\}
\qquad\hbox{and}\qquad
\ell(w) = {\rm Card}(R(w))$$
is the {\it length} of $w$.  If $w\in W$ and $\lambda\in P$ then
$$\ell(wt_\lambda) = \sum_{\alpha\in R^+}
|\langle \lambda,\alpha^\vee\rangle + \chi(w\alpha)|,
\formula$$
where, for a root $\beta\in R$, set  
$\chi(\beta) = 0$, if $\beta\in R^+$, and
$\chi(\beta)= 1$, if $\beta\in R^-$. 

Continuing Example 1.4, we have the picture
$$
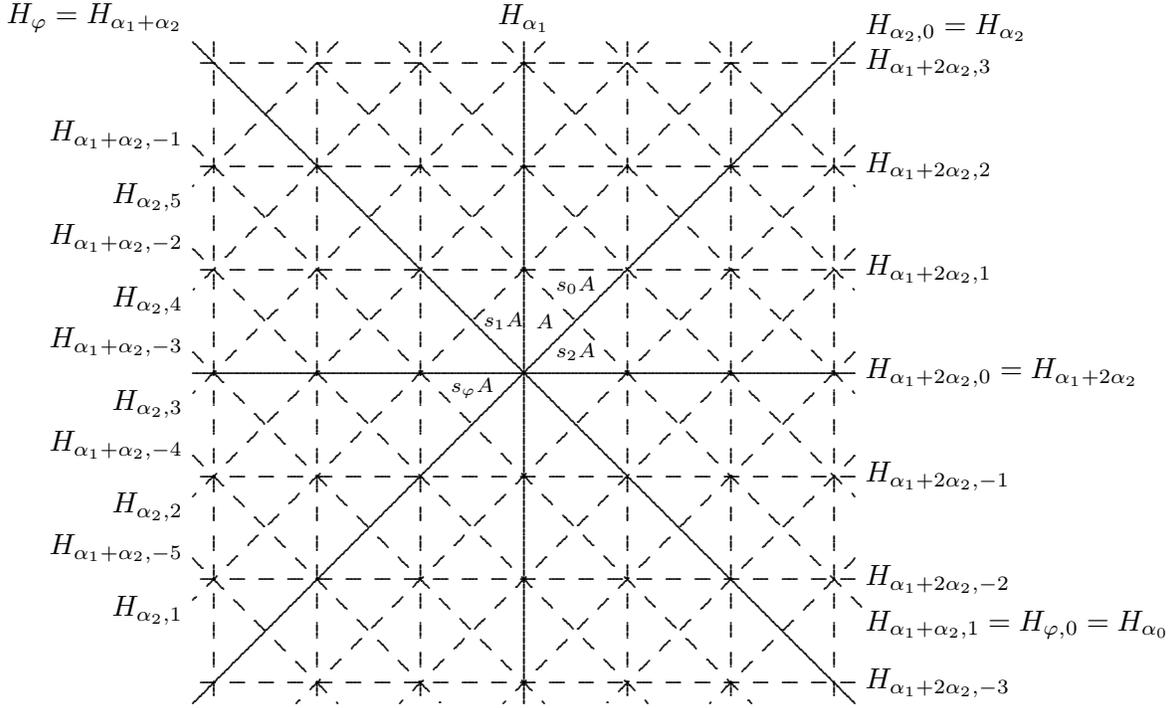
\beginpicture
\setcoordinatesystem units <1.25cm,1.25cm>         
\setplotarea x from -4.5 to 4.5, y from -3.5 to 4    
\put{$H_{\alpha_1}$}[b] at 0 3.3
\put{$H_{\alpha_2,0}=H_{\alpha_2}$}[bl] at 3.3 3.2
\put{$H_{\alpha_2,1}$}[r] at -3.3 -2.3
\put{$H_{\alpha_2,2}$}[r] at -3.3 -1.3
\put{$H_{\alpha_2,3}$}[r] at -3.3 -0.3
\put{$H_{\alpha_2,4}$}[r] at -3.3 0.7
\put{$H_{\alpha_2,5}$}[r] at -3.3 1.7
\put{$H_\varphi = H_{\alpha_1+\alpha_2}$}[br] at -3.3 3.3
\put{$H_{\alpha_1+\alpha_2,1}=H_{\varphi,0}=H_{\alpha_0}$}[tl] at 3.3 -2.3 
\put{$H_{\alpha_1+\alpha_2,-1}$}[r] at -3.3 2.3
\put{$H_{\alpha_1+\alpha_2,-2}$}[r] at -3.3 1.3
\put{$H_{\alpha_1+\alpha_2,-3}$}[r] at -3.3 0.3
\put{$H_{\alpha_1+\alpha_2,-4}$}[r] at -3.3 -0.7
\put{$H_{\alpha_1+\alpha_2,-5}$}[r] at -3.3 -1.7
\put{$H_{\alpha_1+2\alpha_2,0}=H_{\alpha_1+2\alpha_2}$}[l] at 3.3 0
\put{$H_{\alpha_1+2\alpha_2,1}$}[l] at 3.3 1
\put{$H_{\alpha_1+2\alpha_2,2}$}[l] at 3.3 2
\put{$H_{\alpha_1+2\alpha_2,3}$}[l] at 3.3 3
\put{$H_{\alpha_1+2\alpha_2,-1}$}[l] at 3.3 -1
\put{$H_{\alpha_1+2\alpha_2,-2}$}[l] at 3.3 -2
\put{$H_{\alpha_1+2\alpha_2,-3}$}[l] at 3.3 -3
\put{$\scriptstyle{A}$} at 0.2 0.5
\put{$\scriptstyle{s_1A}$} at -0.2 0.5
\put{$\scriptstyle{s_2A}$} at 0.5 0.2
\put{$\scriptstyle{s_0A}$} at 0.5 0.85
\put{$\scriptstyle{s_\varphi A}$} at -0.5 -0.15
\plot -3.2 -3.2   3.2 3.2 /
\plot  3.2 -3.2  -3.2 3.2 /
\plot  0  3.2   0 -3.2 /
\plot  3.2  0  -3.2  0 /
\setdashes
\plot  -3 3.2   -3 -3.2 /
\plot  -2 3.2   -2 -3.2 /
\plot  -1 3.2   -1 -3.2 /
\plot   1 3.2    1 -3.2 /
\plot   2 3.2    2 -3.2 /
\plot   3 3.2    3 -3.2 /
\plot  3.2 -3   -3.2 -3 /
\plot  3.2 -2   -3.2 -2 /
\plot  3.2 -1   -3.2 -1 /
\plot  3.2  1   -3.2  1 /
\plot  3.2  2   -3.2  2 /
\plot  3.2  3   -3.2  3 /
\plot  3.2 -1.8   1.8 -3.2 /
\plot  3.2 -0.8   0.8 -3.2 /
\plot  3.2 0.2    -0.2 -3.2 /
\plot  3.2 1.2   -1.2 -3.2 /
\plot  3.2 2.2   -2.2 -3.2 /
\plot  2.2 3.2   -3.2 -2.2 /
\plot  1.2 3.2   -3.2 -1.2 /
\plot  0.2 3.2   -3.2 -0.2 /
\plot  -0.8 3.2  -3.2 0.8 /
\plot  -1.8 3.2  -3.2 1.8 /
\plot  -2.2 3.2   3.3 -2.3 /
\plot  -1.2 3.2   3.2 -1.2 /
\plot  -0.2 3.2   3.2 -0.2 /
\plot  0.8 3.2  3.2 0.8 /
\plot  1.8 3.2  3.2 1.8 /
\plot  -3.2 -1.8   -1.8 -3.2 /
\plot  -3.2 -0.8   -0.8 -3.2 /
\plot  -3.2 0.2    0.2 -3.2 /
\plot  -3.2 1.2   1.2 -3.2 /
\plot  -3.2 2.2   2.2 -3.2 /
\endpicture
$$

Let
$$H_{\alpha_0} = H_{\varphi,1}
\qquad\hbox{and}\qquad s_0= s_{\varphi,1} = t_{\phi^\vee}s_\phi
=s_\phi t_{-\phi^\vee},
\formula$$
and let $H_{\alpha_1},\ldots, H_{\alpha_n}$ and 
$s_1,\ldots, s_n$ be as in (1.10).
Then the walls of $A$ are the hyperplanes $H_{\alpha_0},
H_{\alpha_1}, \ldots, H_{\alpha_n}$ and the 
group $W_{\rm aff}$ has a presentation by generators
$s_0, s_1,\ldots, s_n$ and relations
$$\matrix{
s_i^2=1, &\qquad &\hbox{for $0\le i\le n$}, \cr
\underbrace{s_is_js_j\cdots}_{m_{ij}\ {\rm factors}}
=\underbrace{s_js_is_j\cdots}_{m_{ij}\ {\rm factors}},
&&i\ne j, \cr
}
\formula$$
where $\pi/m_{ij}$ is the angle between the hyperplanes
$H_{\alpha_i}$ and $H_{\alpha_j}$.

Let $w_0$ be the longest element of $W$ and let
$w_i$ be the longest element of the subgroup
$W_{\omega_i} = \{ w\in W\ |\ w\omega_i=\omega_i\}$.
Let $\varphi^\vee = c_1\alpha^\vee_1+\cdots c_n\alpha^\vee_n$.  Then
(see [Bou, VI \S 2 no.\ 3 Prop.\ 6])
$$\Omega = \{ g\in \tilde W\ |\ \ell(g)=0\}
= \{ g_i\ |\ c_i=1\},
\qquad\hbox{where}\qquad
g_i = t_{\omega_i}w_iw_0.
\formula$$
Each element $g\in \Omega$ sends the alcove $A$ to itself
and thus permutes the walls $H_{\alpha_0}, H_{\alpha_1},
\ldots, H_{\alpha_n}$ of $A$.  Denote the resulting
permutation of $\{0,1,\ldots, n\}$ also by $g$. Then
$$gs_ig^{-1} = s_{g(i)},\qquad\hbox{for $0\le i\le n$},
\formula$$
and the group $\tilde W$ is presented by the generators
$s_0,s_1,\ldots, s_n$ and $g\in \Omega$ with the relations
(1.18) and (1.20).

\vfill\eject

\subsection The affine Hecke algebra

\medskip
Let $\KK=\ZZ[q,q^{-1}]$.
The affine Hecke algebra $\tilde H$ is the algebra over
$\KK$ given by generators
$T_i$, $1\le i\le n$, and $x^\lambda$, $\lambda\in P$,
and relations
$$\matrix{
\underbrace{T_iT_jT_i\cdots}_{m_{ij} {\rm\ factors}}
=
\underbrace{T_jT_iT_j\cdots}_{m_{ij} {\rm\ factors}}\,,
&\qquad &\hbox{for all $i\ne j$, } \hfill\cr
\cr
T_i^2 =(q-q^{-1})T_i+1,
&&\hbox{for all $1\le i\le n$, }\hfill \cr
\cr
x^\lambda x^\mu = x^\mu x^\lambda = x^{\lambda+\mu}\,,
&&\hbox{for all $\lambda,\mu\in P$,} \hfill \cr
\cr
x^\lambda T_i = 
T_i x^{s_i\lambda} + 
(q-q^{-1})\displaystyle{x^\lambda-x^{s_i\lambda}\over 1-x^{-\alpha_i} }
\,,
&&\hbox{for all $1\le i\le n$, $\lambda\in P$.}\hfill \cr
}
\formula$$
An alternative presentation of $\tilde H$ is by the
generators $T_w$, $w\in \tilde W$,  and relations
$$T_{w_1}T_{w_2} = T_{w_1w_2},
\qquad\hbox{if $\ell(w_1w_2)=\ell(w_1)+\ell(w_2)$,}$$
$$T_{s_i}T_w = (q-q^{-1})T_w + T_{s_iw},
\qquad\hbox{if $\ell(s_iw)<\ell(w)$ \quad $(0\le i\le n)$.}$$
With notations as in (1.12-1.20) 
the conversion between the two presentations is given by
the relations
$$\matrix{
\hfill T_w &= T_{i_1}\cdots T_{i_p}, \hfill
&&\hbox{if $w\in W_{\rm aff}$ and 
$w = s_{i_1}\cdots s_{i_p}$ is a reduced word,} \hfill \cr
\cr
T_{g_i} &= x^{\omega_i}T_{w_0w_i}^{-1}, \hfill
&&\hbox{for $g_i\in \Omega$ as in (1.19),} \hfill\cr
\cr
\hfill x^\lambda &= T_{t_\mu}T_{t_\nu}^{-1}, \hfill
&\qquad &\hbox{if $\lambda=\mu-\nu$ with $\mu,\nu\in P^+$}, \hfill \cr
\cr
\hfill T_{s_0} &= T_{s_\phi}x^{-\phi^\vee}, \hfill
&&\hbox{where $\phi$ is the highest root of $R$,} \hfill
\cr}
\formula$$
\medskip\noindent
{\it The Kazhdan-Lusztig basis}

\bigskip
The algebra $\tilde H$ has bases
$$\{ x^\lambda T_w\ |\ w\in W, \lambda\in P\}
\qquad\hbox{and}\qquad
\{ T_w x^\lambda \ |\ w\in W, \lambda\in P\}.
$$
The Kazhdan-Lusztig basis $\{C_w'\ |\ w\in \tilde W\}$ is
another basis of $\tilde H$ which plays an important role.
It is defined as follows.

The {\it bar involution} on $\tilde H$ is the 
$\ZZ$-linear automorphism 
$\overline{\phantom{P}}\colon \tilde H\to\tilde H$
given by
$$\overline{q} = q^{-1}
\qquad\hbox{and}\qquad
\overline{T_w} = T_{w^{-1}}^{-1},
\qquad\hbox{for $w\in \tilde W$.}$$
For $0\le i\le n$, $\overline{T_i} = T_i^{-1} = T_i-(q-q^{-1})$
and the bar involution is a $\ZZ$-algebra automorphism of 
$\tilde H$.   If $w=s_{i_1}\cdots s_{i_p}$ is a reduced word for $w$ then,
by the definition of the Bruhat order (defined after Lemma 1.11),
$$
\eqalign{
\overline{T_w}
&=\overline{T_{i_1}\cdots T_{i_p}}
=\overline{T_{i_1}}\cdots \overline{T_{i_p}}
=T^{-1}_{i_1}\cdots T^{-1}_{i_p} \cr
&=(T_{i_1}-(q-q^{-1}))\cdots (T_{i_p}-(q-q^{-1}))
=T_w + \sum_{v<w} a_{vw}T_v, \cr
}$$
with $a_{vw}\in \ZZ[(q-q^{-1})]$.

Setting $\tau_i = qT_i$ and $t = q^2$, the second relation in (1.21)
$$T_i^2 = (q-q^{-1})T_i +1
\qquad\hbox{becomes}\qquad
\tau_i^2 = (t-1)\tau_i+t.\formula$$
The {\it Kazhdan-Lusztig basis} $\{C_w'\ |\ \tilde w\in \tilde W\}$
of $\tilde H$ is defined [KL] by
$$\overline C_w' = C_w'
\qquad\hbox{and}\qquad
C_w' = t^{-\ell(w)/2}\left(\sum_{y\le w} P_{yw}\tau_y\right),
\formula$$
subject to
$P_{yw}\in \ZZ[t^{1\over2},t^{-{1\over 2}}],
\quad P_{ww}=1,
\quad\hbox{and}\quad
{\rm deg}_t(P_{yw}) \le \hbox{$1\over2$}(\ell(w)-\ell(y)-1).$
\hfill\break
If
$$p_{yw}= q^{-(\ell(w)-\ell(y))}P_{yw}\formula$$
then
$$ C_w' = q^{-\ell(w)}
\sum_{y\le w} P_{yw} q^{\ell(y)}T_y
= \sum_{y\le w} P_{yw} q^{-(\ell(w)-\ell(y))}T_y
= \sum_{y\le w} p_{yw}T_y,
\formula$$
with
$$p_{yw}\in \ZZ[q,q^{-1}],\quad
p_{ww}=1,
\qquad\hbox{and}\qquad
p_{yw}\in q^{-1}\ZZ[q^{-1}],
\formula$$
since ${\rm deg}_q(P_{yw}(q)q^{-(\ell(w)-\ell(y))})
\le \ell(w)-\ell(y)-1-(\ell(w)-\ell(y)) = -1.$
The following proposition establishes the existence and
uniqueness of the $C_w'$ and the $p_{yw}$.

\prop {\sl Let $(\tilde W,\le)$ be a partially ordered set
such that for any $u,v\in \tilde W$ the interval 
$[u,v] = \{ z\in \tilde W\ |\ u\le z\le v\}$ is finite.
Let $M$ be a free $\ZZ[q,q^{-1}]$-module with basis
$\{T_w\ |\ w\in \tilde W\}$ and with a $\ZZ$-linear involution
$\overline{\phantom{T}}\colon M\to M$ such that 
$$\overline{q} = q^{-1}
\qquad\hbox{and}\qquad \overline{T_w} = T_w + \sum_{v<w} a_{vw} T_v.$$
Then there is a unique basis $\{C'_w\ |w\in \tilde W\}$ of $M$ such that
\smallskip\noindent
\itemitem{(a)} $\overline{C'_w} = C'_w$,
\smallskip\noindent
\itemitem{(b)} $\displaystyle{ C_w' = T_w + \sum_{v<w} p_{vw} T_v }$,
\quad with $p_{vw}\in q^{-1}\ZZ[q^{-1}]$ for $v<w$.
}
\pf
The $p_{vw}$ are determined by induction as follows.
Fix $v,w\in W$ with $v<w$. 
If $v=w$ then $p_{vw}=p_{ww}=1$.  For the induction step
assume that $v<w$ and that $p_{zw}$ are known for all $v<z\le w$.

The matrices $A = (a_{vw})$ and $P = (p_{vw})$ are
upper triangular with $1$'s on the diagonal.
The equations
$$\eqalign{
T_w &= \overline{\overline{T_w}} = \sum_v \overline{a_{vw}T_v}
=\sum_{u,v} a_{uv}\overline{a_{vw}}T_u 
\qquad\hbox{and}\cr
\sum_u p_{uw}T_u &= C_w'
=\overline{C_w'} = \sum_v \overline{p_{vw}T_v}
=\sum_v \overline{p_{vw}}a_{uv}T_u, \cr
} $$
imply $A\overline{A} = {\rm Id}$ and $P = A\overline{P}$.  Then 
$$f = \sum_{u<z\le w} a_{uz}\overline{p}_{zw}
=
((A-1)\overline{P})_{uw} = (A\overline P-\overline{P})_{uw}
=(P-\overline{P})_{uw} = p_{uw}-\overline{p}_{uw},$$
is a known element of $\ZZ[q,q^{-1}]$;
$$f = \sum_{k\in \ZZ} f_k q^k
\qquad\hbox{such that}\qquad
\overline{f} 
= \overline{(p_{uw}-\overline{p}_{uw})}
= \overline{p}_{uw}-p_{uw} = - f.$$
Hence $f_k = -f_{-k}$ for all $k\in \ZZ$ and $p_{uw}$ is given by
$p_{uw} =\displaystyle{\sum_{k\in \ZZ_{<0}} f_k q^k }$.
\endpf

The {\it finite Hecke algebra} $H$ and the 
{\it group algebra of $P$} are the subalgebras of $\tilde H$ 
given by
$$\eqalign{
H &= (\hbox{subalgebra of $\tilde H$ generated by $T_1,\ldots, T_n$)},
\qquad\hbox{and}\cr
\KK[P]&=\hbox{$\KK$-span}~\{x^\lambda\ |\ \lambda\in P\},
\qquad\hbox{where $\KK = \ZZ[q,q^{-1}]$},\cr
} \formula$$
respectively.  The Weyl group $W$ acts on $\KK[P]$ by
$$wf = \sum_{\mu\in P} c_\mu x^{w\mu},
\qquad\hbox{for $w\in W$ and
$\displaystyle{f = \sum_{\mu\in P} c_\mu x^\mu\in \KK[P]}$.}
\formula$$

\thm {\sl The center of the affine Hecke algebra is the ring
$$Z(\tilde H) = 
\KK[P]^W = \{ f\in \KK[P]\ |\ wf=f \hbox{\ for all\ } w\in W\}$$
of symmetric functions in $\KK[P]$.
}
\pf
If $z\in \KK[P]^W$ then by the fourth relation in (1.23)
$T_iz = (s_iz)T_i + (q-q^{-1})(1-x^{-\alpha_i})^{-1}(z-s_iz)
=zT_i+0$, for $1\le i\le n$, and by the third relation in (1.23)
$zx^\lambda=x^\lambda z$, for all $\lambda\in P$.  Thus
$z$ commutes with all the generators of $\tilde H$ and so $z\in Z(\tilde H)$.

Assume
$$
z=\sum_{\lambda\in P,w\in W} c_{\lambda,w}x^\lambda T_w\in Z(\tilde H).
$$
Let $m\in W$ be maximal in Bruhat order subject to $c_{\gamma,m}\neq0$ for
some $\gamma\in P$.
If $m\ne 1$ there exists a dominant $\mu\in P$ such that
$c_{\gamma+\mu-m\mu,m}=0$ (otherwise  $c_{\gamma+\mu-m\mu,m}\neq 0$ for every
dominant $\mu\in P$, which is impossible since $z$ is a finite linear
combination of $x^\lambda T_w$). Since $z\in Z(\tilde H)$
we have
$$
z = x^{-\mu} z x^\mu =
\sum_{\lambda\in P,w\in W} c_{\lambda,w} x^{\lambda-\mu} T_w x^\mu.
$$
Repeated use of the third relation in (1.21) yields
$$
T_w x^\mu=\sum_{\nu\in P,v\in W} d_{\nu,v}x^\nu T_v
$$
where $d_{\nu,v}$ are constants such that
$d_{w\mu,w}=1$, $d_{\nu,w}=0$ for
$\nu\ne w\mu$, and $d_{\nu,v}=0$ unless $v\le w$.
So
$$
z=\sum_{\lambda\in P,w\in W} c_{\lambda,w}x^\lambda T_w
= \sum_{\lambda\in P,w\in W} \sum_{\nu\in P,v\in W}
c_{\lambda,w} d_{\nu,v} x^{\lambda-\mu+\nu} T_v
$$
and comparing the coefficients of $x^\gamma T_m$ gives
$
c_{\gamma,m}=c_{\gamma+\mu-m\mu,m} d_{m\mu,m}.
$
Since $c_{\gamma+\mu-m\mu,m}=0$ it follows that
$c_{\gamma,m}=0$, which is a contradiction.  Hence
$z=\sum_{\lambda\in P} c_\lambda x^\lambda\in \KK[P]$.

The fourth relation in (1.23) gives
$$
zT_i=T_iz=(s_iz)T_i+(q-q^{-1})z'
$$
where $z'\in\KK[P]$. Comparing coefficients of
$x^\lambda$  on both sides yields $z' = 0$. Hence
$zT_i=(s_iz)T_i$, and therefore $z=s_iz$ for $1\leq i\leq n$.  So
$z\in \KK[P]^W$.
\endpf

\vfill\eject

\section 2. Symmetric and alternating functions and their q-analogues

\bigskip
Let ${\bf 1}_0$ and $\varepsilon_0$
be the elements of the finite Hecke algebra $H$ 
which are determined by
$$\matrix{
{\bf 1}_0^2= {\bf 1}_0\qquad &\hbox{and}\qquad
&T_i{\bf 1}_0=q{\bf 1}_0, \hfill \quad
&\hbox{for all $1\le i\le n$,} \hfill \cr$$
\varepsilon_0^2= \varepsilon_0\qquad &\hbox{and}\qquad
&T_i\varepsilon_0=(-q^{-1})\varepsilon_0,\hfill 
&\hbox{for all $1\le i\le n$.} \hfill \cr
}$$
In terms of the basis $\{ T_w\ |\ w\in W\}$ of $H$
these elements have the explicit formulas
$${\bf 1}_0 = {1\over W_0(q^2)}\sum_{w\in W} q^{\ell(w)}T_w\,,
\quad\hbox{and}\quad
\varepsilon_0 = {1\over W_0(q^{-2})}\sum_{w\in W} (-q)^{-\ell(w)}T_w\,,
\formula
$$
where
$W_0(t) = \sum_{w\in W} t^{\ell(w)}$.
(To define these elements one should adjoin
the element $W_0(q^2)^{-1}$ to $\KK$ or to $H$.)
The elements ${\bf 1}_0$ and $\varepsilon_0$ are $q$-analogues
of the elements in the group algebra of $W$ given by
$${\bf 1} = {1\over |W|}\sum_{w\in W} w
\qquad\hbox{and}\qquad
\varepsilon = {1\over |W|} \sum_{w\in W} (-1)^{\ell(w)}w,
\formula
$$
and the vector spaces ${\bf 1}_0\tilde H{\bf 1}_0$ and
and $\varepsilon_0\tilde H{\bf 1}_0$ are $q$-analogues of the
vector spaces (more precisely, free $\KK=\ZZ[q,q^{-1}]$-modules)
of {\it symmetric functions} and {\it alternating functions},
$$\eqalign{
\KK[P]^W &= 
\{ f\in \KK[P]\ |\ wf=f \hbox{\ for all\ } w\in W\} = {\bf 1}\cdot\KK[P], \cr
{\cal A} &= 
\{ f\in \KK[P]\ |\ wf=(-1)^{\ell(w)}f \hbox{\ for all\ } w\in W\} 
= \varepsilon\cdot\KK[P], \cr
}
\formula$$
respectively.

For $\mu\in P$ let $W\mu = \{ w\mu\ |\ w\in W\}$ be the orbit of $\mu$
and $W_\mu = \{ w\in W\ |\ w\mu=\mu\}$ the stabilizer of $\mu$ and define
$$\matrix{
\displaystyle{m_\mu =\sum_{\gamma\in W\mu} x^\gamma
= {1\over |W_\mu|}\,{\bf 1}\cdot x^\mu\,,} \hfill
&\qquad\qquad
&\displaystyle{
a_\mu = \sum_{w\in W} (-1)^{\ell(w)} wx^\mu = \varepsilon\cdot x^\mu\,,} \cr
\cr
\cr
M_\mu = {\bf 1}_0 x^\mu {\bf 1}_0\,, \hfill
&&A_\mu = \varepsilon_0 x^\mu {\bf 1}_0\,. \hfill \cr
}
\formula$$
Theorem 2.7 below shows that the
elements in (2.4) which are indexed by elements of $P^+$ and $P^{++}$
form bases (over $\KK$) of 
$\KK[P]^W$, ${\cal A}$, ${\bf 1}_0\tilde H{\bf 1}_0$,
and $\varepsilon_0\tilde H{\bf 1}_0$.  This will be a consequence
of the following {\it straightening rules}.  The straightening
law for the $M_\mu$ given in the following Proposition is a 
generalization of [Mac, III \S 2 Ex.\ 2].

\prop {\sl For $\gamma\in P$ let
$m_\gamma$, $a_\gamma$, $M_\gamma$,  and $A_\gamma$
be as defined in (2.4).
Let $\alpha_i$ be a simple root and let $\mu\in P$ be such that 
$d=\langle \mu,\alpha_i^{\vee}\rangle \ge 0$.  Then
$$m_{s_i\mu}=m_\mu,
\qquad
a_{s_i\mu}=-a_\mu, \qquad\hbox{and}\qquad
A_{s_i\mu}=-A_{\mu}. 
$$
Letting $t=q^{-2}$, $M_\mu = M_{s_i\mu}$ if $d=0$, and if $d>0$ then
$$M_{s_i\mu} = tM_{\mu}
+\Big(\sum_{j=1}^{\lfloor d/2-1\rfloor}
(t^2-1)t^{j-1}M_{\mu-j\alpha_i}\Big)
+\cases{(t-1)t^{d/2-1}
M_{\mu-(d/2)\alpha_i}, &if $d$ is even, \cr
0, &if $d$ is odd.\cr
}$$
}
\pf  The first two equalities follow from 
the definitions of $m_\lambda$ and $a_\mu$
and the fact that $\ell(s_i)=1$.

Let $\mu\in P$ such that
$d=\langle \mu,\alpha_i^{\vee}\rangle \ge 0$.  
Since $x^\mu+x^{s_i\mu}$ is in the center of the tiny
little affine Hecke algebra generated by
$T_i$ and the $x^\gamma$, $\gamma\in P$,
$$\eqalign{
A_\mu+A_{s_i\mu} 
&= \varepsilon_0 (x^\mu+x^{s_i\mu}) {\bf 1}_0
=q^{-1} \varepsilon_0 (x^\mu+x^{s_i\mu})T_i {\bf 1}_0 \cr
&=q^{-1}\varepsilon_0 T_i(x^\mu+x^{s_i\mu}) {\bf 1}_0
=-q^{-2}\varepsilon_0 (x^\mu+x^{s_i\mu}) {\bf 1}_0
=-q^{-2}(A_\mu+A_{s_i\mu}). \cr
}$$
Thus $A_\mu+A_{s_i\mu}=0$ which establishes the third statement.

If $d=0$ then $s_i\mu = \mu$ and the fourth relation in (1.23) is
$x^{s_i\mu} T_i - T_i x^\mu = 0$.  Multiplying by ${\bf 1}_0$ on both the
left and the right (and dividing by $q$) gives 
${\bf 1}_0x^{s_i\mu}{\bf 1}_0 - {\bf 1}_0x^\mu{\bf 1}_0$ as desired.
If $d>0$ then multiplying the fourth relation in (1.23)
by ${\bf 1}_0$ on both the left and the right
(and then multiplying by $q^{-1}$) gives
$${\bf 1}_0(x^{s_i\mu}-x^{\mu}){\bf 1}_0
= q^{-1}(q-q^{-1}){\bf 1}_0
\left({x^{s_i\mu}-x^\mu\over 1-x^{-\alpha_i}} \right){\bf 1}_0.$$
Subtracting the same relation with $\mu$ replaced by $\mu-\alpha_i$
gives
$$\eqalign{
{\bf 1}_0(x^{s_i\mu}-x^\mu){\bf 1}_0
-{\bf 1}_0(x^{s_i\mu+\alpha_i}-x^{\mu-\alpha_i}){\bf 1}_0
&=(1-q^{-2}){\bf 1}_0\left(
{x^{s_i\mu}-x^\mu-x^{s_i\mu+\alpha_i}+x^{\mu-\alpha_i}
\over 1-x^{-\alpha_i} }\right){\bf 1}_0 \cr
&= (1-q^{-2}) {\bf 1}_0 (-x^{s_i\mu+\alpha_i}-x^\mu){\bf 1}_0. \cr
}$$
So
$${\bf 1}_0x^{s_i\mu}{\bf 1}_0
= q^{-2}{\bf 1}_0x^\mu{\bf 1}_0-{\bf 1}_0x^{\mu-\alpha_i}{\bf 1}_0
+q^{-2}{\bf 1}_0x^{s_i\mu+\alpha_i}{\bf 1}_0.
$$
Inductively applying this relation yields the result.
The first cases are
$$M_{s_i\mu} = \cases{
M_\mu, &if $\langle\mu,\alpha_i^\vee\rangle=0$, \cr
q^{-2}M_\mu, &if $\langle\mu,\alpha_i^\vee\rangle=1$, \cr
q^{-2}M_\mu+(q^{-2}-1)M_{\mu-\alpha_i}, 
&if $\langle\mu,\alpha_i^\vee\rangle=2$, \cr
q^{-2}M_\mu+(q^{-4}-1)M_{\mu-\alpha_i}, 
&if $\langle\mu,\alpha_i^\vee\rangle=3$, \cr
q^{-2}M_\mu+(q^{-4}-1)M_{\mu-\alpha_i} + q^{-2}(q^{-2}-1)M_{\mu-2\alpha_i}, 
&if $\langle\mu,\alpha_i^\vee\rangle=4$.\qquad\hbox{\qed} \cr
}
$$

\medskip\noindent
Proposition 2.5 implies that, 
for all $\mu\in P$ and $w\in W$,
$$m_{w\mu} = m_\mu,
\qquad
a_{w\mu} = (-1)^{\ell(w)}a_\mu,
\qquad\hbox{and}\qquad
A_{w\mu} = (-1)^{\ell(w)}A_\mu.
\formula
$$

\thm   {\sl  Let $\KK = \ZZ[q,q^{-1}]$.  As free $\KK$-modules
$$\matrix{
\KK[P]^W \hfill 
&&\hbox{has basis} 
~~\{ m_\lambda\ |\ \lambda\in P^+\}, \hfill 
&\qquad\qquad
&{\bf 1}_0\tilde H{\bf 1}_0 \hfill
&&\hbox{has basis} 
~~\{ M_\lambda \ |\ \lambda\in P^+\}, \hfill $$
\cr
{\cal A} \hfill
&&\hbox{has basis} 
~~\{ a_\mu\ |\ \mu\in P^{++}\}, \hfill 
&&\varepsilon_0\tilde H{\bf 1}_0 \hfill
&&\hbox{has basis} 
~~\{ A_\mu \ |\ \mu\in P^{++}\}. \hfill
\cr}$$
}
\pf
Since $\{x^\mu T_w\ |\ \mu\in P, w\in W\}$ form a basis
of $\tilde H$ the elements $M_\mu = {\bf 1}_0x^\mu{\bf 1}_0
=q^{-\ell(w)}{\bf 1}_0 x^\mu T_w {\bf 1}_0$, $\mu\in P$, span 
${\bf 1}_0\tilde H{\bf 1}_0$.
By Proposition 2.5, if $\mu$ is on the negative side
of a hyperplane $H_{\alpha_i}$, i.e. if 
$\langle \mu,\alpha_i^\vee\rangle<0$, then $M_\mu$ can be
rewritten as a linear combination of $M_\gamma$ such that all
terms have $\gamma$ on the nonnegative side of $H_{\alpha_i}$.  
By repeatedly applying the relation in Proposition 2.5
$M_\mu$ can be rewritten as a linear combination
of $M_\gamma$ such that all terms have $\gamma$ on the
nonnegative side of $H_{\alpha_1},\ldots, H_{\alpha_n}$, i.e.
$\gamma\in P^+=P\cap \overline{C}$, where
$\overline C = \{ x\in \RR^n\ |\ \langle x,\alpha_i^\vee\rangle\ge 0
\hbox{\ for all $1\le i\le n$}\}$.

The proof for the cases of $m_\mu$, $a_\mu$ and $A_\mu$ is easier,
it follows directly from (2.6),  the fact that 
$C = \{ x\in \RR^n\ |\ \langle x,\alpha_i^\vee\rangle>0
\hbox{\ for all $1\le i\le n$}\}$ is a fundamental chamber for the
action of $W$, and that if
$\mu\in P^+\backslash P^{++}$ then $\langle\mu,\alpha_i^\vee\rangle=0$
and $a_\mu=-a_{s_i\mu}=-a_\mu$, in which case $a_\mu=0$ 
(similarly for $A_\mu$).
\endpf

\medskip
For $\lambda\in P$ define the {\it Schur function}, 
or {\it Weyl character}, by
$$s_\lambda = {a_{\lambda+\rho}\over a_\rho},
\qquad\hbox{where}\qquad
\rho = \hbox{$1\over2$}\sum_{\alpha\in R^+}\alpha.
\formula$$
The straightening law for $a_\mu$ in (2.6) implies
the following straightening law for the Schur functions.
If $\mu\in P$ and $w\in W$ then, by (2.6) and the definition of $s_\mu$, 
$$(-1)^{\ell(w)}s_\mu
= {(-1)^{\ell(w)}a_{\mu+\rho}\over a_\rho} 
= {a_{w(\mu+\rho)-\rho+\rho}\over a_\rho} 
=s_{w\circ \mu}, 
\qquad\hbox{where $w\circ\mu = w(\mu+\rho)-\rho$.}
\formula
$$
The {\it dot action} of the Weyl group $W$ on $\fh_\RR^*$ which is 
appearing here,
$w\circ \mu = t_{-\rho}wt_{\rho}\mu = (t_{\rho}^{-1})wt_\rho\mu,$
is the ordinary action of $W$ on $\fh_\RR^*$ except with
the ``center'' shifted to $-\rho$.  For the root system of type $C_2$,
see Example 1.4, the picture is
$$
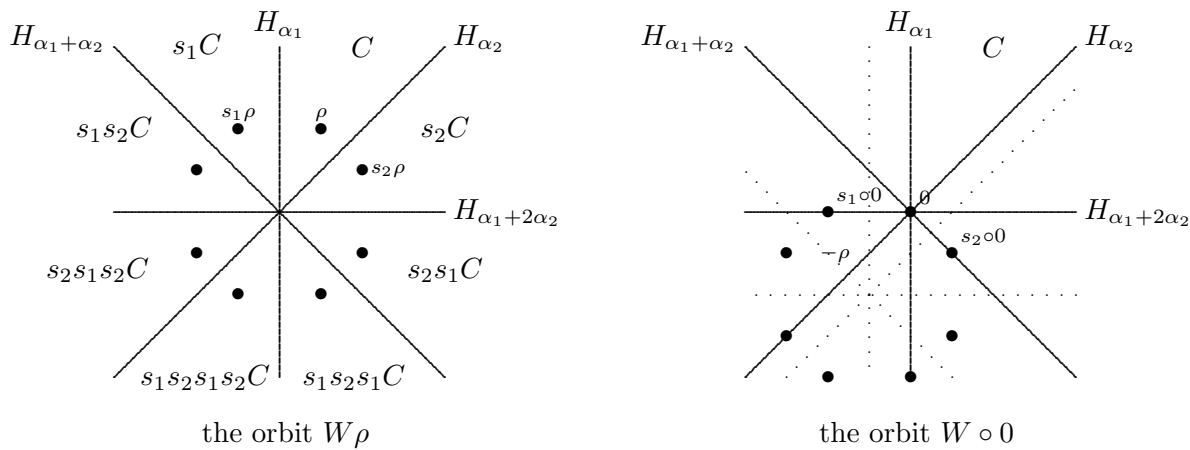
\matrix{
\beginpicture
\setcoordinatesystem units <1cm,1cm>         
\setplotarea x from -2.5 to 2.5, y from -2.5 to 2.5  
\put{$H_{\alpha_1}$}[b] at 0 2.1
\put{$H_{\alpha_2}$}[l] at 2.1 2.1
\put{$H_{\alpha_1+\alpha_2}$}[r] at -2.1 2.1
\put{$H_{\alpha_1+2\alpha_2}$}[l] at 2.1 0
\put{$C$} at 1 2
\put{$s_1C$} at -1 2
\put{$s_2C$} at  2 1
\put{$s_1s_2C$} at  -2 1
\put{$s_2s_1C$} at   2  -0.7 
\put{$s_1s_2s_1C$} at  0.9 -2
\put{$s_2s_1s_2C$} at  -2.2 -0.7
\put{$s_1s_2s_1s_2C$} at  -0.9 -2
\put{$\bullet$} at   0.5  1
\put{$\bullet$} at  -0.5  1
\put{$\bullet$} at   1    0.5
\put{$\bullet$} at  -1    0.5
\put{$\bullet$} at   0.5  -1
\put{$\bullet$} at  -0.5  -1
\put{$\bullet$} at   1 -0.5
\put{$\bullet$} at  -1 -0.5
\put{$\scriptstyle{\rho}$}[b] at 0.5 1.1
\put{$\scriptstyle{s_1\rho}$}[b] at -0.5 1.1
\put{$\scriptstyle{s_2\rho}$}[l] at  1.1 0.5 
\plot -2 -2   2 2 /
\plot  2 -2  -2 2 /
\plot  0  2   0 -2 /
\plot  2  0  -2  0 /
\endpicture
&\phantom{T} &
\beginpicture
\setcoordinatesystem units <1cm,1cm>         
\setplotarea x from -2.5 to 2.5, y from -2.5 to 2.5  
\put{$H_{\alpha_1}$}[b] at 0 2.1
\put{$H_{\alpha_2}$}[l] at 2.1 2.1
\put{$H_{\alpha_1+\alpha_2}$}[r] at -2.1 2.1
\put{$H_{\alpha_1+2\alpha_2}$}[l] at 2.1 0
\put{$C$} at 1 2
\put{$\bullet$} at   0  0
\put{$\bullet$} at  -1  0
\put{$\bullet$} at   0.5  -0.5
\put{$\bullet$} at  -1.5  -0.5
\put{$\bullet$} at   0.5  -1.5
\put{$\bullet$} at  -1.5  -1.5
\put{$\bullet$} at   0 -2
\put{$\bullet$} at  -1 -2
\put{$\scriptstyle{-\rho}$}[tl] at -1.1 -0.4
\put{$\scriptstyle{0}$}[bl] at 0.1 0.1
\put{$\scriptstyle{s_1\circ 0}$}[bl] at -0.9  0.1
\put{$\scriptstyle{s_2\circ 0}$}[bl] at  0.6  -0.4 
\plot -2 -2   2 2 /
\plot  2 -2  -2 2 /
\plot  0  2   0 -2 /
\plot  2  0  -2  0 /
\setdots
\plot -1.5 -2   2 1.5 /
\plot  0.5 -2   -2 0.5 /
\plot  -0.5 2  -0.5 -2 /
\plot  2 -1  -2 -1 / 
\endpicture
\cr
\hbox{the orbit $W\rho$} 
&&\hbox{the orbit $W\circ 0$} \cr
}
$$

The following proposition shows that the Weyl characters
$s_\lambda$ are elements of $\KK[P]^W$.  The equality in part (a)
is the {\it Weyl denominator formula}, a generalization
of the factorization of the Vandermonde determinant
$~\det(x_i^{n-j})=\prod_{1\le i,j\le n} (x_i-x_j)$.
In the remainder of this section
we shall abuse language and use the term
``vector space'' in place of ``free $\KK=\ZZ[q,q^{-1}]$ module''.

\prop {\sl Let $P^+$, $P^{++}$, $\KK[P]^W$ and ${\cal A}$ be as in (1.8)
and (2.4) and let $\rho$ be as in (1.9).
\smallskip\noindent
\item{(a)}  
If $f\in {\cal A}$ then $f$ is divisible by $a_\rho$ and
$$
a_\rho = x^\rho\prod_{\alpha\in R^+} (1-x^{-\alpha})
$$
\smallskip\noindent
\item{(b)} The set $\{ s_\lambda\ |\ \lambda\in P^+\}$ is a basis 
of $\KK[P]^W$.
\smallskip\noindent
\item{(c)}   
The maps
$$\matrix{ 
P^+ &\longrightarrow &P^{++} \cr
\lambda &\longmapsto &\lambda+\rho \cr
}
\qquad\hbox{and} \qquad
\matrix{
\Phi\colon &\KK[P]^W &\longrightarrow &{\cal A} \cr
&f &\longmapsto &a_\rho f \cr
&s_\lambda &\longmapsto &a_{\lambda+\rho} \cr
}
$$
are a bijection and a vector space isomorphism, respectively.
}
\pf 
Since $s_i$ takes $\alpha_i$ to $-\alpha_i$ and permutes the
other elements of $R^+$,
$$\rho-\langle \rho,\alpha_i^\vee\rangle\alpha_i
=s_i\rho = \rho-\alpha_i,
\qquad\hbox{and so}\qquad
\langle \rho,\alpha_i^\vee\rangle = 1,
\quad\hbox{for all $1\le i\le n$.}
$$
Thus the map $P^+\to P^{++}$ given by $\lambda\mapsto \lambda+\rho$
is well defined and it is a bijection since it is invertible.

Let $d = x^\rho\prod_{\alpha\in R^+} (1-x^{-\alpha})
= \prod_{\alpha\in R^+} (x^{\alpha/2}-x^{-\alpha/2})$.
Since $s_i$ takes $\alpha_i$ to $-\alpha_i$ and
permutes the other elements of $R^+$, 
$s_id =-d$ for all $1\le i\le n$ and so
$wd = (-1)^{\ell(w)}d$ for all $w\in W$.  Thus $d$ is an element
of ${\cal A}$.  

If $\alpha\in R^+$ and  
$\displaystyle{
f = \sum_{\mu\in P} c_\mu x^\mu \in {\cal A} }$ then
$$\sum_{\mu\in P} c_\mu x^\mu = 
f = -s_\alpha f = \sum_{\mu\in P} -c_\mu x^{s_\alpha\mu},
\qquad\hbox{and so}\qquad
f = 
\sum_{\langle \mu,\alpha^\vee\rangle \ge 0} c_\mu(x^\mu-x^{s_\alpha\mu}).$$
Since $(1-x^{-\langle \mu,\alpha^\vee\rangle\alpha})$
is divisible by $(1-x^{-\alpha})$ it follows that 
$x^\mu-x^{s_i\mu}=x^\mu(1-x^{-\langle \mu,\alpha^\vee\rangle\alpha})$
is divisible by $(1-x^{-\alpha})$ and thus that
$f$ is divisible by $(1-x^{-\alpha})$ for all $\alpha\in R^+$.
Since the elements $(1-x^{-\alpha})$ are relatively prime in 
the Laurent polynomial ring $\KK[P]$ and $x^\rho$ is a unit 
in $\KK[P]$, $f$ is divisible by $d$.  Since
both $f$ and $d$ are in ${\cal A}$, the quotient $f/d$ is an
element of $\KK[P]^W$.

The monomial $x^\rho$ appears in $a_\rho$ with coefficient $1$
and it is the unique term $x^\mu$ in $a_\rho$ with $\mu\in P^+$.
Since $d$ has highest term $x^\rho$ with coefficient $1$ and
$a_\rho$ is divisible by $d$ it follows that $a_\rho/d=1$.
Thus $a_\rho=d$, the inverse of the map $\Phi$ in (c) is well defined,
and $\Phi$ is an isomorphism.

Since $\{a_{\lambda+\rho}\ |\ \lambda\in P^+\}$ is a basis
of ${\cal A}$ and the map $\Phi$ is an isomorphism it follows that
$\{s_\lambda\ |\ \lambda\in P^+\}$ is a $\KK$-basis of $\KK[P]^W$.
\endpf

\subsection The Satake isomorphism

The following theorem establishes a $q$-analogue
of the isomorphism $\Phi$ from Proposition 2.10(c).
The map $\Phi_1$ in the following theorem is 
the {\it Satake isomorphism}.
We shall continue to abuse language and use the term
``vector space'' in place of ``free $\KK=\ZZ[q,q^{-1}]$ module''.

\thm {\sl 
The vector space isomorphism $\Phi$ in Proposition 2.10(c) generalizes to
a vector space isomorphism
$$\matrix{
\tilde \Phi\colon &Z(\tilde H)=\KK[P]^W &\mapright{\Phi_1}
&Z(\tilde H){\bf 1}_0 = {\bf 1}_0\tilde H{\bf 1}_0 &\mapright{\Phi_2}
&\varepsilon_0\tilde H{\bf 1}_0 \cr
&f &\longmapsto &f{\bf 1}_0 &\longmapsto &A_\rho f {\bf 1}_0 \cr
&s_\lambda &\longmapsto &s_\lambda{\bf 1}_0 &\longmapsto
&A_{\lambda+\rho}. \cr
}
$$
}
\pf 
Using the third equality in (2.6),
$$\varepsilon_0 a_\lambda {\bf 1}_0
=\varepsilon_0\left(\sum_{w\in W} (-1)^{\ell(w)}x^{w\lambda}\right)
{\bf 1}_0
=\sum_{w\in W} (-1)^{\ell(w)} A_{w\lambda} = |W|\,A_{\lambda}.$$
By Proposition 2.10(c) and Theorem 1.33, $s_\lambda\in \KK[P]^W = Z(\tilde H)$,
and so
$$A_\rho s_\lambda {\bf 1}_0
= {1\over |W|} \varepsilon_0 a_\rho {\bf 1}_0 s_\lambda {\bf 1}_0
={1\over |W|} \varepsilon_0 a_\rho s_\lambda {\bf 1}_0^2
={1\over |W|} \varepsilon_0 a_{\lambda+\rho} {\bf 1}_0
= A_{\lambda+\rho}.
$$
Since $\{ s_\lambda\ |\ \lambda\in P^+\}$ is a basis of 
$\KK[P]^W =Z(\tilde H)$
and $\{ A_{\lambda+\rho}\ |\ \lambda\in P^+\}$ is a basis of
$\varepsilon_0 \tilde H {\bf 1}_0$ the composite map
$$
\matrix{
Z(\tilde H) &\mapright{{\bf 1}_0} &Z(\tilde H)\,{\bf 1}_0
&\hookrightarrow &{\bf 1}_0\tilde H{\bf 1}_0
&\mapright{A_\rho} &\varepsilon_0\tilde H{\bf 1}_0 \cr
f &\longmapsto &f\,{\bf 1}_0 &\mapsto &f\,{\bf 1}_0
&\longmapsto &A_\rho\,f\,{\bf 1}_0 \cr
s_\lambda &\longmapsto &s_\lambda\,{\bf 1}_0
&\mapsto &s_\lambda\,{\bf 1}_0 &\longmapsto &A_{\lambda+\rho} \cr
}$$
is a vector space isomorphism.
\endpf

If $\mu\in P$ let
$$W_\mu = \{ w\in W\ |\ w\mu = \mu\}
\qquad\hbox{and}\qquad
W_\mu(t) = \sum_{w\in W_\mu} t^{\ell(w)}.
\formula$$
In particular, if $\mu=0$, then $W_0=W$ and $W_0(t)$ is
the polynomial that appears in (2.1).  

The {\it Hall-Littlewood polynomials}, or
{\it Macdonald spherical functions}, are defined by 
$$
P_\mu(x;t)
= 
{1\over W_\mu(t)}
\sum_{w\in W} w\Big( x^\mu 
\prod_{\alpha\in R^+} {1-tx^{-\alpha}\over 1-x^{-\alpha}}\Big),
\qquad\hbox{for $\mu\in P$}.
\formula$$
Then $m_\mu = P_\mu(x;1)$ and, using the Weyl denominator
formula, 
$$P_\mu(x;0) = \sum_{w\in W} w\Big(
{x^\rho x^\mu\over x^\rho\prod_{\alpha\in R^+} (1-x^{-\alpha})}
\Big)
={1\over a_\rho}\sum_{w\in W} (-1)^{\ell(w)} wx^{\mu+\rho}
={a_{\mu+\rho}\over a_\rho} = s_\mu,
\formula$$
and so, conceptually, the spherical functions
$P_\mu(x;t)$ interpolate between the
Schur functions $s_\mu$ and the monomial symmetric functions
$m_\mu$.

The double cosets in $W\backslash \tilde W/W$ are 
$Wt_\lambda W$, $\lambda\in P^+$.
If $\lambda\in P^+$ let $n_\lambda$ and $m_\lambda$
be the maximal and minimal length elements of $Wt_\lambda W$,
respectively.
Theorem 2.22 below will show that under the Satake isomorphism
the Weyl characters $s_\lambda$
correspond to Kazhdan Lusztig basis elements $C_{n_\lambda}'$
and the polynomials $P_\mu(x,q^{-2})$
correspond to the elements $M_\mu = {\bf 1}_0x^\mu{\bf 1}_0$.
More precisely, we have the following diagram:
$$\matrix{
\Phi_1\colon &Z(\tilde H)=\KK[P]^W &\longrightarrow
&Z(\tilde H){\bf 1}_0 = {\bf 1}_0\tilde H{\bf 1}_0 \cr
\cr
&f &\longmapsto &f{\bf 1}_0 \cr
\cr
&q^{-\ell(w_0)}W_0(q^2)\,s_\lambda\hfill &\longmapsto &C'_{n_\lambda} \cr
\cr
&\displaystyle{ {W_\mu(q^{-2})\over W_0(q^{-2})}\,P_\mu(x;q^{-2})}\hfill 
&\longmapsto &M_\mu  \cr
}
\formula$$
where $w_0$ is the longest element of $W$.
The following three lemmas (of independent interest)
are used in the proof of Theorem 2.22.

\lemma {\sl Let $t_\alpha$, $\alpha\in R^+$, be commuting variables
indexed by the positive roots. 
For $\lambda\in P^+$ let $P_\lambda(x;t)$ be as in (2.13),
$W_\lambda$ as in (2.12), and define
$$R_\lambda(x; t_\alpha) = \sum_{w\in W} w
\left(x^\lambda\prod_{\alpha\in R^+} 
{1-t_\alpha x^{-\alpha}\over 1-x^{-\alpha} }\right)
\qquad\hbox{and}\qquad
W_\lambda(t_\alpha) = 
\sum_{w\in W_\lambda} \Big(\prod_{\alpha\in R(w)} t_\alpha\Big),
$$
where, as in (1.6),
$R(w) = \{ \alpha\in R^+\ |\ w\alpha<0\}$ is the inversion
set of $w$.  Then
\medskip\noindent
\item{(a)}
$\displaystyle{
R_\lambda(x; t_\alpha) = \sum_{\mu\in P^+} u_{\lambda\mu} s_\mu }$,
\quad
with $u_{\lambda\mu}\in \ZZ[t_\alpha]$,\quad 
$u_{\lambda\mu}=0$ unless $\mu\le \lambda$,\quad and \quad
$u_{\lambda\lambda} = W_\lambda(t_\alpha)$. 
\smallskip\noindent
\item{(b)}
$\displaystyle{
P_\lambda(x; t) 
= \sum_{\mu\in P^+} c_{\lambda\mu} s_\mu }$,
\quad
with $c_{\lambda\mu}\in \ZZ[t]$, \quad
$c_{\lambda\mu}=0$ unless $\mu\le \lambda$,\quad and \quad
$c_{\lambda\lambda} = 1.$
}
\pf
(a)  If $E\subseteq R^+$ let 
$$t_E = \prod_{\alpha\in E} t_\alpha 
\qquad\hbox{and}\qquad \alpha_E = \sum_{\alpha\in E} \alpha,$$  
and let $a_\mu$ be as defined in (2.4).
Using the Weyl denominator formula, Proposition 2.10(a),
$$\eqalign{
R_\lambda &=
\sum_{w\in W} w\left( x^\lambda\prod_{\alpha\in R^+} 
{1-t_\alpha x^{-\alpha}\over 1-x^{-\alpha} } \right) 
= \sum_{w\in W} w\left(
{ x^{\lambda+\rho}\prod_{\alpha\in R^+} (1-t_\alpha x^{-\alpha})
\over x^{\rho} \prod_{\alpha\in R^+} (1-x^{-\alpha}) } \right) \cr
&= {1\over a_\rho}\sum_{w\in W} (-1)^{\ell(w)}w\left(
x^{\lambda+\rho}\prod_{\alpha\in R^+} (1-t_\alpha x^{-\alpha})\right)
\cr
&= {1\over a_{\rho}}\sum_{w\in W} (-1)^{\ell(w)} w\left(
\sum_{E\subseteq R^+} (-1)^{|E|} t_E
x^{\lambda+\rho-\alpha_E}\right)
\cr
&= {1\over a_{\rho}}
\sum_{E\subseteq R^+} (-1)^{|E|} t_E
a_{\lambda+\rho-\alpha_E}
=\sum_{E\subseteq R^+} (-1)^{|E|} t_E s_{\lambda+\rho-\alpha_E}, \cr
}$$
which shows that $R_\lambda$ is a symmetric function
and $u_{\lambda\mu}\in \ZZ[t_\alpha]$. 

By the straightening law for Weyl characters (2.9),
$s_{\lambda+\rho-\alpha_E}=0$ or
$s_{\lambda+\rho-\alpha_E}=(-1)^{\ell(v)} s_{\mu+\rho}$ 
with 
$$\hbox{$v\in W$ and $\mu\in P^+$ such that} \quad
\mu+\rho = v^{-1}(\lambda+\rho-\alpha_E).$$
Let $E^c$ denote the complement of $E$ in $R^+$.
Since $v$ permutes the elements of $R^+$,
$$\eqalign{
v^{-1}(\lambda+\rho-\alpha_E)
&= v^{-1}\lambda+v^{-1} \left( \hbox{$1\over2$}\sum_{\alpha\in E^c} \alpha 
- \hbox{$1\over2$} \sum_{\alpha\in E}\alpha \right) \cr
&=v^{-1}\lambda +
\left(\hbox{$1\over2$}\sum_{\alpha\in F^c}\alpha - \hbox{$1\over2$}
\sum_{\alpha\in F} \alpha\right)
=v^{-1}\lambda+\rho-\alpha_F, \cr
}$$
for some subset $F\subseteq R^+$ (which could be determined
explicitly in terms of $E$ and $v$).  
Hence
$$\mu = v^{-1}\lambda + \rho-\alpha_F - \rho = v^{-1}\lambda - \alpha_F
\le v^{-1}\lambda\le \lambda.
\eqno{(*)}$$
This proves that $u_{\lambda\mu}=0$ unless $\mu\le \lambda$.

In ($*$), $\mu=\lambda$ only if $v^{-1}\lambda = \lambda$ and 
$\rho = \rho-\alpha_F = v^{-1}(\rho-\alpha_E)$ in which case
$$\rho-\alpha_E = v\left(\hbox{$1\over 2$}\sum_{\alpha\in R^+} \alpha\right)
=\rho-\sum_{\alpha\in R(v)} \alpha
\qquad\hbox{and}\qquad
E= R(v).$$
Thus 
$$u_{\lambda\lambda}(t_\alpha) 
= \sum_{v^{-1}\in W_\lambda} t_{R(v)}.$$
(b)
By applying (a) to $\lambda=0$,
$$R_0(x;t_\alpha) = 
\sum_{w\in W} w\left( \prod_{\alpha\in R^+} 
{1-t_\alpha x^{-\alpha}\over 1-x^{-\alpha} } \right) 
= W_0(t_\alpha).
\eqno(*)$$
Let $W^\lambda$ be a set of minimal length coset representatives
of the cosets in $W/W_\lambda$.  
Every element $w\in W$ can be written uniquely as $w=uv$ with
$u\in W^\lambda$ and $v\in W_\lambda$ (see [Bou, IV \S 1 Ex.\ 3]).
Let 
$$Z(\lambda) = \{ \alpha\in R^+\ |\ \langle\lambda,\alpha^\vee\rangle =0\},$$
and let $Z(\lambda)^c$ be the complement of $Z(\lambda)$ in $R^+$.  
Then $v\in W_\lambda$ permutes the elements of $Z(\lambda)^c$ among
themselves and so
$$\eqalign{
R_\lambda(x;t_\alpha)
&= \sum_{u\in W^\lambda}
u\left(x^\lambda \prod_{\alpha\in Z(\lambda)^c}
{1-t_\alpha x^{-\alpha}\over 1-x^{-\alpha}}
\sum_{v\in W_\lambda} v\Big(
\prod_{\alpha\in Z(\lambda)} {1-t_\alpha x^{-\alpha}\over 1-x^{-\alpha}}
\Big) \right)
\cr
&= \sum_{u\in W^\lambda}
u\left(x^\lambda 
\prod_{\alpha\in Z(\lambda)^c} {1-t_\alpha x^{-\alpha}\over 1-x^{-\alpha}} 
W_\lambda(t_\alpha) \right), \cr
}$$
where the last equality follows from ($*$).
Thus there is an element
$P_\lambda(x;t_\alpha)\in \FF[P]$ where $\FF$ is the
field of fractions of $\ZZ[t_\alpha]$ such that
$$R_\lambda(x;t_\alpha)
=W_\lambda(t_\alpha) 
\sum_{u\in W^\lambda}
u\left(x^\lambda 
\prod_{\alpha\in Z(\lambda)^c} 
{1-t_\alpha x^{-\alpha}\over 1-x^{-\alpha}}\right) 
=W_\lambda(t_\alpha) P_\lambda(x;t_\alpha).
$$
Since $R_\lambda$ is a symmetric polynomial, i.e. an element
of $\ZZ[t_\alpha][P]^W$, $P_\lambda(x;t_\alpha)\in \FF[P]^W$. 
Since the $t_\alpha$ occur only in the numerators of the terms 
in the sum defining $P_\lambda$
in fact $P_\lambda$ is a symmetric polynomial with coefficients
in $\ZZ[t_\alpha]$.
It follows that all the $u_{\lambda\mu}$ appearing in part (a)
are divisible by $W_\lambda(t_\alpha)$ and 
$$P_\lambda(x;t_\alpha)
= \sum_{\mu\in P} c_{\lambda\mu} s_\mu,
\qquad\hbox{where}\quad
c_{\lambda\mu} =
{1\over W_\lambda(t_\alpha)}\,u_{\lambda\mu}$$
are polynomials in $\ZZ[t_\alpha]$ such that 
$c_{\lambda\lambda}=1$ and 
$c_{\lambda\mu}=0$ unless $\mu\le \lambda$.
The result in (b) now follows by specializing $t_\alpha=t$ for 
all $\alpha\in R^+$.
\endpf

\medskip
Lemma 2.16 has the following interesting (and useful) corollary,
see [Mac3].

\cor {\sl  Let $\rho$ and $\alpha^\vee$ be as in (1.9) and (1.1),
respectively and let $W_0(t)$ be as defined in (2.12).
\smallskip\noindent
\item{(a)}
$\displaystyle{
~~\sum_{w\in W} w\left(
\prod_{\alpha\in R^+} {1-tx^{-\alpha}\over 1-x^{-\alpha}}\right)
=W_0(t).
}$
\smallskip\noindent
\item{(b)} 
$\displaystyle{
~~\prod_{\alpha\in R^+} {1-t^{\langle \rho,\alpha^\vee\rangle+1}
\over 1-t^{\langle \rho,\alpha^\vee\rangle} }
=W_0(t). 
}$
}
\pf
(a) follows from Lemma 2.16 part (a) by setting $\lambda=0$
and specializing $t_\alpha=t$ for all $\alpha\in R^+$.
\smallskip\noindent
(b)  Applying the homomorphism
$$\matrix{
\ZZ[t^{\pm1}][P] &\longrightarrow &\ZZ[t^{\pm1}] \cr
e^{\lambda} &\longmapsto &t^{\langle -\rho,\lambda\rangle} \cr
}$$
to both sides of the identity (c) for the root system
$R^\vee = \{\alpha^\vee\ |\ \alpha\in R\}$ gives
$$W_0(t)
= \sum_{w\in W} \prod_{\alpha\in R^+} \left(
{1-t^{\langle \rho,w\alpha^\vee\rangle+1}\over 
1-t^{\langle \rho,w\alpha^\vee\rangle} } \right).
\eqno{(*)}$$
If $w\in W$, $w\ne 1$,  and $w=s_{i_1}\cdots s_{i_p}$ is a reduced
word for $w$ then $w^{-1}(-\alpha_{i_1})=
(s_{i_1}w)^{-1}\alpha_{i_1}\in R(w)$ and so 
$$\hbox{there is an $\alpha\in R^+$ such that 
$w\alpha^\vee = -\alpha_{i_1}^\vee$.}$$
Then
$$\eqalign{
\prod_{\alpha\in R^+} { 1-t^{\langle \rho,w\alpha^\vee\rangle+1}
\over 1-t^{\langle \rho,w\alpha^\vee\rangle} }
&=
{ 1-t^{\langle \rho,-\alpha_{i_1}^\vee\rangle+1}
\over 1-t^{\langle \rho,-\alpha_{i_1}^\vee\rangle} }
\prod_{\alpha\in R^+\atop w\alpha\ne -\alpha_{i_1}} 
{1-t^{\langle \rho,w\alpha^\vee\rangle+1}
\over 1-t^{\langle \rho,w\alpha^\vee\rangle} } \cr
&=
{ 1-t^{-1+1}\over 1-t}
\prod_{\alpha\in R^+\atop w\alpha\ne -\alpha_{i_1}} 
{1-t^{\langle \rho,w\alpha^\vee\rangle+1}
\over 1-t^{\langle \rho,w\alpha^\vee\rangle} }
=0. \cr
}$$
Thus the only nonzero term on the right hand side of ($*$) occurs for
$w=1$.
\endpf

\lemma {\sl For $\lambda\in P^+$ let $t_\lambda\in \tilde W$ 
be the translation in $\lambda$ and let $n_\lambda$ be the 
maximal length element in the double coset $Wt_\lambda W$.  
Let $M_\lambda = {\bf 1}_0x^\lambda{\bf 1}_0,$
as in (2.4).  Then
$$
q^{-\ell(w_0)} W_0(q^2)
\cdot {W_0(q^{-2})\over W_\lambda(q^{-2})}\cdot
M_\lambda = 
\sum_{x\in Wt_\lambda W} q^{\ell(x)-\ell(n_\lambda)} T_x,$$
in the affine Hecke algebra $\tilde H$.
}
\pf
Let $\lambda\in P^+$.  
Let $W_\lambda={\rm Stab}(\lambda)$ and let
$w_0$ and $w_\lambda$ be the maximal length elements in $W$ and 
$W_\lambda$, respectively.  Let $m_\lambda$ (resp. $n_\lambda$) 
be the minimal (resp. maximal)
length element in the double coset $Wt_\lambda W$.
For each positive root $\alpha$
the hyperplanes $H_{\alpha,i}$, 
$1\le i\le \langle\lambda,\alpha^\vee\rangle$,
are between the fundamental alcove $A$ and the alcove $t_\lambda A$ and so
$$\ell(t_\lambda) 
= \sum_{\alpha\in R^+} \langle \lambda,\alpha^\vee\rangle
=2\langle\lambda,\rho^\vee\rangle,
\qquad\hbox{where}\quad 
\rho^\vee = \hbox{$1\over 2$}\sum_{\alpha\in R^+} \alpha^\vee,
\formula
$$
and since
$n_\lambda = t_{w_0\lambda}w_0$ and 
$m_\lambda = t_\lambda(w_\lambda w_0)$, 
$$\eqalign{
\ell(m_\lambda) 
&= \ell(t_\lambda)-\ell(w_0w_\lambda)
=\ell(t_\lambda)-(\ell(w_0)-\ell(w_\lambda)),
\qquad\hbox{and} \cr
\ell(n_\lambda) &= \ell(t_\lambda)+\ell(w_0) = \ell(m_\lambda)
+\ell(w_0)-\ell(w_\lambda)+\ell(w_0). \cr
}
\formula$$
For example, in the setting of Example 1.4,
if $\lambda = 2\omega_2$ in type $C_2$, then $W_{\lambda} = \{1, s_1\}$,
$w_\lambda = s_1$, $w_0=s_1s_2s_1s_2$, $\ell(t_\lambda) = 6$,
$\ell(m_\lambda) = 4$, and $\ell(n_\lambda)=10$.  Labeling the alcove
$wA$ by the element $w$ the double coset $Wt_\lambda W$ consists of the 
elements in the four shaded diamonds.
$$
\matrix{
\beginpicture
\setcoordinatesystem units <1.25cm,1.25cm>         
\setplotarea x from -4.5 to 4.5, y from -3.5 to 4    
\put{$H_{\alpha_1}$}[b] at 0 3.3
\put{$H_{\alpha_2}$}[bl] at 3.3 3.2
\put{$H_\varphi = H_{\alpha_1+\alpha_2}$}[br] at -3.3 3.3
\put{$H_{\alpha_1+\alpha_2,1}=H_{\varphi,0}=H_{\alpha_0}$}[tl] at 3.3 -2.3 
\put{$H_{\alpha_1+2\alpha_2}$}[l] at 3.3 0
\put{$\scriptstyle{t_\lambda}$} at 0.2 2.5
\put{$\scriptstyle{m_\lambda}$} at 0.2 1.5
\put{$\scriptstyle{n_\lambda}$} at -0.2 -2.5
\put{$\scriptstyle{w_0}$} at -0.2 -0.5
\put{$\scriptstyle{w_\lambda}$} at -0.2 0.5
\plot -3.2 -3.2   3.2 3.2 /
\plot  3.2 -3.2  -3.2 3.2 /
\plot  0  3.2   0 -3.2 /
\plot  3.2  0  -3.2  0 /
\plot -1 2  0 3 /
\plot -1 2  0 1 /
\plot  0 1  1 2 /
\plot  0 3  1 2 /
\plot -1 -2  0 -3 /
\plot -1 -2  0 -1 /
\plot  0 -1  1 -2 /
\plot  0 -3  1 -2 /
\plot  1 0  2 1 /
\plot  2 1  3 0 /
\plot  1 0  2 -1 /
\plot  2 -1  3 0 /
\plot  -1 0  -2 1 /
\plot  -2 1  -3 0 /
\plot  -1 0  -2 -1 /
\plot  -2 -1  -3 0 /
\setdashes
\plot  -3 3.2   -3 -3.2 /
\plot  -2 3.2   -2 -3.2 /
\plot  -1 3.2   -1 -3.2 /
\plot   1 3.2    1 -3.2 /
\plot   2 3.2    2 -3.2 /
\plot   3 3.2    3 -3.2 /
\plot  3.2 -3   -3.2 -3 /
\plot  3.2 -2   -3.2 -2 /
\plot  3.2 -1   -3.2 -1 /
\plot  3.2  1   -3.2  1 /
\plot  3.2  2   -3.2  2 /
\plot  3.2  3   -3.2  3 /
\plot  3.2 -1.8   1.8 -3.2 /
\plot  3.2 -0.8   0.8 -3.2 /
\plot  3.2 0.2    -0.2 -3.2 /
\plot  3.2 1.2   -1.2 -3.2 /
\plot  3.2 2.2   -2.2 -3.2 /
\plot  2.2 3.2   -3.2 -2.2 /
\plot  1.2 3.2   -3.2 -1.2 /
\plot  0.2 3.2   -3.2 -0.2 /
\plot  -0.8 3.2  -3.2 0.8 /
\plot  -1.8 3.2  -3.2 1.8 /
\plot  -2.2 3.2   3.3 -2.3 /
\plot  -1.2 3.2   3.2 -1.2 /
\plot  -0.2 3.2   3.2 -0.2 /
\plot  0.8 3.2  3.2 0.8 /
\plot  1.8 3.2  3.2 1.8 /
\plot  -3.2 -1.8   -1.8 -3.2 /
\plot  -3.2 -0.8   -0.8 -3.2 /
\plot  -3.2 0.2    0.2 -3.2 /
\plot  -3.2 1.2   1.2 -3.2 /
\plot  -3.2 2.2   2.2 -3.2 /
\vshade 1 0 0   2 -1 1 /
\vshade 2 -1 1  3 0 0 /
\vshade -2 -1 1   -1 0 0 /
\vshade -3 0 0  -2 -1 1 /
\vshade -1 2 2   0 1 3 /
\vshade  0 1 3   1 2 2 /
\vshade -1 -2 -2   0 -3 -1 /
\vshade  0 -3 -1   1 -2 -2 /
\endpicture
\cr
\cr
\hbox{The double coset $Wt_\lambda W$}\qquad\qquad \cr
}
$$

Then
$$\eqalign{
{\bf 1}_0x^\lambda {\bf 1}_0 
&= 
{\bf 1}_0 T_{t_\lambda} {\bf 1}_0
= 
{\bf 1}_0 T_{m_\lambda w_0 w_\lambda} {\bf 1}_0
= 
{\bf 1}_0 T_{m_\lambda} T_{w_0 w_\lambda} {\bf 1}_0
= 
q^{\ell(w_0)-\ell(w_\lambda)}
{\bf 1}_0 T_{m_\lambda} {\bf 1}_0
\cr
&= { q^{\ell(w_0)-\ell(w_\lambda)-\ell(m_\lambda)}\over W(q^2)}
\left(\sum_{w\in W} q^{\ell(w)}T_w\right)q^{\ell(m_\lambda)}
T_{m_\lambda}{\bf 1}_0.\cr
}$$
Let $W^\lambda$ be a set of minimal length coset representatives
of the cosets in $W/W_\lambda$.  Every element $w\in W$ has a unique
expression $w=uv$ with $u\in W^\lambda$ and $v\in W_\lambda$.
If $v\in W_\lambda$ then
$$
vm_\lambda = vt_\lambda w_\lambda w_0 =t_\lambda vw_\lambda w_0 
= m_\lambda (w_\lambda w_0)^{-1} v w_\lambda w_0
=m_\lambda (w_0^{-1} w_\lambda^{-1} v w_\lambda w_0). 
$$
Since conjugation by $w_\lambda$ (resp. conjugation by $w_0$) 
is an automorphism of $W_\lambda$ (resp. $W$) which takes simple reflections
to simple reflections, 
$\ell(v) = \ell(w_0^{-1}w_\lambda^{-1} v w_\lambda w_0)$.
Thus
$$\eqalign{
{\bf 1}_0 x^\lambda{\bf 1}_0
&= { q^{\ell(w_0)-\ell(w_\lambda)-\ell(m_\lambda)}\over W_0(q^2)}
\sum_{u\in W^\lambda} q^{\ell(u)}T_u
\sum_{v\in W_\lambda} q^{\ell(v)}T_v 
q^{\ell(m_\lambda)}T_{m_\lambda}{\bf 1}_0 \cr
&= { q^{2\ell(w_0)-2\ell(w_\lambda)-\ell(t_\lambda)}\over W_0(q^2)}
\left( \sum_{u\in W^\lambda} q^{\ell(u)}T_u 
      q^{\ell(m_\lambda)}T_{m_\lambda}\right)
\left(\sum_{v\in w_0^{-1}w_\lambda^{-1}
W_{\lambda}w_\lambda w_0} q^{\ell(v)}T_v\right){\bf 1}_0 \cr
&= { q^{-2\ell(w_\lambda)-\ell(t_\lambda)}\over W_0(q^{-2})}
\left(\sum_{u\in W^\lambda} q^{\ell(u)}T_u \right)
      q^{\ell(m_\lambda)}T_{m_\lambda}
W_\lambda(q^2){\bf 1}_0\cr
&= { q^{-2\ell(w_\lambda)-\ell(t_\lambda)}
W_\lambda(q^2)
\over W_0(q^2)W_0(q^{-2})}
\left(\sum_{u\in W^\lambda} q^{\ell(u)}T_u \right)
      q^{\ell(m_\lambda)}T_{m_\lambda}
\left(\sum_{w\in W} q^{\ell(w)}T_w\right)\cr
&= { q^{-\ell(t_\lambda)}
W_\lambda(q^{-2}) \over W_0(q^2)W_0(q^{-2})}
\sum_{x\in Wt_\lambda W} q^{\ell(x)} T_x \cr
&= {q^{-\ell(t_\lambda)+\ell(n_\lambda)}
W_\lambda(q^{-2})\over W_0(q^2)W_0(q^{-2})}
\left(\sum_{x\in Wt_\lambda W} q^{\ell(x)-\ell(n_\lambda)}T_x\right) \cr
&= {q^{\ell(w_0)}\over W_0(q^2)}{W_\lambda(q^{-2})\over W_0(q^{-2})}
\left(\sum_{x\in Wt_\lambda W} q^{\ell(x)-\ell(n_\lambda)}T_x\right).
\qquad\hbox{\qed}\cr
}$$

\lemma  Let $w_0$ be the longest element of $W$ and let $\lambda\in P$.
\smallskip\noindent
\itemitem{(a)} $\overline{x^\lambda} = T_{w_0}x^{w_0\lambda}T_{w_0}^{-1}$.
\smallskip\noindent
\itemitem{(b)} $\overline{{\bf 1}_0} = {\bf 1}_0$ and 
$\overline{\varepsilon_0}=\varepsilon_0$.
\smallskip\noindent
\itemitem{(c)} If $z\in \ZZ[P]^W$ then $\bar z =z$. 
\smallskip\noindent
\itemitem{(d)} $\overline{q^{-\ell(w_0)}A_{\lambda+\rho}} 
= q^{-\ell(w_0)}A_{\lambda+\rho}$.
\pf
(a)
If $\lambda\in P^+$ then
$w_0t_\lambda = t_{w_0\lambda}w_0$,
$\ell(w_0t_\lambda) = \ell(w_0)+\ell(t_\lambda)$ and
$\ell(t_{w_0\lambda})w_0) = \ell(t_{w_0\lambda}) + \ell(w_0)$.
Thus, 
$$T_{w_0}T_{t_\lambda}=T_{w_0t_\lambda}=T_{t_{w_0\lambda}w_0}
=T_{t_{w_0\lambda}}T_{w_0},
\qquad\hbox{for $\lambda\in P^+$}.$$
Let $\lambda\in P$ and write $\lambda=\mu-\nu$ with $\mu,\nu\in P^+$.
Since $-w_0\mu\in P^+$ and $-w_0\nu\in P^+$,
$$\overline{x^\lambda} 
= \overline{ T_{t_\mu} T_{t_{\nu}}^{-1} }
=T_{t_{-\mu}}^{-1} T_{t_{-\nu}}
=T_{w_0}T_{t_{-w_0\mu}}^{-1}T_{t_{-w_0\nu}}T_{w_0}^{-1}
=T_{w_0}(x^{-w_0\lambda})^{-1}T_{w_0}^{-1},
=T_{w_0}x^{w_0\lambda}T_{w_0}^{-1}.$$
\smallskip\noindent
(b)  For $1\le i\le n$,
$$\matrix{
\overline{{\bf 1}_0^2} = \overline{{\bf 1}_0}^2
&\qquad &\hbox{and} &\qquad 
&T_i \overline{{\bf 1}_0} = \overline{T_i^{-1}{\bf 1}_0}
=\overline{q^{-1}{\bf 1}_0} = q\overline{{\bf 1}_0}, \hfill \cr
\cr
\overline{\varepsilon_0^2} = \overline{\varepsilon_0}^2
&\qquad &\hbox{and} &\qquad 
&T_i\overline{\varepsilon_0} = \overline{T_i^{-1}\varepsilon_0}
=\overline{-q\varepsilon_0} = -q^{-1}\overline{\varepsilon_0}. \hfill \cr
}$$
These are the defining properties (2.1) of ${\bf 1}_0$ and 
$\varepsilon_0$ and so $\overline{{\bf 1}_0}={\bf 1}_0$ and 
$\overline{\varepsilon_0} = \varepsilon_0$.
\smallskip\noindent
(c)  If $z = \sum_{\mu\in P} c_\mu x^\mu \in \ZZ[P]^W$, then,
since $c_\mu\in \ZZ$, $\overline{c_\mu}=c_\mu$ and, by (a), 
$$
\overline{z} = \sum_{\mu\in P} \overline{c_\mu}\overline{x^\mu}
=\sum_{\mu\in P} c_\mu T_{w_0}x^{w_0\mu} T_{w_0}^{-1}
=T_{w_0}\left( \sum_{\mu\in P} c_\mu x^{w_0\mu}\right) T_{w_0}^{-1}
=T_{w_0} z T_{w_0}^{-1},$$
since $z\in \ZZ[P]^W$ is $W$-invariant.  Finally, since
$\ZZ[P]^W\subseteq Z(\tilde H)$, $z$ is central,
and $\overline{z} = T_{w_0}zT_{w_0}^{-1}=z$.
\smallskip\noindent
(d)  By (a), (b) and the third equality in (2.6),
$$\eqalign{
\overline{ q^{-\ell(w_0)} A_{\lambda+\rho} }
&=q^{\ell(w_0)}\overline{\varepsilon_0 x^{\lambda+\rho} {\bf 1}_0 } 
=q^{\ell(w_0)}
\varepsilon_0 T_{w_0} x^{w_0(\lambda+\rho)} T_{w_0}^{-1} {\bf 1}_0 \cr
&=q^{\ell(w_0)} (-q^{-1})^{\ell(w_0)} \varepsilon_0 x^{w_0(\lambda+\rho)}
{\bf 1}_0 q^{-\ell(w_0)}
=(-q^{-1})^{\ell(w_0)} A_{w_0(\lambda+\rho)}  \cr
&= q^{-\ell(w_0)} A_{\lambda+\rho}. 
\qquad\hbox{\qed}\cr
}$$

The following theorem is due to Lusztig [Lu].  Part (a) was originally
proved in a different, but equivalent, formulation
by Macdonald [Mac2, (4.1.2)].

\thm {\sl  If $\mu\in P$ let $W_\mu$ be the stabilizer of
$\mu$ and let $W_\mu(t)$ be as in (2.12). 
\smallskip\noindent
\item{(a)} Let $\mu\in P$. 
Let $P_\mu(x;t)$ be the Macdonald spherical function defined in (2.13)
and define $M_\mu={\bf 1}_0x^\mu{\bf 1}_0$ as in (2.4).
In the affine Hecke algebra $\tilde H$,
$$
{W_\mu(q^{-2})\over W_0(q^{-2})}\cdot
P_\mu(x;q^{-2}){\bf 1}_0 = M_\mu.  
$$
\item{(b)}
For $\lambda\in P^+$ let $t_\lambda\in \tilde W$ 
be the translation in $\lambda$ and let $n_\lambda$ be the 
maximal length element in the double coset $Wt_\lambda W$.  
Let $s_\lambda$ be the Weyl character and let $C_{n_\lambda}'$
be the Kazhdan-Lusztig basis element as defined in (2.8) and (1.28),
respectively.  In the affine Hecke algebra $\tilde H$,
$$q^{-\ell(w_0)}W_0(q^2)\cdot s_\lambda{\bf 1}_0 = C_{n_\lambda}'.$$
}
\pf
(a) By Theorem 2.11 there is an element $\tilde P_\lambda\in \KK[P]^W$ 
such that $\tilde P_\lambda {\bf 1}_0 = {\bf 1}_0 x^\lambda {\bf 1}_0$.
To find $\tilde P_\lambda$ first do a rank 1 calculation,
$$\eqalign{
(q^{-1}+T_i)x^\lambda {\bf 1}_0
&=\left(q^{-1}x^\lambda+x^{s_i\lambda}T_i+(q-q^{-1})\Big(
{x^\lambda-x^{s_i\lambda} \over 1-x^{-\alpha_i}}\Big)\right){\bf 1}_0 \cr
&={1\over 1-x^{-\alpha_i} }
\left(\eqalign{
q^{-1}x^\lambda(1-x^{-\alpha_i})
+qx^{s_i\lambda}(1-x^{-\alpha_i} \cr
+qx^\lambda-qx^{s_i\lambda}-q^{-1}x^\lambda+q^{-1}x^{s_i\lambda} \cr
}\right){\bf 1}_0 \cr
&=(1-x^{-\alpha_i})^{-1}
(-q^{-1}x^{\lambda-\alpha_i}-qx^{s_i\lambda-\alpha_i}
+qx^\lambda+q^{-1}x^{s_i\lambda})
{\bf 1}_0 \cr
&=(1-x^{-\alpha_i})^{-1}(
x^\lambda(q-q^{-1}x^{-\alpha_i})
+x^{s_i\lambda}(q^{-1}-qx^{-\alpha_i})){\bf 1}_0 \cr
&=\left(
{q-q^{-1}x^{-\alpha_i}\over 1-x^{-\alpha_i}}\cdot x^\lambda
+
{x^{-\alpha_i}\over x^{-\alpha_i}}\cdot
{q^{-1}x^{\alpha_i}-q\over x^{\alpha_i}-1}\cdot
x^{s_i\lambda}
\right)
{\bf 1}_0 \cr
&=(1+s_i)\left(
{q-q^{-1}x^{-\alpha_i}\over 1-x^{-\alpha_i}}\,x^\lambda\right){\bf 1}_0. \cr
}$$
Since ${\bf 1}_0$ is a linear combination of products of $T_i$
it can also be written as a linear combination of
products of $q^{-1}+T_i$.  Thus
${\bf 1}_0 x^\lambda {\bf 1}_0$ 
can be written as a linear combination of terms of the form
$$(1+s_{i_1})
\left({q-q^{-1}x^{-\alpha_{i_1}}\over 1-x^{-\alpha_{i_1}}}\right)
\cdots (1+s_{i_p})
\left({q-q^{-1}x^{-\alpha_{i_p}}\over 1-x^{-\alpha_{i_p}}}\right)
x^\lambda.$$
Thus
$${\bf 1}_0x^\lambda {\bf 1}_0 = \tilde P_\lambda {\bf 1}_0,
\qquad\hbox{where}\qquad
\tilde P_\lambda = \sum_{w\in W} x^{w\lambda} wc_w,
$$
and the $c_w$ are some linear combinations of products of
terms of the form $(q-q^{-1}x^{\alpha})/(1-x^\alpha)$ for 
roots $\alpha\in R$.
Since $\tilde P_\lambda$ is an element of $\KK[P]^W$,
$$\tilde P_\lambda = \sum_{w\in W} w(x^{w_0\lambda} w_0 c_{w_0}),$$
where $w_0$ is the longest element of $W$.  
The coefficient $w_0c_{w_0}$ comes from the highest term in the expansion 
of 
$${\bf 1}_0 = {1\over W_0(q^2)}(q^{2\ell(w_0)}T_{w_0} + \hbox{lower terms})
$$ 
in terms of linear combination of products of the $(q^{-1}+T_i)$.
If $w_0=s_{i_1}\cdots s_{i_p}$ is a reduced word for $w_0$ then
$$\eqalign{
w_0c_{w_0}&={q^{\ell(w_0)}\over W_0(q^2)}
s_{i_1}
\left({q-q^{-1}x^{-\alpha_{i_1}}\over 1-x^{-\alpha_{i_1}}}\right)
\cdots
s_{i_p}
\left({q-q^{-1}x^{-\alpha_{i_p}}\over 1-x^{-\alpha_{i_p}}}\right) \cr
&=
s_{i_1}\cdots s_{i_p}
\left({q-q^{-1}x^{-s_{i_p}\cdots s_{i_2}\alpha_{i_1}}\over 1-x^{-s_{i_p}\cdots s_{i_2}\alpha_{i_1}}}
\right)
\left({q-q^{-1}x^{-s_{i_p}\cdots s_{i_3}\alpha_{i_2}}\over 1-x^{-s_{i_p}\cdots s_{i_2}\alpha_{i_2}}}
\right)
\cdots
\left({q-q^{-1}x^{-\alpha_{i_p}}\over 1-x^{-\alpha_{i_p}}}\right) \cr
&={q^{\ell(w_0)}\over W_0(q^2)}w_0
\prod_{\alpha\in R^+} {q-q^{-1}x^{-\alpha}\over 1-x^{-\alpha} }
={q^{2\ell(w_0)}\over W_0(q^2)}w_0
\prod_{\alpha\in R^+} {1-q^{-2}x^{-\alpha}\over 1-x^{-\alpha} }, \cr
}$$
by Lemma 1.11 and the fact that $\ell(w_0)={\rm Card}(R^+)$.  Thus,
since $q^{-2\ell(w_0)}W_0(q^2)=W_0(q^{-2})$,
$$\tilde P_\lambda = {1\over W_0(q^{-2})}
\sum_{w\in W} w\left(
x^\lambda\prod_{\alpha\in R^+} 
{1-q^{-2}x^{-\alpha}\over 1-x^{-\alpha} }\right).
$$

(b)
Since $W_0(q^{-2})=q^{-2\ell(w_0)}W_0(q^2)$, Lemma 2.21 gives
$$\overline{q^{-\ell(w_0)}W_0(q^2)s_\lambda {\bf 1}_0 }
= q^{\ell(w_0)} W_0(q^{-2}) \overline{s_\lambda} {\bf 1}_0
= q^{-\ell(w_0)}W_0(q^2)s_\lambda{\bf 1}_0.$$
By Lemma 2.16(b),
$$s_\lambda = 
\sum_{\mu\in P^+} K_{\lambda\mu}(t) P_\mu(x;t),$$
where $K_{\lambda\mu}(t)\in \ZZ[t]$,
$K_{\lambda\mu}(t)=0$ unless $\mu\le \lambda$ and 
$K_{\lambda\lambda}(t)=1$.
Thus, by part (a) and Lemma 2.18
$$\eqalign{
q^{-\ell(w_0)}W_0(q^2)s_\lambda{\bf 1}_0
&=
\sum_{\mu\in P^+} q^{-\ell(w_0)}W_0(q^2)
K_{\lambda\mu}(q^{-2})P_\mu(x;q^{-2}) {\bf 1}_0 \cr
&= \sum_{\mu\in P^+} \sum_{x\in Wt_\mu W}
q^{\ell(x)-\ell(n_\mu)}K_{\lambda\mu}(q^{-2})T_x, \cr
}$$
where the polynomials $K_{\lambda\mu}(q^{-2})\in \ZZ[q^{-2}]$
are 0 unless $\mu\le \lambda$ and $K_{\lambda\lambda}(q^{-2})=1$.
Hence $q^{-\ell(w_0)}W(q^2)s_\lambda {\bf 1}_0$ is a
bar invariant element of $\tilde H$ such that its expansion 
in terms of the basis
$\{ T_w\ |\ w\in \tilde W\}$ is triangular with coefficient
of $T_{n_\lambda}=1$ and all other coefficients in $q^{-1}\ZZ[q^{-1}]$.
These are the defining properties (1.28-9) of $C_{n_\lambda}'$.
\endpf

\vfill\eject

\section 3. Orthogonality and formulas for Kostka-Foulkes poylnomials

\medskip
Let $\KK = \ZZ[t]$.
If $f = \sum_{\mu\in P} f_\mu x^\mu\in \KK[P]$ let
$$\bar f = \sum_{\mu\in P} f_\mu x^{-\mu},
\qquad\hbox{and}\qquad
[f]_1 = f_0 = (\hbox{coefficient of $1$ in $f$}).
\formula$$
Define a symmetric bilinear form
$$\langle, \rangle_t\colon \KK[P]\times \KK[P] \to \KK
\qquad\hbox{by}\qquad
\langle f,g\rangle_t 
= {1\over |W|} \left[f\bar g
\prod_{\alpha\in R} {1-x^\alpha\over 1-tx^\alpha}
\right]_1.
\formula$$
``Specializing'' $t$ at the values $0$ and $1$ gives inner products
$\langle, \rangle_0\colon \KK[P]\times \KK[P] \to \KK$
and 
$\langle, \rangle_0\colon \KK[P]\times \KK[P] \to \KK$ with
$$\langle f,g\rangle_0 
= {1\over |W|} 
\left[f\bar g \prod_{\alpha\in R} 1-x^\alpha \right]_1
\qquad\hbox{and}\qquad
\langle f,g\rangle_1 = {1\over |W|} [f\bar g]_1.
\formula$$

\prop {\sl  Let $\lambda$ and $\mu\in P^+$.
Then
$$\langle m_\lambda, m_\mu\rangle_1 
= {1\over |W_\lambda|}\delta_{\lambda\mu},
\qquad
\langle s_\lambda, s_\mu\rangle_0 
= \delta_{\lambda\mu},
\qquad\hbox{and}\qquad
\langle P_\lambda, P_\mu\rangle_t 
= {1\over W_\lambda(t)}\delta_{\lambda\mu}.
$$
}
\pf
The first equality follows from
$$|W_\lambda|\langle m_\lambda,m_\mu\rangle_1 = 
{|W_\lambda|\over |W|} 
\sum_{\gamma\in W\lambda\atop \nu\in W\mu}
[x^\gamma x^{-\nu}]_1
=\delta_{\lambda\mu}{|W_\lambda|\over |W|}
\sum_{\gamma\in W\lambda} 1 = \delta_{\lambda\mu}.$$
If $\lambda,\mu\in P^+$, 
$$\eqalign{
\langle s_\lambda,s_\mu\rangle_0
&= {1\over |W|} [\overline{a_\rho s_\lambda}a_\rho s_\mu]_1
={1\over |W|} [\overline{a_{\lambda+\rho}} a_{\mu+\rho}]_1 \cr
&={1\over |W|} \sum_{v,w\in W}
(-1)^{\ell(v)}(-1)^{\ell(w)} 
[x^{v(\lambda+\rho)}x^{-w(\mu+\rho)}]_1 \cr
&=\delta_{\lambda\mu} {1\over |W|} \sum_{v\in W}
(-1)^{\ell(v)}(-1)^{\ell(v)} = \delta_{\lambda\mu}, \cr
}$$
giving the second statement.  

\smallskip
By Lemma 2.16(b) the matrix $K^{-1}$ given by the values
$(K^{-1})_{\lambda\mu}$ in the equation
$$P_\lambda(x;t) = \sum_{\mu} (K^{-1})_{\lambda\mu}s_\mu,$$
has entries in $\ZZ[t]$ and is upper triangular with 1's on
the diagonal, i.e.  $(K^{-1})_{\lambda\lambda}=1$ 
and $(K^{-1})_{\lambda\mu}=0$ unless $\mu\le \lambda$.  
Since $P_\lambda(x;1)=m_\lambda$ the matrix $k^{-1}$ 
describing the change of basis
$$m_\lambda = \sum_{\mu} (k^{-1})_{\lambda\mu}s_\mu,$$
is the specialization of $K^{-1}$ at $t=1$ and so
$k^{-1}$ has entries in $\ZZ$ and
is upper triangular with 1's on the diagonal.
Hence the matrix $A = K^{-1}k^{-1}$ giving the change of basis
$$P_\lambda(x;t) = \sum_{\nu\le \lambda} A_{\lambda\nu}m_\mu,
\formula$$
has $A_{\lambda\mu}\in \ZZ[t]$, $A_{\lambda\lambda}=1$, 
and $A_{\lambda\mu}=0$ unless $\mu\le \lambda$.  

Let $Q^+$ be the set of nonnegative integral linear combinations
of positive roots.
$$\eqalign{
P_\mu(x;t) W_\mu(t) 
\Big(\prod_{\alpha\in R} {1-x^\alpha\over 1-tx^\alpha}\Big)
&= \sum_{w\in W} w\left( x^\mu \prod_{\alpha\in R^+}
{1-x^\alpha\over 1-tx^\alpha}\right)  \cr
&= \sum_{w\in W} w\left(
x^\mu\prod_{\alpha\in R^+} (1+\sum_{r>0}t^{r-1}(t-1)x^{r\alpha})
\right)
\cr
&= \sum_{w\in W} w\left(\sum_{\nu\in Q^+} c_\nu x^{\mu+\nu}\right) 
= \sum_{\nu\in Q^+} c_\nu\left(
\sum_{w\in W} wx^{\mu+\nu} \right), \cr
}$$
where $c_\nu\in \ZZ[t]$ and $c_0=1$.  Hence
$$
P_\mu(x;t) W_\mu(t)\prod_{\alpha\in R}
{1-x^\alpha\over 1-tx^\alpha}
= |W_\mu|m_\mu + \sum_{\gamma> \mu} B_{\mu\gamma}m_\gamma
= \sum_{\gamma\ge \mu} B_{\mu\gamma}m_\gamma,
\eqno{(*)}$$
with $B_{\mu\gamma}\in \ZZ[t]$ and $B_{\mu\mu}=|W_\mu|$.

Assume that $\lambda\le \mu$ if $\lambda$ and $\mu$ are 
comparable.  Then, by using (3.5) and ($*$),
$$
\langle P_\lambda, P_\mu\rangle_t
={1\over W_\mu(t)}
\left\langle P_\lambda, P_\mu W_\mu(t)
\prod_{\alpha\in R} {1-x^\alpha\over 1-tx^\alpha}\right\rangle_1 
={1\over W_\mu(t)} \left\langle
\sum_{\nu\le \lambda} A_{\lambda\nu} m_\nu, 
\sum_{\gamma\ge \mu} B_{\mu\gamma} m_\gamma\right\rangle_1.
$$
Since $A_{\lambda\lambda}=1$ and $B_{\mu\mu}=|W_\mu|$ the result
follows from $\langle m_\lambda, m_\mu\rangle_1 
= |W_\lambda|^{-1}\delta_{\lambda\mu}$.
\endpf

The following theorem shows that the spherical functions 
$P_\lambda(x,t)$ are uniquely determined by the 
triangularity in (3.5) and the orthogonality in the third
equality of Proposition 3.4.

\thm {\sl  Let $\KK = \ZZ[t]$.
The spherical functions $P_\lambda(x;t)$ are the unique 
elements of $\KK[P]^W$ such that 
\smallskip\noindent
\item{(a)} 
$\displaystyle{
P_\lambda = m_\lambda + \sum_{\mu<\lambda} A_{\lambda\mu} m_\mu},$
\smallskip\noindent
\item{(b)} 
$\displaystyle{
\langle P_\lambda, P_\mu\rangle_t=0}$ if $\lambda\ne \mu$.
}
\pf
Assume that the $P_\mu$ are determined for $\mu<\lambda$.  Then
the condition in (a) can be rewritten as 
$$P_\lambda = m_\lambda + \sum_{\mu<\lambda} C_{\lambda\mu} P_\mu,$$
for some constants $C_{\lambda\mu}$.  Take the inner
product on each side with $P_\nu$, $\nu<\lambda$, 
and use property (b) to get the system of equations
$$
0=\langle m_\lambda,P_\nu\rangle_t + \sum_{\mu<\lambda}
C_{\lambda\mu}\langle P_\mu,P_\nu\rangle_t 
= \langle m_\lambda, P_\nu\rangle_t + 
C_{\lambda\nu}\langle P_\nu,P_\nu\rangle_t. 
$$
Hence
$$
C_{\lambda\nu} = {-\langle m_\lambda, P_\nu\rangle_t \over
\langle P_\nu, P_\nu\rangle_t}\,,
\qquad
\hbox{for each $\nu<\lambda$,}$$
and this determines $P_\lambda$.
\endpf

\medskip\noindent
\remark
(a) The inner product $\langle,\rangle_t$ arises naturally
in the context of $p$-adic groups.   Let $S^1 = \{ z\in \CC\ |\ |z|=1\}$
and view the $x^\lambda$, $\lambda\in P$, 
as characters of 
$$T = \Hom(P,S^1) \qquad\hbox{via}\qquad
\matrix{x^\lambda\colon &T &\longrightarrow &\CC^* \cr
&s &\longmapsto &s(\lambda). \cr}
\formula$$
Let $ds$ be the Haar measure on $T$ normalized so that
$$\langle x^\lambda, x^\mu\rangle = 
\int_T x^\lambda(s)\overline{x^\mu(s)} ds = \delta_{\lambda\mu}.
\formula$$
Letting $\QQ_p$ be the field of $p$-adic numbers,
Macdonald [Mac2, (5.1.2)] showed that the Plancherel measure for the
$p$-adic Chevalley group $G(\QQ_p)$ corresponding to the root system
$R$ is given by
$$d\mu(s) = {W_0(p^{-1})\over |W|} \prod_{\alpha\in R} {1-x^\alpha(s)
\over 1-p^{-1}x^\alpha(s)}.\formula$$
The corresponding inner product is 
$$W_0(p^{-1})\langle f,g\rangle_{p^{-1}} 
= \int_T f(s)\overline{g(s)} d\mu(s),
\qquad\hbox{for $f,g\in C(T)$,}$$
where $C(T)$ is the vector space of continuous functions on $T$.

\smallskip\noindent
(b)  The inner product $\langle,\rangle_t$ arises naturally
in another representation theoretic context.  The
complex semisimple Lie algebra $\fg$ corresponding to the root system
$R$ acts on $S(\fg^*)$, the ring of polynomials on $\fg$,
by the (co-)adjoint action and 
as graded $\fg$-modules the characters of 
$S(\fg^*)$ and the subring of invariants $S(\fg^*)^{\fg}$ are 
$$\eqalign{
{\rm grch}(S(\fg^*))
&=\left(\prod_{i=1}^r {1\over 1-t}\right)
\left(\prod_{\alpha\in R} {1\over 1-tx^\alpha}\right)
\qquad\hbox{and}\cr
{\rm grch}(S(\fg^*)^{\fg})
&=\prod_{i=1}^r{1\over 1-t^{d_i}}
=\left(\prod_{i=1}^r{1-t\over 1-t^{d_i}}\right)
\left(\prod_{i=1}^r {1\over 1-t}\right)
= {1\over W_0(t)} \prod_{i=1}^r {1\over 1-t}, \cr}
\formula$$
where $r$ is the rank of $\fg$ and $d_1,\ldots, d_r$ are the
{\it degrees} of the Weyl group $W$.  Let ${\cal H}$ denote
the vector space of harmonic polynomials.  An important theorem of 
Kostant [Ks, Theorem 0.2] gives that 
$$S(\fg^*) \cong S(\fg^*)^\fg\otimes {\cal H},
\qquad\hbox{and thus,}\qquad
{\rm grch}({\cal H}) 
= W_0(t)\prod_{\alpha\in R} {1\over 1-tx^\alpha}.
\formula$$
If $L(\lambda)$ denotes the finite dimensional irreducible
$\fg$-module of highest weight $\lambda\in P^+$ then
$L(\lambda)$ has character $s_\lambda$ and using the 
notation of (3.2),
$$\eqalign{
\sum_{k\ge 0} 
\dim(\Hom_{\fg}(L(\lambda), L(\mu)\otimes {\cal H}^k)t^k
&=\left\langle s_\lambda,
s_\mu W_0(t)\prod_{\alpha\in R} {1\over 1-tx^\alpha}
\right\rangle_0 \cr
&=W_0(t) \left[s_\lambda\overline{s_\mu}\prod_{\alpha\in R}
{1-x^\alpha\over 1-tx^\alpha}\right]_1
= W_0(t)\langle s_\lambda,s_\mu\rangle_t, \cr}
\formula$$
where ${\cal H}^k$ is the vector space of degree $k$ harmonic
polynomials.

\bigskip
\bigskip\noindent
{\it Formulas for Kostka-Foulkes polynomials}

\bigskip
For $\lambda\in P$ let $s_\lambda$ denote the Weyl character,
as defined in (2.8).  The {\it Kostka-Foulkes polynomials},
or {\it $q$-weight multiplicities}, 
$K_{\lambda\mu}(t)$, $\lambda,\mu\in P^+$, are defined by
the change of basis formula
$$s_\lambda = \sum_{\mu\in P^+} K_{\lambda\mu}(t) P_\mu(x;t),
\formula$$
where the Macdonald spherical functions $P_\mu(x;t)$ are
as in (2.13).

For each $\alpha\in R^+$ define the {\it raising operator}
$R_\alpha\colon P\to P$ by 
$$R_\alpha\lambda = \lambda+\alpha,
\qquad\hbox{and define}\qquad
(R_{\beta_1}\cdots R_{\beta_\ell})s_\lambda 
= s_{R_{\beta_1}\cdots R_{\beta_\ell}\lambda},
\formula$$
for any sequence $\beta_1,\ldots, \beta_\ell$ of positive roots.
Using the straightening law for Weyl characters (2.9),
$$s_\mu = (-1)^{\ell(w)} s_{w\circ\mu},\qquad
\hbox{where}\qquad w\circ\mu = w(\mu+\rho)-\rho,$$
any $s_\mu$ is equal to 0 or to $\pm s_\lambda$ with $\lambda\in P^+$.  
Composing the action of raising operators on Weyl characters
should be avoided.
For example, if $\alpha_i$ is a simple root then
(since $\langle \rho,\alpha_i^\vee\rangle=1$)
$s_{-\alpha_i} = - s_{s_i\circ(-\alpha_i)}
=-s_{s_i(\rho-\alpha_i)-\rho} = -s_{-\alpha_i}$
giving that $s_{-\alpha_i}=0$ and so
$$R_{\alpha_i}(R_{\alpha_i}s_{-2\alpha_i})=
R_{\alpha_i}s_{-\alpha_i}= R_{\alpha_i}\cdot 0 = 0,
\qquad\hbox{whereas}\qquad
(R_{\alpha_i}R_{\alpha_i})s_{-2\alpha_i} = s_0 =1.
$$

Let $Q^+$ be the set of nonnegative integral linear combinations
of positive roots.  Define the {\it $q$-analogue
of the partition function} $F(\gamma;t)$, $\gamma\in P$, by
$$\prod_{\alpha\in R^+} {1\over 1-tx^\alpha}
= \sum_{\gamma\in Q^+} F(\gamma;t)x^\gamma,
\qquad\hbox{and}\qquad
\hbox{$F(\gamma;t)=0$, \quad if $\gamma\not\in Q^+$.}
\formula$$

\thm {\sl Let $\lambda,\mu\in P^+$.  Let $t_\mu$ be the translation
in $\mu$ as defined in (1.12) and let $n_\mu$ be the 
longest element of the double coset $Wt_\mu W$.  Let
$W_\mu(t)$ be as in (2.12), $P_\mu(x;t)$ as in (2.13) and 
let $\langle,\rangle_t$ be the inner product defined in (3.2).
For $y,w\in\tilde W$
let $P_{yw}\in \ZZ[t^{\pm {1\over2}}]$ denote the 
Kazhdan-Lusztig polynomial defined in (1.28-9) and let
$\rho^\vee = {1\over 2}\sum_{\alpha\in R^+} \alpha^\vee$.
\medskip\noindent
\item{(a)} $K_{\lambda,\mu}(t) = 
W_\mu(t)\,\langle s_\lambda, P_\mu(x;t)\rangle_t$.
\smallskip\noindent
\item{(b)} $K_{\lambda\mu}(t)
= \hbox{coefficient of $s_\lambda$ in}\ \ 
\displaystyle{
\left(\prod_{\alpha\in R^+}{1\over 1-tR_\alpha}\right) s_\mu.} $
\smallskip\noindent
\item{(c)} $\displaystyle{
K_{\lambda\mu}(t) = 
\sum_{w\in W} (-1)^{\ell(w)} F(w(\lambda+\rho)-(\mu+\rho);t) }$.
\smallskip\noindent
\item{(d)} $K_{\lambda\mu}(t) 
= t^{\langle \lambda-\mu,\rho^\vee\rangle}P_{x,n_\lambda}(t^{-1})$,
for any $x\in Wt_\mu W$.
}
\pf
(a) This follows from the third equality in Proposition 3.4 and
the definition of $K_{\lambda\mu}(t)$.

\smallskip\noindent
(b)
Since
$$\eqalign{
P_\mu(x;t) W_\mu(t)\prod_{\alpha\in R} {1\over 1-tx^\alpha}
&= \sum_{w\in W} w\left( x^\mu
\prod_{\alpha\in R^+} {1-tx^{-\alpha}\over 1-x^{-\alpha}}\right)
\prod_{\alpha\in R} {1\over 1-tx^\alpha} \cr
&= \sum_{w\in W} w\left( x^{\mu+\rho}{1\over
x^\rho\prod_{\alpha\in R^+}(1-x^{-\alpha})(1-tx^{-\alpha}) }
\right) \cr
&= {1\over a_\rho}\sum_{w\in W} (-1)^{\ell(w)}
w\left( \prod_{\alpha\in R^+} \Big({1\over 1-tx^\alpha}\Big)
x^{\mu+\rho}\right)\,. \cr
}$$
Then
$$\eqalign{
K_{\lambda\mu}(t)
&= (\hbox{coefficient of $P_\mu(x;t)$ in $s_\lambda$})
= \langle s_\lambda, W_\mu(t)P_\mu(x;t) \rangle_t \cr
& = \left\langle s_\lambda,
W_\mu(t)P_\mu(x;t) \prod_{\alpha\in R}{1\over 1-tx^\alpha}
\right\rangle_0 \cr
&= \hbox{coefficient of $s_\lambda$ in } \ \ 
{1\over a_\rho}\sum_{w\in W} (-1)^{\ell(w)}
w\left( \prod_{\alpha\in R^+} \Big({1\over 1-tx^\alpha}\Big)
x^{\mu+\rho}\right) \cr
&= \hbox{coefficient of $s_\lambda$ in }\ \ 
\left(\prod_{\alpha\in R^+}{1\over 1-tR_\alpha}\right) s_\mu\,. \cr
}$$
(c)  
$$\eqalign{
K_{\lambda\mu}(t)
&= \hbox{coefficient of $s_\lambda$ in } \ \ 
{1\over a_\rho}\sum_{w\in W} (-1)^{\ell(w)}
w\left( \prod_{\alpha\in R^+} \Big({1\over 1-tx^\alpha}\Big)
x^{\mu+\rho}\right) \cr
&= \hbox{coefficient of $a_{\lambda+\rho}$ in }\ \ 
\sum_{w\in W} (-1)^{\ell(w)}
w\left( \Big(\sum_{\gamma\in Q^+} F(\gamma;t)x^\gamma\Big)
x^{\mu+\rho}\right) \cr
&= \hbox{coefficient of $x^{\lambda+\rho}$ in }\ \ 
\sum_{w\in W} (-1)^{\ell(w)} w\left( \sum_{\gamma\in Q^+} 
F(\gamma;t)x^{\gamma+\mu+\rho}\right)
\cr
&=\sum_{w\in W}(-1)^{\ell(w)} F(w(\lambda+\rho)-(\mu+\rho);t),
\cr
}$$
since $w^{-1}(\gamma+(\mu+\rho))=\lambda+\rho$ implies
$\gamma = w(\lambda+\rho)-(\mu+\rho)$.

\smallskip\noindent
(d) Let $\lambda\in P^+$.  By Theorem 2.22 and Lemma 2.18
$$\eqalign{
\sum_{x\le n_\lambda} q^{-(\ell(n_\lambda)-\ell(x))}
P_{x,n_\lambda}(q^2)T_x
&= C_{n_\lambda}' 
= q^{-\ell(w_0)}W_0(q^2) s_\lambda {\bf 1}_0 \cr
& = q^{-\ell(w_0)}W_0(q^2) \sum_{\mu\le \lambda} K_{\lambda\mu}(q^{-2})
P_\mu(x;q^{-2}) {\bf 1}_0 \cr
&= q^{-\ell(w_0)}W_0(q^2) \sum_{\mu\le \lambda} K_{\lambda\mu}(q^{-2})
{W_0(q^{-2})\over W_\mu(q^{-2})} M_\mu \cr
&= \sum_{\mu\le \lambda} K_{\lambda\mu}(q^{-2})
\sum_{x\in Wt_\mu W} q^{\ell(x)-\ell(n_\mu)} T_x.\cr
}$$
Hence, for $\mu\le \lambda$ and $x\in Wt_\mu W$,
$$K_{\lambda\mu}(q^{-2})
=q^{\ell(n_\mu)-\ell(n_\lambda)}P_{x,n_\lambda}(q^2).$$
By (2.19) and (2.20),
$$\ell(n_\mu)-\ell(n_\lambda) = 
\ell(t_\mu)+\ell(w_0)-(\ell(t_\lambda)+\ell(w_0))
=2\langle \mu,\rho^\vee\rangle - 2\langle \lambda,\rho^\vee\rangle,
$$
and the result follows by replacing $t=q^2$.
\endpf

With notations as in Remark 3.7(b), Theorem 3.17(a) together with
the fact that $s_0=P_0(x;t)$ give the following important
formula for the Kostka-Foulkes polynomial in the case that $\mu=0$,
$$K_{\lambda,0}(t)
=W_0(t)\langle s_\lambda, P_0(x;t)\rangle_t
=W_0(t) \langle s_\lambda, s_0\rangle_t
=\sum_{k\ge 0} \dim(\Hom_{\fg}(L(\lambda),{\cal H}^k)t^k.
\formula$$

\section 4.  The positive formula

In the type A case Lascoux and Sch\"utzenberger [LS]
have used the theory of column strict tableaux to
give a positive formula for the Kostka-Foulkes
polynomial.   In this section we give a proof of this formula. 
Versions of this proof have appeared previously in [Sch] and in [Bt].

The starting point is the formula for $K_{\lambda\mu}(t)$
in Theorem 3.17(a).  To match the setup in [Mac]
we shall work in a slightly different setting 
(corresponding to the Weyl group $W$ and the
weight lattice of the reductive group $GL_n(\CC)$).
In this case the vector space $\fh_\RR^*= \RR^n$ has orthonormal basis
$\varepsilon_1,\ldots, \varepsilon_n$, where
$\varepsilon_i=(0,\ldots,0,1,0,\ldots,0)$ with the 1 in the $i$th spot,  
the Weyl group is the symmetric group $S_n$ acting on $\RR^n$ 
by permuting the coordinates, the weight lattice $P$ is 
replaced by the lattice
$$\ZZ^n = \{ (\gamma_1,\ldots,\gamma_n)\ |\ \gamma_i\in \ZZ\}
\qquad\hbox{and}\qquad
\delta = (n-1,n-2,\ldots, 2,1,0)
\formula$$
replaces the element $\rho$.  The 
positive roots are given by
$R^+= \{ \varepsilon_i-\varepsilon_j\ |\ 1\le i<j\le n\}$
and the Schur functions (defined as in (2.8)) are viewed as
(Laurent) polynomials in the variables $x_1,\ldots, x_n$,
where $x_i = x^{\varepsilon_i}$ and 
the symmetric group $S_n$ acts by permuting the variables.
If $w\in S_n$ then $(-1)^{\ell(w)}=\det(w)$ is the {\it sign} of 
the permutation $w$ and
the straightening law for Schur functions 
(see (2.9) and [Mac, I paragraph after (3.1)])
is
$$s_\mu = (-1)^{\ell(w)}s_{w\circ \mu},
\qquad\hbox{where}\qquad
w\circ\mu = w(\mu+\delta)-\delta,\formula$$
for $w\in S_n$ and $\mu\in \ZZ^n$.
The set of {\it partitions}
$${\cal P} = \{ (\lambda_1, \ldots,\lambda_n)\in \ZZ^n\ |\ 
\lambda_1\ge \cdots\ge \lambda_n\ge 0\}
\formula$$
takes the role played by the set $P^+$.
Conforming to the conventions in [Mac] so that gravity goes
up and to the left, each partition $\mu = 
(\mu_1,\ldots, \mu_n)\in {\cal P}$ is identified
with the collection of boxes in a corner which has
$\mu_i$ boxes in row $i$,
where, as for matrices,
the rows and columns of $\mu$ are indexed from
top to bottom and left to right, respectively.
For example,
$$
\beginpicture
\setcoordinatesystem units <0.25cm,0.25cm>         
\setplotarea x from -7 to 7, y from 0 to 6    
\linethickness=0.5pt                          
\putrule from 0 6 to 5 6          %
\putrule from 0 5 to 5 5          
\putrule from 0 4 to 5 4          %
\putrule from 0 3 to 3 3          %
\putrule from 0 2 to 3 2          %
\putrule from 0 1 to 1 1          %
\putrule from 0 0 to 1 0          %
\putrule from 0 0 to 0 6        %
\putrule from 1 0 to 1 6        %
\putrule from 2 2 to 2 6        %
\putrule from 3 2 to 3 6        
\putrule from 4 4 to 4 6        %
\putrule from 5 4 to 5 6        %
\put{$(5,5,3,3,1,1)=$} at -6 3.5
\put{.} at 6 3
\endpicture
$$

For each pair $1\le i<j\le n$ define the {\it raising operator}
$R_{ij}\colon \ZZ^n \to \ZZ^n$  (see (3.15) and [Mac, I \S 1 (1.14)])
by
$$R_{ij}\mu = \mu+\varepsilon_i-\varepsilon_j
\qquad\hbox{and define}\qquad
(R_{i_1j_1}\cdots R_{i_\ell j_\ell})s_\mu
=s_{R_{i_1j_1}\cdots R_{i_\ell j_\ell}\mu},
\formula$$
for a sequence of pairs $i_1<j_1, \ldots, i_\ell<j_\ell$. 
Using the straightening law (4.2) any
Schur function $s_\mu$ indexed by $\mu\in \ZZ^n$ 
with $\mu_1+\cdots+\mu_n\ge 0$
is either equal to 0 or to $\pm s_\lambda$ for some 
$\lambda\in {\cal P}$.  Composing the action of raising operators on
Schur functions $s_\lambda$ should be avoided.  For example, if
$n=2$ and $s_1$ denotes the transposition in the symmetric group $S_2$ 
then, by the straightening law, $s_{(0,1)} = - s_{s_1((0,1)+(1,0))-(1,0)}
=-s_{(1,1)-(1,0)} = -s_{(0,1)}$ giving that
$s_{(0,1)}=0$ and so
$$R_{12}(R_{12}s_{(-1,2)}) = R_{12}s_{(0,1)}= R_{12}\cdot 0=0,
\qquad\hbox{whereas}\qquad
(R_{12}^2)s_{(-1,2)} = s_{(1,0)}=x_1+x_2.$$

With notation as in (4.2) and (4.4) we may define the 
{\it Hall-Littlewood polynomials} for this type A case
by (see Theorem 3.17(b) and [Mac, III (4.6)])
$$Q_\mu = \left(\prod_{1\le i<j\le n} {1\over 1-tR_{ij}}\right) s_\mu,
\qquad\hbox{for all $\mu\in \ZZ^n$,}
\formula$$
and the {\it Kostka-Foulkes polynomials} $K_{\lambda\mu}(t)$, 
$\lambda,\mu\in {\cal P}$, by 
$$Q_\mu = \sum_{\lambda\in {\cal P}} K_{\lambda\mu}(t)s_\lambda.
\formula$$

\subsection Insertion and Pieri rules

Let $\lambda$ and $\mu = (\mu_1,\ldots, \mu_n)$ be partitions.
A {\it column strict tableau of shape $\lambda$ and
weight $\mu$} is a filling of the boxes of $\lambda$
with $\mu_1$ ~1s, $\mu_2$ ~2s, $\ldots$,
$\mu_n$ ~$n$s, such that 
\smallskip\noindent
\itemitem{(a)} the rows are weakly increasing from left to right,
\smallskip\noindent
\itemitem{(b)} the columns are strictly increasing from 
top to bottom.
\smallskip\noindent
If $T$ is a column strict tableau write
${\rm shp}(T)$ and ${\rm wt}(T)$ for the shape and the
weight of $T$ so that 
$$\eqalign{
{\rm shp}(T) &= (\lambda_1,\ldots, \lambda_n),
\qquad\hbox{where}\quad
\lambda_i = \hbox{number of boxes in row $i$ of $T$,\quad and} \cr
{\rm wt}(T) &= (\mu_1,\ldots,\mu_n),
\qquad\hbox{where}\quad
\mu_i = \hbox{number of $i\,$s in $T$}. \cr
}$$
For example,
\smallskip\noindent
$$
\beginpicture
\setcoordinatesystem units <0.5cm,0.5cm>         
\setplotarea x from -2 to 15, y from 0 to 3    
\linethickness=0.5pt                          
\put{$T=$} at -1.5 4.5
\put{has\quad ${\rm shp}(T) = (9,7,7,4,2,1,0)$\quad and}[l] at 11 3.5
\put{\phantom{has}\quad ${\rm wt}(T) = (7,6,5,5,3,2,2)$.}[l] at  11 2.5
\put{7} at 0.5 0.5
\put{6} at 0.5 1.5
\put{4} at 0.5 2.5
\put{3} at 0.5 3.5
\put{2} at 0.5 4.5
\put{1} at 0.5 5.5
\put{7} at 1.5 1.5
\put{5} at 1.5 2.5
\put{3} at 1.5 3.5
\put{2} at 1.5 4.5
\put{1} at 1.5 5.5
\put{5} at 2.5 2.5
\put{3} at 2.5 3.5
\put{2} at 2.5 4.5
\put{1} at 2.5 5.5
\put{6}  at 3.5 2.5
\put{4}  at 3.5 3.5
\put{2}  at 3.5 4.5
\put{1}  at 3.5 5.5
\put{4} at 4.5 3.5
\put{3} at 4.5 4.5
\put{1} at 4.5 5.5
\put{4}  at 5.5 3.5
\put{3}  at 5.5 4.5
\put{1}  at 5.5 5.5
\put{5}  at 6.5 3.5
\put{4}  at 6.5 4.5
\put{1}  at 6.5 5.5
\put{2}  at 7.5 5.5
\put{2}  at 8.5 5.5
\putrule from 0 6 to 9 6          %
\putrule from 0 5 to 9 5          
\putrule from 0 4 to 7 4          %
\putrule from 0 3 to 7 3          %
\putrule from 0 2 to 4 2          %
\putrule from 0 2 to 2 2          %
\putrule from 0 1 to 2 1          %
\putrule from 0 0 to 1 0          %
\putrule from 0 0 to 0 6        %
\putrule from 1 0 to 1 6        %
\putrule from 2 1 to 2 6        %
\putrule from 3 2 to 3 6        
\putrule from 4 2 to 4 6        %
\putrule from 5 3 to 5 6        %
\putrule from 6 3 to 6 6        %
\putrule from 7 3 to 7 6        %
\putrule from 8 5 to 8 6        %
\putrule from 9 5 to 9 6        %
\endpicture
$$
For partitions $\lambda$ and $\mu$ and, more generally,
for any two sets ${\cal S}, {\cal W}\subseteq {\cal P}$ 
write
$$\eqalign{
B(\lambda) &= \{ 
\hbox{column strict tableaux $T$}\ |\ 
{\rm shp}(T) = \lambda\}, \cr
B(\lambda)_\mu &= \{ 
\hbox{column strict tableaux $T$}\ |\ 
{\rm shp}(T) = \lambda\ \hbox{and}\ {\rm wt}(T)=\mu\}, 
\cr
B({\cal S})_{\cal W} &=\{
\hbox{column strict tableaux $T$}\ |\ 
{\rm shp}(T) \in {\cal S}\ \hbox{and}\ {\rm wt}(T)\in {\cal W}\}. \cr
}
\formula$$ 

Let $\lambda$ and $\gamma$ be partitions such that $\gamma\subseteq \lambda$
(as collections of boxes in a corner, i.e. $\gamma_i\le \lambda_i$
for $1\le i\le n$).  The {\it skew shape} $\lambda/\gamma$
is the collection of boxes of $\lambda$ which are not in $\gamma$.
The {\it jeu de taquin} reduces
a column strict filling of a skew shape $\lambda/\gamma$ to a
column strict tableau of partition shape.  
At each step ``gravity'' moves one box up or to the left
without violating the column strict condition (weakly increasing in rows,
strictly increasing in columns).
The jeu de taquin is most easily illustrated by example:
$$
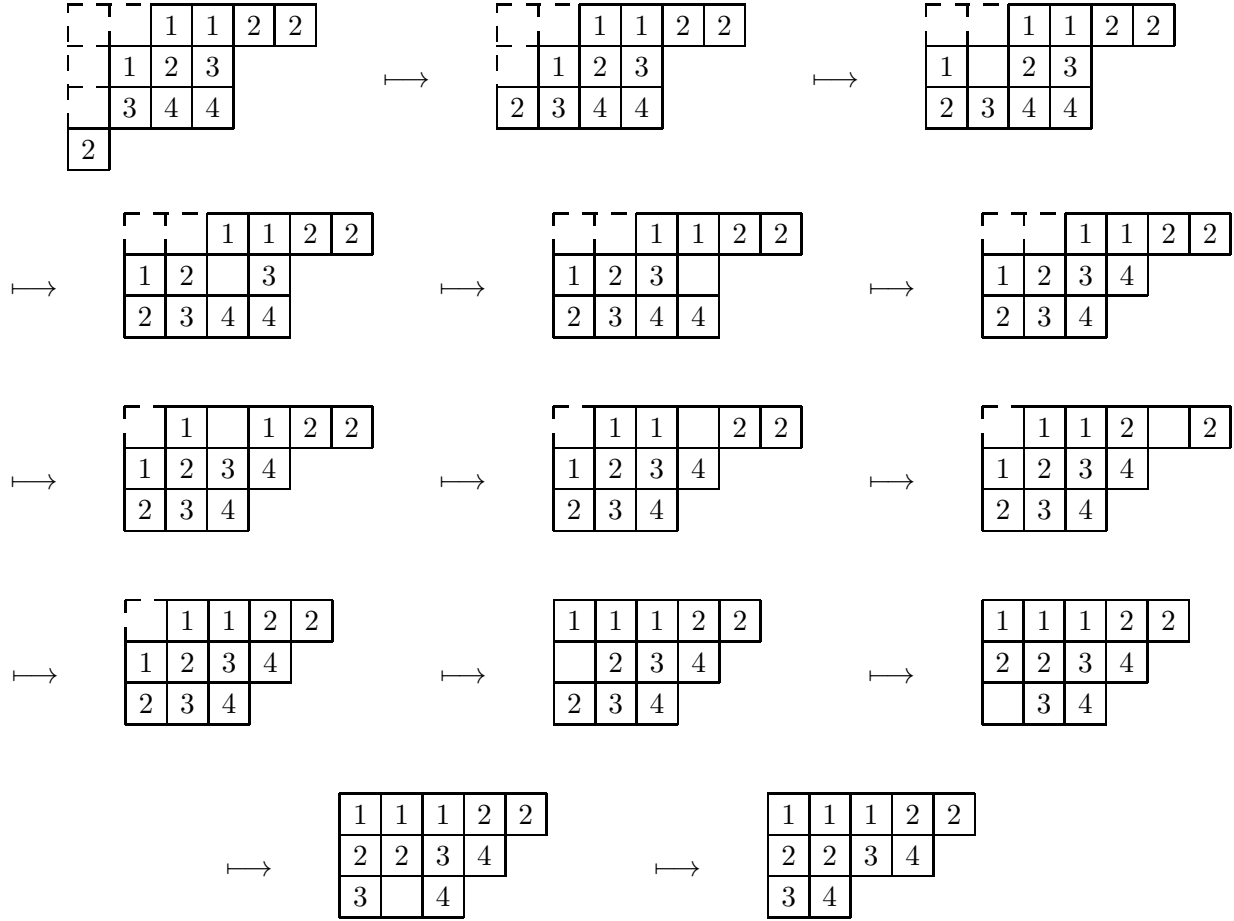
\beginpicture
\setcoordinatesystem units <0.5cm,0.5cm>         
\setplotarea x from 0 to 6, y from -2 to 2    
\linethickness=0.5pt                          
\put{2}  at 0.5 -1.5
\put{3} at 1.5 -0.5
\put{1} at 1.5 0.5
\put{4}  at 2.5 -0.5
\put{2}  at 2.5 0.5
\put{1}  at 2.5 1.5
\put{4}  at 3.5 -0.5
\put{3}  at 3.5 0.5
\put{1}  at 3.5 1.5
\put{2}  at 4.5 1.5
\put{2}  at 5.5 1.5
\putrule from 2 2 to 6 2          %
\putrule from 1 1 to 6 1          
\putrule from 1 0 to 4 0          %
\putrule from 0 -1 to 4 -1          %
\putrule from 0 -2 to 1 -2          %
\putrule from 0 -2 to 0 -1        
\putrule from 1 -2 to 1 1        %
\putrule from 2 -1 to 2 2        %
\putrule from 3 -1 to 3 2        %
\putrule from 4 -1 to 4 2        %
\putrule from 5 1 to 5 2        %
\putrule from 6 1 to 6 2        %
\setdashes
\putrule from 0 2 to 2 2          %
\putrule from 0 1 to 1 1          
\putrule from 0 0 to 1 0          %
\putrule from 0 -1 to 0 2       
\putrule from 1 1 to 1 2        %
\endpicture
\qquad\longmapsto\qquad
\beginpicture
\setcoordinatesystem units <0.5cm,0.5cm>         
\setplotarea x from 0 to 6, y from -2 to 2    
\linethickness=0.5pt                          
\put{2}  at 0.5 -0.5
\put{3} at 1.5 -0.5
\put{1} at 1.5 0.5
\put{4}  at 2.5 -0.5
\put{2}  at 2.5 0.5
\put{1}  at 2.5 1.5
\put{4}  at 3.5 -0.5
\put{3}  at 3.5 0.5
\put{1}  at 3.5 1.5
\put{2}  at 4.5 1.5
\put{2}  at 5.5 1.5
\putrule from 2 2 to 6 2          %
\putrule from 1 1 to 6 1          
\putrule from 0 0 to 4 0          %
\putrule from 0 -1 to 4 -1        %
\putrule from 0 -1 to 0 0      
\putrule from 1 -1 to 1 1      %
\putrule from 2 -1 to 2 2      %
\putrule from 3 -1 to 3 2      %
\putrule from 4 -1 to 4 2      %
\putrule from 5 1 to 5 2       %
\putrule from 6 1 to 6 2       %
\setdashes
\putrule from 0 2 to 2 2          %
\putrule from 0 1 to 1 1          
\putrule from 0 -1 to 0 2       
\putrule from 1 1 to 1 2        %
\endpicture
\qquad\longmapsto\qquad
\beginpicture
\setcoordinatesystem units <0.5cm,0.5cm>         
\setplotarea x from 0 to 6, y from -2 to 2    
\linethickness=0.5pt                          
\put{2}  at 0.5 -0.5
\put{3} at 1.5 -0.5
\put{1} at 0.5 0.5
\put{4}  at 2.5 -0.5
\put{2}  at 2.5 0.5
\put{1}  at 2.5 1.5
\put{4}  at 3.5 -0.5
\put{3}  at 3.5 0.5
\put{1}  at 3.5 1.5
\put{2}  at 4.5 1.5
\put{2}  at 5.5 1.5
\putrule from 2 2 to 6 2          %
\putrule from 0 1 to 6 1          
\putrule from 0 0 to 4 0          %
\putrule from 0 -1 to 4 -1        %
\putrule from 0 -1 to 0 1      
\putrule from 1 -1 to 1 1      %
\putrule from 2 -1 to 2 2      %
\putrule from 3 -1 to 3 2      %
\putrule from 4 -1 to 4 2      %
\putrule from 5 1 to 5 2       %
\putrule from 6 1 to 6 2       %
\setdashes
\putrule from 0 2 to 2 2          
\putrule from 0 1 to 0 2        
\putrule from 1 1 to 1 2        %
\endpicture
$$
$$
\longmapsto\qquad
\beginpicture
\setcoordinatesystem units <0.5cm,0.5cm>         
\setplotarea x from 0 to 6, y from -2 to 2    
\linethickness=0.5pt                          
\put{2}  at 0.5 -0.5
\put{3} at 1.5 -0.5
\put{1} at 0.5 0.5
\put{4}  at 2.5 -0.5
\put{2}  at 1.5 0.5
\put{1}  at 2.5 1.5
\put{4}  at 3.5 -0.5
\put{3}  at 3.5 0.5
\put{1}  at 3.5 1.5
\put{2}  at 4.5 1.5
\put{2}  at 5.5 1.5
\putrule from 2 2 to 6 2          %
\putrule from 0 1 to 6 1          
\putrule from 0 0 to 4 0          %
\putrule from 0 -1 to 4 -1        %
\putrule from 0 -1 to 0 1      
\putrule from 1 -1 to 1 1      %
\putrule from 2 -1 to 2 2      %
\putrule from 3 -1 to 3 2      %
\putrule from 4 -1 to 4 2      %
\putrule from 5 1 to 5 2       %
\putrule from 6 1 to 6 2       %
\setdashes
\putrule from 0 2 to 2 2          
\putrule from 0 1 to 0 2        
\putrule from 1 1 to 1 2        %
\endpicture
\qquad\longmapsto\qquad
\beginpicture
\setcoordinatesystem units <0.5cm,0.5cm>         
\setplotarea x from 0 to 6, y from -2 to 2    
\linethickness=0.5pt                          
\put{2}  at 0.5 -0.5
\put{3} at 1.5 -0.5
\put{1} at 0.5 0.5
\put{4}  at 2.5 -0.5
\put{2}  at 1.5 0.5
\put{1}  at 2.5 1.5
\put{4}  at 3.5 -0.5
\put{3}  at 2.5 0.5
\put{1}  at 3.5 1.5
\put{2}  at 4.5 1.5
\put{2}  at 5.5 1.5
\putrule from 2 2 to 6 2          %
\putrule from 0 1 to 6 1          
\putrule from 0 0 to 4 0          %
\putrule from 0 -1 to 4 -1        %
\putrule from 0 -1 to 0 1      
\putrule from 1 -1 to 1 1      %
\putrule from 2 -1 to 2 2      %
\putrule from 3 -1 to 3 2      %
\putrule from 4 -1 to 4 2      %
\putrule from 5 1 to 5 2       %
\putrule from 6 1 to 6 2       %
\setdashes
\putrule from 0 2 to 2 2          
\putrule from 0 1 to 0 2        
\putrule from 1 1 to 1 2        %
\endpicture
\qquad\longmapsto\qquad
\beginpicture
\setcoordinatesystem units <0.5cm,0.5cm>         
\setplotarea x from 0 to 6, y from -2 to 2    
\linethickness=0.5pt                          
\put{2}  at 0.5 -0.5
\put{3} at 1.5 -0.5
\put{1} at 0.5 0.5
\put{4}  at 2.5 -0.5
\put{2}  at 1.5 0.5
\put{1}  at 2.5 1.5
\put{4}  at 3.5 0.5
\put{3}  at 2.5 0.5
\put{1}  at 3.5 1.5
\put{2}  at 4.5 1.5
\put{2}  at 5.5 1.5
\putrule from 2 2 to 6 2          %
\putrule from 0 1 to 6 1          
\putrule from 0 0 to 4 0          %
\putrule from 0 -1 to 3 -1        %
\putrule from 0 -1 to 0 1      
\putrule from 1 -1 to 1 1      %
\putrule from 2 -1 to 2 2      %
\putrule from 3 -1 to 3 2      %
\putrule from 4 0 to 4 2       %
\putrule from 5 1 to 5 2       %
\putrule from 6 1 to 6 2       %
\setdashes
\putrule from 0 2 to 2 2          
\putrule from 0 1 to 0 2        
\putrule from 1 1 to 1 2        %
\endpicture
$$
$$
\longmapsto\qquad
\beginpicture
\setcoordinatesystem units <0.5cm,0.5cm>         
\setplotarea x from 0 to 6, y from -2 to 2    
\linethickness=0.5pt                          
\put{2}  at 0.5 -0.5
\put{3} at 1.5 -0.5
\put{1} at 0.5 0.5
\put{4}  at 2.5 -0.5
\put{2}  at 1.5 0.5
\put{1}  at 1.5 1.5
\put{4}  at 3.5 0.5
\put{3}  at 2.5 0.5
\put{1}  at 3.5 1.5
\put{2}  at 4.5 1.5
\put{2}  at 5.5 1.5
\putrule from 1 2 to 6 2          %
\putrule from 0 1 to 6 1          
\putrule from 0 0 to 4 0          %
\putrule from 0 -1 to 3 -1        %
\putrule from 0 -1 to 0 1      
\putrule from 1 -1 to 1 2      %
\putrule from 2 -1 to 2 2      %
\putrule from 3 -1 to 3 2      %
\putrule from 4 0 to 4 2       %
\putrule from 5 1 to 5 2       %
\putrule from 6 1 to 6 2       %
\setdashes
\putrule from 0 2 to 1 2          
\putrule from 0 1 to 0 2        
\endpicture
\qquad\longmapsto\qquad
\beginpicture
\setcoordinatesystem units <0.5cm,0.5cm>         
\setplotarea x from 0 to 6, y from -2 to 2    
\linethickness=0.5pt                          
\put{2}  at 0.5 -0.5
\put{3} at 1.5 -0.5
\put{1} at 0.5 0.5
\put{4}  at 2.5 -0.5
\put{2}  at 1.5 0.5
\put{1}  at 1.5 1.5
\put{4}  at 3.5 0.5
\put{3}  at 2.5 0.5
\put{1}  at 2.5 1.5
\put{2}  at 4.5 1.5
\put{2}  at 5.5 1.5
\putrule from 1 2 to 6 2          %
\putrule from 0 1 to 6 1          
\putrule from 0 0 to 4 0          %
\putrule from 0 -1 to 3 -1        %
\putrule from 0 -1 to 0 1      
\putrule from 1 -1 to 1 2      %
\putrule from 2 -1 to 2 2      %
\putrule from 3 -1 to 3 2      %
\putrule from 4 0 to 4 2       %
\putrule from 5 1 to 5 2       %
\putrule from 6 1 to 6 2       %
\setdashes
\putrule from 0 2 to 1 2          
\putrule from 0 1 to 0 2        
\endpicture
\qquad\longmapsto\qquad
\beginpicture
\setcoordinatesystem units <0.5cm,0.5cm>         
\setplotarea x from 0 to 6, y from -2 to 2    
\linethickness=0.5pt                          
\put{2}  at 0.5 -0.5
\put{3} at 1.5 -0.5
\put{1} at 0.5 0.5
\put{4}  at 2.5 -0.5
\put{2}  at 1.5 0.5
\put{1}  at 1.5 1.5
\put{4}  at 3.5 0.5
\put{3}  at 2.5 0.5
\put{1}  at 2.5 1.5
\put{2}  at 3.5 1.5
\put{2}  at 5.5 1.5
\putrule from 1 2 to 6 2          %
\putrule from 0 1 to 6 1          
\putrule from 0 0 to 4 0          %
\putrule from 0 -1 to 3 -1        %
\putrule from 0 -1 to 0 1      
\putrule from 1 -1 to 1 2      %
\putrule from 2 -1 to 2 2      %
\putrule from 3 -1 to 3 2      %
\putrule from 4 0 to 4 2       %
\putrule from 5 1 to 5 2       %
\putrule from 6 1 to 6 2       %
\setdashes
\putrule from 0 2 to 1 2          
\putrule from 0 1 to 0 2        
\endpicture
$$
$$\longmapsto\qquad
\beginpicture
\setcoordinatesystem units <0.5cm,0.5cm>         
\setplotarea x from 0 to 6, y from -2 to 2    
\linethickness=0.5pt                          
\put{2}  at 0.5 -0.5
\put{3} at 1.5 -0.5
\put{1} at 0.5 0.5
\put{4}  at 2.5 -0.5
\put{2}  at 1.5 0.5
\put{1}  at 1.5 1.5
\put{4}  at 3.5 0.5
\put{3}  at 2.5 0.5
\put{1}  at 2.5 1.5
\put{2}  at 3.5 1.5
\put{2}  at 4.5 1.5
\putrule from 1 2 to 5 2          %
\putrule from 0 1 to 5 1          
\putrule from 0 0 to 4 0          %
\putrule from 0 -1 to 3 -1        %
\putrule from 0 -1 to 0 1      
\putrule from 1 -1 to 1 2      %
\putrule from 2 -1 to 2 2      %
\putrule from 3 -1 to 3 2      %
\putrule from 4 0 to 4 2       %
\putrule from 5 1 to 5 2       %
\setdashes
\putrule from 0 2 to 1 2          
\putrule from 0 1 to 0 2        
\endpicture
\qquad\longmapsto\qquad
\beginpicture
\setcoordinatesystem units <0.5cm,0.5cm>         
\setplotarea x from 0 to 6, y from -2 to 2    
\linethickness=0.5pt                          
\put{2}  at 0.5 -0.5
\put{3} at 1.5 -0.5
\put{1} at 0.5 1.5
\put{4}  at 2.5 -0.5
\put{2}  at 1.5 0.5
\put{1}  at 1.5 1.5
\put{4}  at 3.5 0.5
\put{3}  at 2.5 0.5
\put{1}  at 2.5 1.5
\put{2}  at 3.5 1.5
\put{2}  at 4.5 1.5
\putrule from 0 2 to 5 2          %
\putrule from 0 1 to 5 1          
\putrule from 0 0 to 4 0          %
\putrule from 0 -1 to 3 -1        %
\putrule from 0 -1 to 0 2      
\putrule from 1 -1 to 1 2      %
\putrule from 2 -1 to 2 2      %
\putrule from 3 -1 to 3 2      %
\putrule from 4 0 to 4 2       %
\putrule from 5 1 to 5 2       %
\endpicture
\qquad\longmapsto\qquad
\beginpicture
\setcoordinatesystem units <0.5cm,0.5cm>         
\setplotarea x from 0 to 6, y from -2 to 2    
\linethickness=0.5pt                          
\put{2}  at 0.5 0.5
\put{3} at 1.5 -0.5
\put{1} at 0.5 1.5
\put{4}  at 2.5 -0.5
\put{2}  at 1.5 0.5
\put{1}  at 1.5 1.5
\put{4}  at 3.5 0.5
\put{3}  at 2.5 0.5
\put{1}  at 2.5 1.5
\put{2}  at 3.5 1.5
\put{2}  at 4.5 1.5
\putrule from 0 2 to 5 2          %
\putrule from 0 1 to 5 1          
\putrule from 0 0 to 4 0          %
\putrule from 0 -1 to 3 -1        %
\putrule from 0 -1 to 0 2      
\putrule from 1 -1 to 1 2      %
\putrule from 2 -1 to 2 2      %
\putrule from 3 -1 to 3 2      %
\putrule from 4 0 to 4 2       %
\putrule from 5 1 to 5 2       %
\endpicture
$$
$$\longmapsto\qquad
\beginpicture
\setcoordinatesystem units <0.5cm,0.5cm>         
\setplotarea x from 0 to 6, y from -2 to 2    
\linethickness=0.5pt                          
\put{2}  at 0.5 0.5
\put{3} at 0.5 -0.5
\put{1} at 0.5 1.5
\put{4}  at 2.5 -0.5
\put{2}  at 1.5 0.5
\put{1}  at 1.5 1.5
\put{4}  at 3.5 0.5
\put{3}  at 2.5 0.5
\put{1}  at 2.5 1.5
\put{2}  at 3.5 1.5
\put{2}  at 4.5 1.5
\putrule from 0 2 to 5 2          %
\putrule from 0 1 to 5 1          
\putrule from 0 0 to 4 0          %
\putrule from 0 -1 to 3 -1        %
\putrule from 0 -1 to 0 2      
\putrule from 1 -1 to 1 2      %
\putrule from 2 -1 to 2 2      %
\putrule from 3 -1 to 3 2      %
\putrule from 4 0 to 4 2       %
\putrule from 5 1 to 5 2       %
\endpicture
\qquad\longmapsto\qquad
\beginpicture
\setcoordinatesystem units <0.5cm,0.5cm>         
\setplotarea x from 0 to 6, y from -2 to 2    
\linethickness=0.5pt                          
\put{2}  at 0.5 0.5
\put{3} at 0.5 -0.5
\put{1} at 0.5 1.5
\put{4}  at 1.5 -0.5
\put{2}  at 1.5 0.5
\put{1}  at 1.5 1.5
\put{4}  at 3.5 0.5
\put{3}  at 2.5 0.5
\put{1}  at 2.5 1.5
\put{2}  at 3.5 1.5
\put{2}  at 4.5 1.5
\putrule from 0 2 to 5 2          %
\putrule from 0 1 to 5 1          
\putrule from 0 0 to 4 0          %
\putrule from 0 -1 to 2 -1        %
\putrule from 0 -1 to 0 2      
\putrule from 1 -1 to 1 2      %
\putrule from 2 -1 to 2 2      %
\putrule from 3 0 to 3 2      %
\putrule from 4 0 to 4 2       %
\putrule from 5 1 to 5 2       %
\endpicture
$$
The result of the jeu de taquin is independent
of the choice of order of the moves ([Fu, \S 1.2 Claim 2] which
is proved in [Fu, \S 2 and \S 3]).

The {\it plactic monoid} is the set $B({\cal P})$ of column strict
tableaux with product given by
$$T_1*T_2 = 
\hbox{jeu de taquin reduction of } \quad 
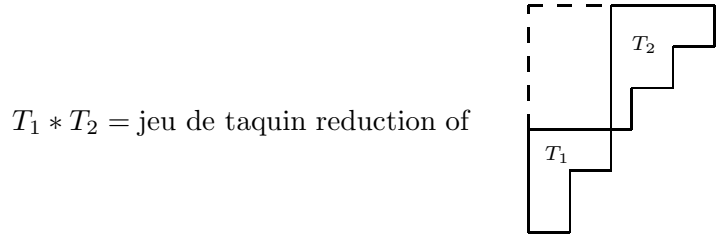
\beginpicture
\setcoordinatesystem units <0.25cm,0.25cm>         
\setplotarea x from -5 to 10, y from -6 to 7    
\linethickness=0.5pt                          
\put{$\scriptstyle{T_1}$}[tl] at -3.2 -0.8
\put{$\scriptstyle{T_2}$}[tl] at 1 4.5
\putrule from 0 6 to 5 6          
\putrule from 3 4 to 5 4          
\putrule from 1 2 to 3 2          %
\putrule from 0 0 to 1 0          %
\putrule from 0 0 to 0 6        %
\putrule from 1 0 to 1 2        %
\putrule from 3 2 to 3 4        
\putrule from 5 4 to 5 6        %
\putrule from -4 0 to 0 0          %
\putrule from -2 -2 to 0 -2        
\putrule from -4 -5 to -2 -5       
\putrule from 0 -2 to 0 0       
\putrule from -2 -5 to -2 -2    
\putrule from -4 -5 to -4 0     %
\setdashes
\putrule from -4 6 to 4 6
\putrule from -4 0 to -4 6
\endpicture
$$
This is an associative monoid ([Fu, \S 1.1 Claim 1] which is 
proved in [Fu, \S 2 and \S 3]).

If $x$ is a ``letter'', i.e. a column strict tableau of shape $(1)=\square$,
then
$$\eqalign{
&\hbox{$x*T$ is the {\it column insertion} of $x$ into $T$,}
\qquad\hbox{and}\cr
&\hbox{$T*x$ is the {\it row insertion} of $x$ into $T$.} \cr
}\formula
$$
The shape $\lambda$ of $P=T*x$ differs from the shape $\gamma$
of $T$ by single box and so if $\gamma$ and $P$ are given
then the pair $(T,x)$ can be recovered by ``uninserting'' the
box $\lambda/\gamma$ from $P$.  The tableaux $P$ and $T$ differ
by at most one entry in each row. The entries where $P$ and $T$ differ 
form the {\it bumping path of $x$}.  The bumping path begins with $x$
in the first row of $P$ and ends at the entry in the box $\lambda/\gamma$.
For example,
$$
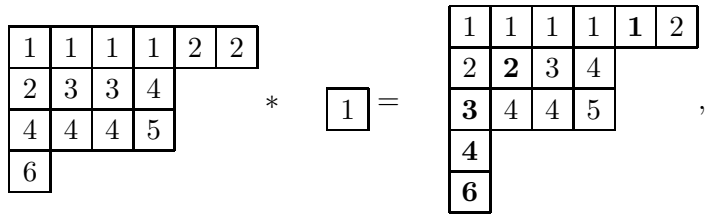
\beginpicture
\setcoordinatesystem units <0.5cm,0.5cm>         
\setplotarea x from -1 to 6, y from -2 to 2    
\linethickness=0.5pt                          
\put{6}  at 0.5 -1.5
\put{4}  at 0.5 -0.5
\put{2}  at 0.5 0.5
\put{1}  at 0.5 1.5
\put{4} at 1.5 -0.5
\put{3} at 1.5 0.5
\put{1} at 1.5 1.5
\put{4}  at 2.5 -0.5
\put{3}  at 2.5 0.5
\put{1}  at 2.5 1.5
\put{5}  at 3.5 -0.5
\put{4}  at 3.5 0.5
\put{1}  at 3.5 1.5
\put{2}  at 4.5 1.5
\put{2}  at 5.5 1.5
\putrule from 0 2 to 6 2          %
\putrule from 0 1 to 6 1          
\putrule from 0 0 to 4 0          %
\putrule from 0 -1 to 4 -1          %
\putrule from 0 -2 to 1 -2          %
\putrule from 0 -2 to 0 2        
\putrule from 1 -2 to 1 2        %
\putrule from 2 -1 to 2 2        %
\putrule from 3 -1 to 3 2        %
\putrule from 4 -1 to 4 2        %
\putrule from 5 1 to 5 2        %
\putrule from 6 1 to 6 2        %
\endpicture
*
\beginpicture
\setcoordinatesystem units <0.5cm,0.5cm>         
\setplotarea x from -1 to 1, y from -2 to 2    
\linethickness=0.5pt                          
\put{1}  at 0.5 0
\putrule from 0 -0.5 to 1 -0.5          
\putrule from 0 0.5 to 1 0.5          %
\putrule from 0 -0.5 to 0 0.5        
\putrule from 1 -0.5 to 1 0.5        %
\endpicture
=
\beginpicture
\setcoordinatesystem units <0.5cm,0.5cm>         
\setplotarea x from -1 to 6, y from -2 to 2    
\linethickness=0.5pt                          
\put{{\bf 6}}  at 0.5 -2
\put{{\bf 4}}  at 0.5 -1
\put{{\bf 3}}  at 0.5 0
\put{2}  at 0.5 1
\put{1}  at 0.5 2
\put{4} at 1.5  0
\put{{\bf 2}} at 1.5 1
\put{1} at 1.5 2
\put{4}  at 2.5 0
\put{3}  at 2.5 1
\put{1}  at 2.5 2
\put{5}  at 3.5 0
\put{4}  at 3.5 1
\put{1}  at 3.5 2
\put{{\bf 1}}  at 4.5 2
\put{2}  at 5.5 2
\putrule from 0 2.5 to 6 2.5          %
\putrule from 0 1.5 to 6 1.5          
\putrule from 0 0.5 to 4 0.5          %
\putrule from 0 -0.5 to 4 -0.5          %
\putrule from 0 -1.5 to 1 -1.5          %
\putrule from 0 -2.5 to 1 -2.5        %
\putrule from 0 -2.5 to 0 2.5        
\putrule from 1 -2.5 to 1 2.5        %
\putrule from 2 -0.5 to 2 2.5        %
\putrule from 3 -0.5 to 3 2.5        %
\putrule from 4 -0.5 to 4 2.5        %
\putrule from 5 1.5 to 5 2.5        %
\putrule from 6 1.5 to 6 2.5        %
\endpicture
,
$$
where the bold face entries form the bumping path.

The {\it monoid of words} is the free monoid $B^*$ generated
by $\{1,2,\ldots, n\}$.  The {\it weight} ${\rm wt}(w)$ of
a word $w=w_1\cdots w_n$ is 
$${\rm wt}(w)={\rm wt}(w_1\cdots w_n)=(\mu_1,\ldots,\mu_n)
\qquad\hbox{where}\qquad
\hbox{$\mu_i$ is the number of $i$'s in $w$.}$$
For example,
$w=3214566532211$
is a word of weight ${\rm wt}(w) = (3,3,2,1,2,2)$.
The {\it insertion} map
$$\matrix{
B^* &\longrightarrow &B({\cal P}) \cr
w_1\cdots w_n
&\longmapsto &w_1*\cdots*w_n \cr
}
\formula$$
is a weight preserving homomorphism of monoids.

A {\it horizontal strip} is a skew shape which
contains at most one box in each column.  
The {\it length} of a horizontal strip $\lambda/\gamma$ 
is the number of boxes in $\lambda/\gamma$.  
The boxes containing $\times$ in the picture
$$
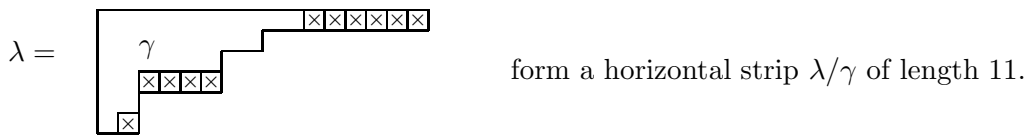
\beginpicture
\setcoordinatesystem units <0.5cm,0.25cm>         
\setplotarea x from -2 to 20, y from 0 to 6    
\linethickness=0.5pt                          
\put{$\lambda = $}[r] at -1 4
\putrule from 0 6 to 8 6          %
\putrule from 4 5 to 8 5          
\putrule from 3 4 to 4 4          %
\putrule from 1 3 to 3 3          %
\putrule from 1 2 to 3 2          %
\putrule from 0.5 1 to 1 1        %
\putrule from 0 0 to 1 0          %
\putrule from 0 0 to 0 6        %
\putrule from 0.5 0 to 0.5 1    %
\putrule from 1 0 to 1 3        %
\putrule from 1.5 2 to 1.5 3    %
\putrule from 2 2 to 2 3        %
\putrule from 2.5 2 to 2.5 3    %
\putrule from 3 2 to 3 4        
\putrule from 4 4 to 4 5        %
\putrule from 5 5 to 5 6        %
\putrule from 5.5 5 to 5.5 6    %
\putrule from 6 5 to 6 6        %
\putrule from 6.5 5 to 6.5 6    %
\putrule from 7 5 to 7 6        %
\putrule from 7.5 5 to 7.5 6    %
\putrule from 8 5 to 8 6        %
\put{$\gamma$}[tl] at 1 4.5
\put{$\scriptstyle{\times}$} at 0.75 0.5
\put{$\scriptstyle{\times}$} at 1.25 2.5
\put{$\scriptstyle{\times}$} at 1.75 2.5
\put{$\scriptstyle{\times}$} at 2.25 2.5
\put{$\scriptstyle{\times}$} at 2.75 2.5
\put{$\scriptstyle{\times}$} at 5.25 5.5
\put{$\scriptstyle{\times}$} at 5.75 5.5
\put{$\scriptstyle{\times}$} at 6.25 5.5
\put{$\scriptstyle{\times}$} at 6.75 5.5
\put{$\scriptstyle{\times}$} at 7.25 5.5
\put{$\scriptstyle{\times}$} at 7.75 5.5
\put{form a horizontal strip $\lambda/\gamma$ of length $11$.}[l]
at 10 3
\endpicture
$$
For a partitions $\mu$ and $\gamma$ and a nonnegative integer $r$ let
$$\eqalign{
\gamma\otimes (r) 
= (r)\otimes \gamma
&= \{ \hbox{partitions $\lambda$}
\ |\ \hbox{$\lambda/\gamma$ is a horizontal strip of length $r$}\},
\cr
(B(r)\otimes B(\gamma))_\mu &= \{
\hbox{pairs $v\otimes T$}\ |\ \hbox{$v\in B(r), T\in B(\gamma)$
such that ${\rm wt}(v)+{\rm wt}(T)=\mu$} \},
\cr
(B(\gamma)\otimes B(r))_\mu &= \{
\hbox{pairs $T\otimes v$}\ |\ \hbox{$v\in B(r), T\in B(\gamma)$
such that ${\rm wt}(v)+{\rm wt}(T)=\mu$} \}. \cr
}
\formula$$
The following lemma gives tableau versions of the Pieri rule
[Mac, I (5.16)].  The second bijection of the lemma is proved
in [Fu, \S 1.1 Proposition] and the proof of the first bijection is 
similar (see also [Bt, Propositions 2.3.4 and 2.3.11]).

\lemma  {\sl Let $\gamma,\mu,\tau\in {\cal P}$ be partitions and 
let $r,s\in \ZZ_{\ge 0}$.  There are bijections
$$\matrix{
(B(r)\otimes B(\gamma))_\mu &\longleftrightarrow
& B(\gamma\otimes (r))_\mu \cr
v\otimes T &\longrightarrow &v*T \cr
}
\qquad\hbox{and}\qquad
\matrix{
(B(\gamma)\otimes B(s))_\tau &\longleftrightarrow
& B(\gamma\otimes (s))_\tau \cr
T\otimes u &\longrightarrow &T*u \cr
}$$
}
\endlemma

\vfill\eject

\subsection Charge

\medskip
Let 
$\displaystyle{
B({\cal P})_{\ge} = \bigcup_{1\le i\le n} B({\cal P})_{\ge i} }$, where
$$B({\cal P})_{\ge i}
=\left\{ \hbox{column strict tableaux $b$}\ \Big\vert\ 
\matrix{
\hbox{${\rm wt}(b) = (\mu_1,\ldots,\mu_n)$ has} \hfill \cr
\hbox{$\mu_1=\cdots=\mu_{i-1}=0$ and $\mu_i\ge\cdots\ge\mu_n\ge 0$}
\hfill \cr}
\right\}.$$
Let 
$i^k=
\beginpicture
\setcoordinatesystem units <0.25cm,0.25cm>         
\setplotarea x from 0 to 6, y from -1 to 1    
\linethickness=0.5pt                          
\put{$\scriptstyle{i}$} at 0.5 0.5
\put{$\scriptstyle{i}$} at 1.5 0.5
\put{$\cdots$} at 3 0.5
\put{$\scriptstyle{i}$} at 4.5 0.5
\putrule from 0 0 to 5 0          %
\putrule from 0 1 to 5 1            
\putrule from 0 -0 to 0 1        %
\putrule from 1 -0 to 1 1        %
\putrule from 2 -0 to 2 1        %
\putrule from 4 -0 to 4 1        
\putrule from 5 -0 to 5 1        %
\endpicture
$
be the unique column strict tableau of shape
$(k)$ and weight $(0,\ldots,k,0,\ldots,0)$, where the $k$
appears in the $i\,$th entry.
{\it Charge} is the function
${\rm ch}\colon B({\cal P})_{\ge} \longrightarrow \ZZ_{\ge 0}$
such that
\smallskip\noindent
\itemitem{(a)}  ${\rm ch}(\emptyset) = 0$,
\smallskip\noindent
\itemitem{(b)}  if $T\in B({\cal P})_{\ge (i+1)}$ and 
$T*i^{\mu_i}\in B({\cal P})_{\ge i}$
then ${\rm ch}(T*i^{\mu_i}) = {\rm ch}(T)$,
\smallskip\noindent
\itemitem{(c)} if $T\in B({\cal P})_{\ge i}$ and 
$x$ is a letter not equal to $i$ then
${\rm ch}(x*T)={\rm ch}(T*x)+1$.
\medskip\noindent
The proof of the existence and uniqueness of the function ${\rm ch}$
is presented beautifully in [Ki].

\bigskip
\thm (Lascoux-Sch\"utzenberger [LS], [Sch])
{\sl For partitions $\lambda$ and $\mu$,
$$K_{\lambda\mu}(t) = \sum_{b\in B(\lambda)_\mu}
t^{{\rm ch}(b)}\,,$$
where the sum is over all column strict tableaux $b$ of shape $\lambda$
and weight $\mu$.
}
\pf
The proof is by induction on $n$.
Assume that the statement of the theorem holds for all partitions
$\mu = (\mu_1,\ldots, \mu_n)$.
We shall prove that, for all partitions 
$(\mu_0,\mu)=(\mu_0,\mu_1,\ldots,\mu_n)$,
$Q_{(\mu_0,\mu)}$ has an expansion 
$$Q_{(\mu_0,\mu)} = \sum_{p\in B(\nu)_{(\mu_0,\mu)}} t^{{\rm ch}(p)}s_\nu,
\formula$$

Beginning with the expression (4.5),
$$
Q_{(\mu_0,\mu)} 
= \left(\prod_{0\le i<j\le n} {1\over 1-tR_{ij}}\right)
s_{(\mu_0,\mu)} 
= \left(\prod_{j=1}^n {1\over 1-tR_{0j}}\right)
\left(\prod_{1\le i<j\le n} {1\over 1-tR_{ij}}\right) s_{(\mu_0,\mu)}. 
$$
By the definition of the Kostka-Foulkes polynomials (4.6) this is equal to
$$\eqalign{
Q_{(\mu_0,\mu)} 
&= \left(\prod_{j=1}^n {1\over 1-tR_{0j}}\right) 
\sum_{\lambda\in {\cal P}} K_{\lambda\mu}(t) 
s_{(\mu_0,\lambda)} \cr
&=\sum_{\lambda\in {\cal P}} 
K_{\lambda\mu}(t) 
\sum_{r\in \ZZ_{\ge 0}} t^r 
\sum_{k_1,\ldots, k_n\in \ZZ_{\ge 0} \atop
k_1+\cdots+k_n=r} R_{01}^{k_1}\cdots R_{0n}^{k_n}
s_{(\mu_0,\lambda)} \cr
&=\sum_{\lambda\in {\cal P}} 
K_{\lambda\mu}(t) 
\sum_{r\in \ZZ_{\ge 0}} t^r 
\sum_{k_1,\ldots, k_n\in \ZZ_{\ge 0} \atop
k_1+\cdots+k_n=r} 
s_{(\mu_0+r,\lambda-(k_1,\ldots,k_n))} \cr
}$$
Let $\gamma = \lambda-(k_1,\ldots, k_n)$ be such that 
$\lambda/\gamma$ is not a horizontal strip (usually
$\gamma$ isn't even a partition).  Let $m$ be the
first place a violation to being a horizontal strip occurs, i.e.
$$\hbox{let $m$ be minimal such that $\lambda_m-k_m< \lambda_{m+1}$.}$$
For example, in the followng picture, $\gamma = \lambda-(3,1,2,2,1,0)$
and $m=3$.
$$
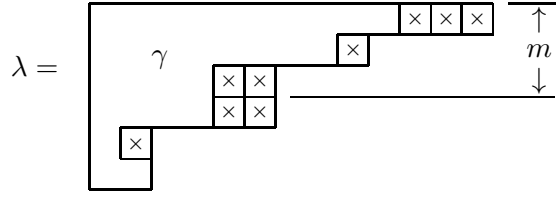
\beginpicture
\setcoordinatesystem units <.75cm,0.375cm>         
\setplotarea x from -2 to 6.5, y from 0 to 7    
\linethickness=0.5pt                          
\put{$\lambda = $}[r] at -0.5 4
\putrule from 0 6 to 6.5 6          %
\putrule from 4 5 to 6.5 5          
\putrule from 2 4 to 4.5 4          %
\putrule from 2 3 to 3 3          %
\putrule from 0.5 2 to 3 2          %
\putrule from 0.5 1 to 1 1        %
\putrule from 0 0 to 1 0          %
\putrule from 0 0 to 0 6        %
\putrule from 0.5 1 to 0.5 2    %
\putrule from 1 0 to 1 2        %
\putrule from 2 2 to 2 4        %
\putrule from 2.5 2 to 2.5 4    %
\putrule from 3 2 to 3 4        
\putrule from 4 4 to 4 5        %
\putrule from 4.5 4 to 4.5 5        %
\putrule from 5 5 to 5 6        %
\putrule from 5.5 5 to 5.5 6    %
\putrule from 6 5 to 6 6        %
\putrule from 6.5 5 to 6.5 6    %
\put{$\gamma$}[tl] at 1 4.5
\put{$\scriptstyle{\times}$} at 0.75 1.5
\put{$\scriptstyle{\times}$} at 2.25 2.5
\put{$\scriptstyle{\times}$} at 2.75 2.5
\put{$\scriptstyle{\times}$} at 2.25 3.5
\put{$\scriptstyle{\times}$} at 2.75 3.5
\put{$\scriptstyle{\times}$} at 4.25 4.5
\put{$\scriptstyle{\times}$} at 5.25 5.5
\put{$\scriptstyle{\times}$} at 5.75 5.5
\put{$\scriptstyle{\times}$} at 6.25 5.5
\putrule from 3.25 3 to 7.5 3
\putrule from 6.75 6 to 7.5 6
\put{$\uparrow$}[t] at 7.25 5.9
\put{$\downarrow$}[b] at 7.25 3.1
\put{$m$} at 7.25 4.5
\endpicture
$$
Let $s_m$ be the simple transposition in the symmetric group
which switches $m$ and $m+1$ and define
$$\tilde \gamma = s_m\circ \gamma,
\qquad\hbox{so that}\qquad 
s_{(\mu_0+r,\gamma)}=-s_{(\mu_0+r,\tilde \gamma)}.$$
Then $\tilde\gamma = \lambda-(\tilde k_1,\ldots, \tilde k_n)$ with 
$
\lambda_i-\tilde k_i = \lambda_i-k_i$, for $i\ne m,m+1$,
and
$$\lambda_m-\tilde k_m = \lambda_{m+1}-k_{m+1}-1, 
\qquad\hbox{and}\qquad
\lambda_{m+1}-\tilde k_{m+1} = \lambda_m-k_m+1. 
$$
Thus $\tilde \gamma_m = \lambda_{m+1}- k_{m+1}-1 <\lambda_{m+1}$
and so $\lambda/\tilde \gamma$ is not a horizontal strip.
This pairing $\gamma\leftrightarrow \tilde \gamma$ provides a cancellation
in the expression for $Q_{(\mu_0,\mu)}$ and thus
$$
Q_{(\mu_0,\mu)}
=\sum_{\lambda\in {\cal P}} 
\sum_{r\in \ZZ_{\ge 0}} t^r K_{\lambda\mu}(t) 
\sum_{\gamma\in {\cal P}
\atop \lambda\in \gamma\otimes (r)} s_{(\mu_0+r,\gamma)} 
=
\sum_{\gamma,r} 
\sum_{\lambda\in {\cal P} \atop \lambda\in \gamma\otimes (r)} 
t^r K_{\lambda\mu}(t) s_{(\mu_0+r,\gamma)}\,,
$$
where $\gamma\otimes (r)$ is as defined in (4.10).
By the induction assumption this is equal to
$$
Q_{(\mu_0,\mu)} =
\sum_{\gamma, r} 
\sum_{\lambda\in {\cal P} \atop \lambda\in \gamma\otimes (r)} 
\sum_{b\in B(\lambda)_\mu}
t^r t^{{\rm ch}(b)} s_{(\mu_0+r,\gamma)} 
=
\sum_{\gamma,r} 
\sum_{b\in B(\gamma\otimes(r))_\mu}
t^{r+{\rm ch}(b)} s_{(\mu_0+r,\gamma)}, 
$$
with $B(\gamma\otimes (r))_\mu$ as in (4.7).
By the first bijection in Lemma 4.11 this can be rewritten as
$$
\eqalign{
Q_{(\mu_0,\mu)}
&=
\sum_{\gamma, r} 
\sum_{v\otimes T\in (B(r)\otimes B(\gamma))_\mu}
t^{r+{\rm ch}(v*T)} s_{(\mu_0+r,\gamma)} \cr
&=
\sum_{\gamma, r}
\sum_{v\otimes T\in (B(r)\otimes B(\gamma))_\mu}
t^{r+{\rm ch}(v*T*0^{\mu_0}))} s_{(\mu_0+r,\gamma)} \cr
&=
\sum_{\gamma, r}
\sum_{v\otimes T\in (B(r)\otimes B(\gamma))_\mu}
t^{{\rm ch}(T*0^{\mu_0}*v))} s_{(\mu_0+r,\gamma)}, \cr
}
\formula$$
where the last two equalities come from the defining properties
of the charge function ${\rm ch}$.

Let $v\otimes T\in (B(r)\otimes B(\gamma))_\mu$ and let
$$p = T*0^{\mu_0}*v 
\qquad \hbox{and}\qquad 
\nu = {\rm shp}(T*0^{\mu_0}*v).$$
Let $d$ be such that 
$$\mu_0+r+d > \nu_d
\qquad\hbox{and}\qquad
\mu_0+r+d-1\le \nu_{d-1},$$
where, by convention, $\nu_0 = \mu_0+r$.
If $d>1$ define $\tilde \gamma$ and $\tilde r$ by
$$\tilde \gamma = 
(\gamma_1,\ldots, \gamma_{d-2},\mu_0+r+d-1,\gamma_d,\ldots, \gamma_n)
\quad\hbox{and}\quad
\mu_0+\tilde r+d-1 = \gamma_{d-1},$$
so that, if
$s_i$ denotes the transposition $(i,i+1)$ in the symmetric group, then
$(\mu_0+\tilde r,\tilde \gamma)
=(s_0\cdots s_{d-3}s_{d-2}s_{d-3}\cdots s_0)\circ (\mu_0+r,\gamma)$,
and
$$
s_{(\mu_0+r,\gamma)}
=(-1)^{2(d-3)+1}s_{(\mu_0+\tilde r,\tilde\gamma)}
=-s_{(\mu_0+\tilde r,\tilde\gamma)}.
\formula
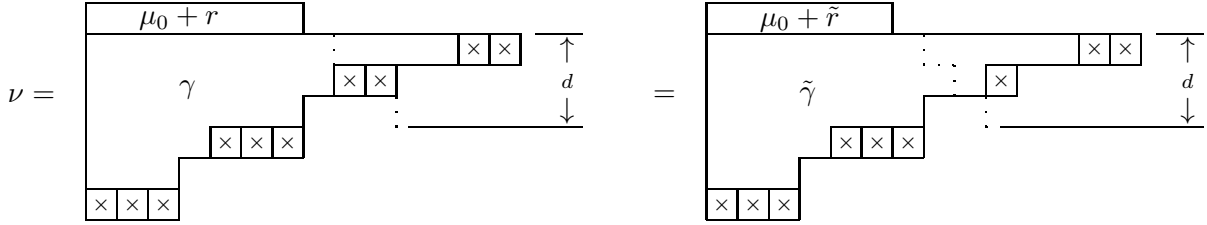
$$
$$
\beginpicture
\setcoordinatesystem units <.75cm,0.375cm>         
\setplotarea x from -2 to 8, y from 0 to 7    
\linethickness=0.5pt                          
\put{$\nu=$}[r] at -1 4
\putrule from -0.5 6 to 6.5 6       %
\putrule from 4 5 to 6.5 5          
\putrule from 3 4 to 4.5 4          %
\putrule from 1.5 3 to 3 3            %
\putrule from 1 2 to 3 2            %
\putrule from -0.5 1 to 1 1         %
\putrule from -0.5 0 to 1 0         %
\putrule from -0.5 0 to -0.5 7  %
\putrule from 0 0 to 0 1        %
\putrule from 0.5 0 to 0.5 1    %
\putrule from 1 0 to 1 2        %
\putrule from 1.5 2 to 1.5 3    %
\putrule from 2 2 to 2 3        %
\putrule from 2.5 2 to 2.5 3    %
\putrule from 3 2 to 3 4        
\putrule from 4 4 to 4 5        %
\putrule from 4.5 4 to 4.5 5    %
\putrule from 5.5 5 to 5.5 6    %
\putrule from 6 5 to 6 6        %
\putrule from 6.5 5 to 6.5 6    %
\putrule from -0.5 7 to 3 7
\putrule from 3 6 to 3 7
\put{$\mu_0+ r$} at 1 6.5
\put{$\gamma$}[tl] at 1 4.5
\put{$\scriptstyle{\times}$} at -0.25 0.5
\put{$\scriptstyle{\times}$} at 0.25 0.5
\put{$\scriptstyle{\times}$} at 0.75 0.5
\put{$\scriptstyle{\times}$} at 1.75 2.5
\put{$\scriptstyle{\times}$} at 2.25 2.5
\put{$\scriptstyle{\times}$} at 2.75 2.5
\put{$\scriptstyle{\times}$} at 4.25 4.5
\put{$\scriptstyle{\times}$} at 5.75 5.5
\put{$\scriptstyle{\times}$} at 6.25 5.5
\putrule from 4.7 3 to 7.5 3
\putrule from 6.75 6 to 7.5 6
\put{$\downarrow$}[b] at 7.25 3.1
\put{$\scriptstyle{d}$} at 7.25 4.5
\put{$\uparrow$}[t] at 7.25 5.9
\put{$\scriptstyle{\times}$} at 3.75 4.5
\putrule from 3.5 5 to 4 5
\putrule from 3.5 4 to 3.5 5
\setdots
\putrule from 3.5 5 to 3.5 6
\putrule from 3.5 5 to 4 5 
\putrule from 4 4 to 4 5
\putrule from 4 4 to 4.5 4
\putrule from 4.5 3 to 4.5 4
\putrule from 4.5 3 to 5 3
\endpicture
\beginpicture
\setcoordinatesystem units <.75cm,0.375cm>         
\setplotarea x from -2 to 6.5, y from 0 to 7    
\linethickness=0.5pt                          
\put{$=$}[r] at -1 4
\putrule from -0.5 6 to 6.5 6       %
\putrule from 4 5 to 6.5 5          
\putrule from 3 4 to 4.5 4          %
\putrule from 1.5 3 to 3 3          %
\putrule from 1 2 to 3 2            %
\putrule from -0.5 1 to 1 1         %
\putrule from -0.5 0 to 1 0         %
\putrule from -0.5 0 to -0.5 7  %
\putrule from  0 0 to 0 1       %
\putrule from 0.5 0 to 0.5 1    %
\putrule from 1 0 to 1 2        %
\putrule from 1.5 2 to 1.5 3    %
\putrule from 2 2 to 2 3        %
\putrule from 2.5 2 to 2.5 3    %
\putrule from 3 2 to 3 4        
\putrule from 4 4 to 4 5        %
\putrule from 4.5 4 to 4.5 5    %
\putrule from 5.5 5 to 5.5 6    %
\putrule from 6 5 to 6 6        %
\putrule from 6.5 5 to 6.5 6    %
\putrule from -0.5 7 to 2.5 7
\putrule from 2.5 6 to 2.5 7
\put{$\mu_0+\tilde r$} at 1 6.5
\put{$\tilde \gamma$}[tl] at 1 4.5
\put{$\scriptstyle{\times}$} at -0.25 0.5
\put{$\scriptstyle{\times}$} at 0.25 0.5
\put{$\scriptstyle{\times}$} at 0.75 0.5
\put{$\scriptstyle{\times}$} at 1.75 2.5
\put{$\scriptstyle{\times}$} at 2.25 2.5
\put{$\scriptstyle{\times}$} at 2.75 2.5
\put{$\scriptstyle{\times}$} at 4.25 4.5
\put{$\scriptstyle{\times}$} at 5.75 5.5
\put{$\scriptstyle{\times}$} at 6.25 5.5
\putrule from 4.25 3 to 7.5 3
\putrule from 6.75 6 to 7.5 6
\put{$\downarrow$}[b] at 7.25 3.1
\put{$\scriptstyle{d}$} at 7.25 4.5
\put{$\uparrow$}[t] at 7.25 5.9
\setdots
\putrule from 3 5 to 3 6
\putrule from 3 5 to 3.5 5 
\putrule from 3.5 4 to 3.5 5
\putrule from 3.5 4 to 4 4
\putrule from 4 3 to 4 4
\putrule from 4 3 to 4.5 3
\endpicture
$$
Note that $\tilde{\tilde\gamma}=\gamma$ and $\tilde{\tilde r}=r$.

\medskip\noindent
{\it Case 1:}  $d>1$ and 
$(\mu_0+r,\gamma) = (\mu_0+\tilde r, \tilde \gamma)$.
\quad
In this case (4.15) implies $s_{(\mu_0+r,\gamma)}=0$.

\medskip\noindent
{\it Case 2:}  $d>1$ and 
$(\mu_0+r,\gamma) \ne (\mu_0+\tilde r, \tilde \gamma)$.
\quad
Then
$$\nu\in \gamma\otimes (\mu_0+r)
\qquad\hbox{and}\qquad
\nu\in \tilde\gamma\otimes (\mu_0+\tilde r).$$
Row uninserting 
the horizontal strips $\nu/\gamma$ and $\nu/\tilde \gamma$
from $p$, i.e. using the second bijection in Lemma 4.11,
produces pairs
$$
T\otimes u=T\otimes (0^{\mu_0}*v)
\in (B(\gamma)\otimes B(\mu_0+r))_{(\mu_0,\mu)}
\qquad\hbox{and}\qquad
\tilde T\otimes \tilde u\in 
(B(\tilde\gamma)\otimes B(\mu_0+\tilde r))_{(\mu_0,\mu)},
$$
respectively.  Consider the $\ell=\mu_0+r$ 
bumping paths in the tableau $p$ which arise from $T*u$.
These begin with the letters $u_1\le \ldots \le u_\ell$ of $u$ and end 
at the boxes of the horizontal strip $\nu/\gamma$.
Similarly, there are $\tilde \ell = \mu_0+\tilde r$ bumping paths
in $p$ arising from $\tilde T*\tilde u$.  Note that
\smallskip\noindent
\item{(a)} since $u=0^{\mu_0}*v$  begins with $\mu_0$ ~$0\,$s 
the leftmost $\mu_0$ bumping paths in $T*u$ travel vertically,
directly down the first $\mu_0$ columns of $p$, and 
\item{(b)} in rows numbered $\ge d$ the bumping paths
for $\tilde T*\tilde u$ coincide exactly with the
bumping paths for $T*u$, since the horizontal strips
$\nu/\gamma$ and $\nu/\tilde \gamma$ coincide exactly
in rows $\ge d$ and these paths are obtained by uninserting
the boxes in this portion of the horizontal strip.
\smallskip\noindent
$$\matrix
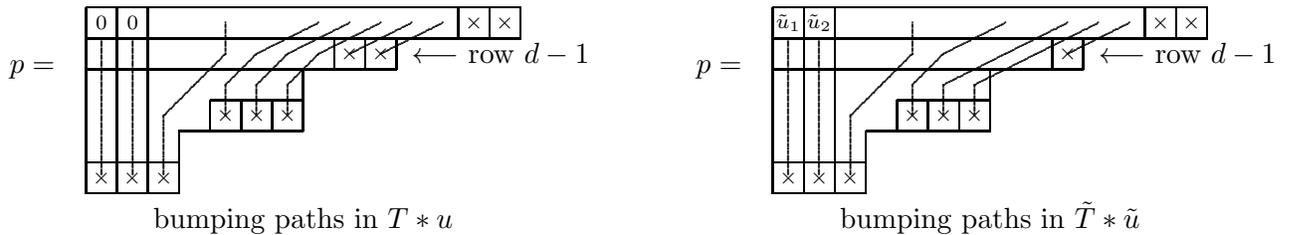
{
\beginpicture
\setcoordinatesystem units <.75cm,0.375cm>         
\setplotarea x from -2 to 8, y from 0 to 7    
\linethickness=0.5pt                          
\put{$p=$}[r] at -1 4
\put{$\scriptstyle{0}$} at -0.25 5.5
\put{$\scriptstyle{0}$} at 0.25 5.5
\putrule from -0.5 6 to 6.5 6       %
\putrule from -0.5 5 to 6.5 5       
\putrule from -0.5 4 to 4.5 4       %
\putrule from 1.5 3 to 3 3          %
\putrule from 1 2 to 3 2            %
\putrule from -0.5 1 to 1 1         %
\putrule from -0.5 0 to 1 0         %
\putrule from -0.5 0 to -0.5 6  %
\putrule from 0 0 to 0 6        %
\putrule from 0.5 0 to 0.5 6    %
\putrule from 1 0 to 1 2        %
\putrule from 1.5 2 to 1.5 3    %
\putrule from 2 2 to 2 3        %
\putrule from 2.5 2 to 2.5 3    %
\putrule from 3 2 to 3 4        
\putrule from 4 4 to 4 5        %
\putrule from 4.5 4 to 4.5 5    %
\putrule from 5.5 5 to 5.5 6    %
\putrule from 6 5 to 6 6        %
\putrule from 6.5 5 to 6.5 6    %
\put{$\scriptstyle{\times}$} at -0.25 0.5
\put{$\scriptstyle{\times}$} at 0.25 0.5
\put{$\scriptstyle{\times}$} at 0.75 0.5
\put{$\scriptstyle{\times}$} at 1.75 2.5
\put{$\scriptstyle{\times}$} at 2.25 2.5
\put{$\scriptstyle{\times}$} at 2.75 2.5
\put{$\scriptstyle{\times}$} at 4.25 4.5
\put{$\scriptstyle{\times}$} at 5.75 5.5
\put{$\scriptstyle{\times}$} at 6.25 5.5
\put{$\longleftarrow$ row $d-1$}[l] at 4.75 4.5
\put{$\scriptstyle{\times}$} at 3.75 4.5
\putrule from 3.5 5 to 4 5
\putrule from 3.5 4 to 3.5 5
\plot -0.25 0.6   -0.25 5    /
\plot 0.25 0.6   0.25 5      /
\plot 0.75 0.6   0.75 2.5    /
\plot 0.75 2.5   1.25 3.5      /
\plot 1.25 3.5   1.75 4.5      /
\plot 1.75 4.5   1.75 5.5      /
\plot 1.75 2.5   1.75 3.5      /
\plot 1.75 3.5   2.25 4.5      /
\plot 2.25 4.5   3.25 5.5      /
\plot 2.25 2.5   2.25 3.5      /
\plot 2.25 3.5   2.75 4.5      /
\plot 2.75 4.5   3.75 5.5      /
\plot 2.75 2.5   2.75 3.5      /
\plot 2.75 3.5   3.25 4.5      /
\plot 3.25 4.5   4.25 5.5      /
\plot 3.75 4.5   4.75 5.5      /
\plot 4.25 4.5   5.25 5.5      /
\endpicture
&\phantom{.}
&
\beginpicture
\setcoordinatesystem units <.75cm,0.375cm>         
\setplotarea x from -2 to 8, y from 0 to 7    
\linethickness=0.5pt                          
\put{$p=$}[r] at -1 4
\put{$\scriptstyle{\tilde u_1}$} at -0.25 5.5
\put{$\scriptstyle{\tilde u_2}$} at 0.25 5.5
\putrule from -0.5 6 to 6.5 6       %
\putrule from -0.5 5 to 6.5 5       
\putrule from -0.5 4 to 4.5 4       %
\putrule from 1.5 3 to 3 3          %
\putrule from 1 2 to 3 2            %
\putrule from -0.5 1 to 1 1         %
\putrule from -0.5 0 to 1 0         %
\putrule from -0.5 0 to -0.5 6  %
\putrule from 0 0 to 0 6        %
\putrule from 0.5 0 to 0.5 6    %
\putrule from 1 0 to 1 2        %
\putrule from 1.5 2 to 1.5 3    %
\putrule from 2 2 to 2 3        %
\putrule from 2.5 2 to 2.5 3    %
\putrule from 3 2 to 3 4        
\putrule from 4 4 to 4 5        %
\putrule from 4.5 4 to 4.5 5    %
\putrule from 5.5 5 to 5.5 6    %
\putrule from 6 5 to 6 6        %
\putrule from 6.5 5 to 6.5 6    %
\put{$\scriptstyle{\times}$} at -0.25 0.5
\put{$\scriptstyle{\times}$} at 0.25 0.5
\put{$\scriptstyle{\times}$} at 0.75 0.5
\put{$\scriptstyle{\times}$} at 1.75 2.5
\put{$\scriptstyle{\times}$} at 2.25 2.5
\put{$\scriptstyle{\times}$} at 2.75 2.5
\put{$\scriptstyle{\times}$} at 4.25 4.5
\put{$\scriptstyle{\times}$} at 5.75 5.5
\put{$\scriptstyle{\times}$} at 6.25 5.5
\put{$\longleftarrow$ row $d-1$}[l] at 4.75 4.5
\putrule from 3.5 5 to 4 5
\plot -0.25 0.6   -0.25 5    /
\plot 0.25 0.6   0.25 5      /
\plot 0.75 0.6   0.75 2.5    /
\plot 0.75 2.5   1.25 3.5      /
\plot 1.25 3.5   1.75 4.5      /
\plot 1.75 4.5   1.75 5.5      /
\plot 1.75 2.5   1.75 3.5      /
\plot 1.75 3.5   2.25 4.5      /
\plot 2.25 4.5   3.25 5.5      /
\plot 2.25 2.5   2.25 3.5      /
\plot 2.25 3.5   3.25 4.5      /
\plot 3.25 4.5   4.25 5.5      /
\plot 2.75 2.5   2.75 3.5      /
\plot 2.75 3.5   3.75 4.5      /
\plot 3.75 4.5   4.75 5.5      /
\plot 4.25 4.5   5.25 5.5      /
\endpicture
\cr
\hbox{bumping paths in $T*u$} 
&&\hbox{bumping paths in $\tilde T*\tilde u$} \cr
}
$$
Suppose there are $k$ bumping paths which
end in rows $\ge d$.  The picture above has $k=6$ and corresponds to
Case 2b below. 
\smallskip\noindent
\itemitem{{\it Case 2a:}}
If $\mu_0+\tilde r>\mu_0+r$ then
the $k$ bumping paths which end in rows $\ge d$ 
are the same or slightly ``more left'' 
in $\tilde T*\tilde u$ than in $T*u$.  Since the first
$\mu_0$ bumping paths cannot be any ``more left'' than
vertical, this forces that the 
first $\mu_0$ entries of $\tilde u$ are $0$, i.e. that 
$\tilde u = 0^{\mu_0}*\tilde v$ for some $v\in B(\tilde r)$.
\smallskip\noindent
\itemitem{{\it Case 2b:}}
If $\mu_0+\tilde r<\mu_0+r$ then
the $k$ bumping paths which end in rows $\ge d$ 
are the same or slightly ``more right'' 
in $\tilde T*\tilde u$ than in $T*u$.  
There are $k+r-\tilde r$ bumping paths of $T*u$
passing through the first $\mu_0+r-(d-1)$ squares of row $d-1$,
namely, the $k$ bumping paths of $T*u$ which end in rows $\ge d$ 
and the $(\mu_0+r)-(\mu_0+\tilde r)$
bumping paths of $T*u$ which end in row $d-1$ and which
do not appear as bumping paths for $\tilde T*\tilde u$.
The first $\mu_0$ of these paths pass through the squares
in positions $(d-1,1),\ldots, (d-1,\mu_0)$
and the last $r-\tilde r$ of them
pass through the squares in positions $(d-1,\mu_0+\tilde r+d-1+1),
\ldots (d-1, \mu_0+r+d-1)$.  Since the remaining number of paths,
$$k+r-\tilde r-\mu_0-(\mu_0+r-\mu_0-\tilde r)
=k-\mu_0 < \mu_0+\tilde r-\mu_0 < \mu_0+\tilde r+(d-1)-\mu_0,$$
there must be a box in position $(d-1,j)$
for some $\mu_0<j<\mu_0+\tilde r+(d-1)$ which does not
have a bumping path for $T*u$ passing through it.
All the bumping paths of $T*u$ which pass through row $d-1$
to the left of this
box remain the same as bumping paths for $\tilde T*\tilde u$ and
the first $\mu_0$ of these begin at an entry $0$ in 
the first row of $p$.  Thus, as in Case 2a, the first
$\mu_0$ entries of $\tilde u$ are $0$, i.e.
$\tilde u = 0^{\mu_0}*\tilde v$ for some $v\in B(\tilde r)$.

\smallskip\noindent
So,
$$\tilde T\otimes \tilde u = \tilde T\otimes (0^{\mu_0}*\tilde v),
\qquad\hbox{with}\quad 
\tilde v\otimes \tilde T\in (B(\tilde r)\otimes B(\tilde \gamma))_\mu,
$$
and the terms in the last line of (4.14) 
corresponding to the pairs $v\otimes T$ and
$\tilde v\otimes \tilde T$ cancel each other because
$$T*0^{\mu_0}*v=\tilde T*0^{\mu_0}*\tilde v
\qquad\hbox{and}\qquad
s_{(\mu_0+r,\gamma)}=-s_{(\mu_0+\tilde r,\tilde \gamma)}.$$

\smallskip\noindent
{\it Case 3:}  $d=1$.
\quad
Since $\mu_0+r+1>\nu_1$ and $\nu\in \gamma\otimes (\mu_0+r)$
the horizontal strip $\nu/\gamma$ has its boxes in each
of the first $\mu_0+r$ columns, i.e. 
$$\nu = (\nu_0,\nu_1,\ldots, \nu_n) = (\mu_0+r,\gamma_1,\ldots, \gamma_n)
=(\mu_0+r,\gamma).$$
Row uninsertion of the horizontal strip $\nu/\gamma$ from the column 
strict tableau $p$, i.e. using the second bijection in Lemma 4.11,
recovers the pair $T\otimes (0^{\mu_0}*v)$ and shows that $0^{\mu_0}*v$ is the 
first row of $p$.

\medskip
In conclusion, in the last line of (4.14) the terms corresponding 
to Case 1 vanish,
the terms corresponding to Case 2 cancel off and the remaining
Case 3 terms give formula (4.13), as desired.
\endpf

\section 5. References

\medskip\noindent
\item{[Bou]} {\smallcaps N.\ Bourbaki},
{\sl Groupes et algebres de Lie}, 
Chapt.\ IV-VI, Masson, Paris, 1981.

\medskip\noindent
\item{[Bt]} {\smallcaps L.M.\ Butler},
{\it Subgroup lattices and symmetric functions},
Memoirs Amer.\ Math.\ Soc.\ No.\ 539, {\bf 112}, 1994.

\medskip\noindent
\item{[Fu]} {\smallcaps W.\ Fulton},
{\sl Young tableaux: with applications to representation
theory and geometry},
London Math.\ Soc.\ Student Texts {\bf 35},
Cambridge Univ.\ Press, 1997.

\medskip\noindent
\item{[JLZ]} {\smallcaps A.\ Joseph, G.\ Letzter and S.\ Zelikson},
{\it On the Brylinski-Kostant filtration}
J.\ Amer.\ Math.\ Soc.\ {\bf 13} no.\ 4 (2000), 945--970.

\medskip\noindent
\item{[Kt]} {\smallcaps S.\ Kato},
{\it Spherical functions and a q-analogue of Kostant's weight multiplicity
formula}, Invent. Math. {\bf 66} (1982), 461--468.

\medskip\noindent
\item{[KL]} {\smallcaps D.\ Kazhdan and G.\ Lusztig},
{\it Representations of Coxeter groups and Hecke algebras},
Invent.\ Math.\ {\bf 53} (1979), 165--184.

\medskip\noindent
\item{[Ki]} {\smallcaps K.\ Killpatrick},
{\it A combinatorial proof of a recursion for the $q$-Kostka
polynomials}, J.\ Comb. Theory Ser.\ A {\bf 92} (2000), 29--53.

\medskip\noindent
\item{[Ks]} {\smallcaps B.\ Kostant},
{\it Lie group representations on polynomial rings},
Amer.\ J.\ of Math.\ {\bf 85} (1963), 327--404.

\medskip\noindent 
\item{[LS]} {\smallcaps A.\ Lascoux and M.P.\ Sch\"utzenberger},
{\it Sur une conjecture de H.O.\ Foulkes}, Compt. Rend. Acad. Sci. Paris
{\bf 286A} (1978), 323--324.

\medskip\noindent
\item{[Lu]} {\smallcaps G.\ Lusztig},
{\it Singularities, character formulas and a $q$-analog of 
weight multiplicities},
in {\sl Analysis and topology on singular spaces, II, III (Luminy, 1981)},
Ast\'erisque, {\bf 101-102}, (1983), 208--229.


\medskip\noindent
\item{[Mac]} {\smallcaps I.G.\ Macdonald},
{\sl Symmetric functions and Hall polynomials},
Second edition, Oxford Mathematical Monographs, Oxford University Press, 
New York, 1995.

\medskip\noindent
\item{[Mac2]} {\smallcaps I.G.\ Macdonald},
{\sl Spherical functions on a group of $p$-adic type},
Publ. Ramanujan Institute No.\ 2, Madras, 1971.

\medskip\noindent
\item{[Mac3]} {\smallcaps I.G.\ Macdonald},
{\it The Poincar\'e series of a Coxeter group},
Math.\ Ann.\ {\bf 199} (1972), 161--174.

\medskip\noindent
\item{[Mac4]} {\smallcaps I.G.\ Macdonald},
{\it Affine Hecke algebras and orthogonal polynomials},
S\'eminaire Bourbaki, 47\`eme ann\'ee, ${\rm n}^{\rm o}$ 797, 1994--95,
Ast\'erisque {\bf 237} (1996), 189--207.

\medskip\noindent
\item{[Nak1]} {\smallcaps H.\ Nakajima},
{\sl Lectures on Hilbert schemes of points on surfaces},
Univ.\  Lect.\ Ser. {\bf 18}, Amer.\ Math.\ Soc., Providence, RI, 1999.

\medskip\noindent
\item{[Nak2]} {\smallcaps H.\ Nakajima},
{\it Quiver varieties and finite dimensional representations of 
quantum affine algebras}, J.\ Amer.\ Math.\ Soc.\ {\bf 14} no.\ 1
(2001), 145--238.

\medskip\noindent
\item{[Sch]} {\smallcaps M.P.\ Sch\"utzenberger},
{\it Propri\'et\'es nouvelles des tableaux de Young},
S\'eminaire Delange-Pisot-Poitou, $19^{\rm e}$ ann\'ee 1977-1978, no.\ 26,
Secr\'eteriat Math\'ematique, Paris.

\medskip\noindent
\item{[Shj1]} {\smallcaps T.\ Shoji}, 
{\it Green functions of reductive groups over a finite field},
Proc.\ Symp.\ Pure Math.\ {\bf 47} (1987), 289--301.

\medskip\noindent
\item{[Shj2]} {\smallcaps T.\ Shoji}, 
{\it Green functions associated to complex reflection groups I and II},
J.\ Algebra {\bf 245} no. 2 (2001), 650--694; and to appear in J.\ Algebra.

\medskip\noindent
\item{[St]} {\smallcaps R.\ Steinberg},
{\sl Lectures on Chevalley groups}, 
Yale University,  New Haven, CT, 1968.

\vfill\eject
\end